\documentclass{amsart}
  \usepackage{amscd,amssymb,epsfig}

                     \numberwithin{equation}{subsection}

                     \newtheorem{propo}{Proposition}[subsection]
                     \newtheorem{corol}[propo]{Corollary}
                     \newtheorem{theor}[propo]{Theorem}
                     \newtheorem{lemma}[propo]{Lemma}
                     \theoremstyle{definition}

                     \theoremstyle{remark}

                     \newcommand{\ZZ}{\mathbb{Z}}
                     \newcommand{\RR}{\mathbb{R}}

\newcommand{\A}{\mathcal{A}}
 \newcommand{\B}{\mathcal{B}}

\newcommand{\Ker}{\operatorname{Ker}}

\newcommand{\sign}{\operatorname{sign}}

\newcommand{\rk}{\operatorname{rk}}

\newcommand{\R}{\mathcal{R}}

\newcommand{\C}{\mathcal{C}}

\newcommand{\N}{\mathcal{N}}

\newcommand{\inc}{\operatorname{in}}

\newcommand{\card}{\operatorname{card}}

\newcommand{\ar}{\operatorname{g}}

                     \newcommand{\id}{\operatorname{id}}

\newcommand{\modu}{\operatorname{mod}}
\newcommand{\Int}{\operatorname{Int}}

\begin{document}
      \title{Cobordisms of   words}
                     \author[Vladimir Turaev]{Vladimir Turaev}
                     \address{%
              Department of Mathematics \newline
\indent  Indiana University \newline
                     \indent Bloomington IN47405 \newline
                     \indent USA \newline
\indent e-mail: vtouraev@indiana.edu }
                     \begin{abstract}  We introduce an equivalence  relation, called cobordism,  for words and produce cobordism invariants of words.
                     \end{abstract}
                     \maketitle

   \section {Introduction}

 Finite sequences of elements  of a given set $\alpha$ are called
words in the alphabet~$\alpha$.  Words have been extensively
studied by algebraic and
  combinatorial means, see \cite{lo1}, \cite{lo2}.  Gauss \cite{ga}
        used     words to encode closed  plane curves, viewed up to homeomorphism.
         For further
        work on Gauss words of curves,
         see  \cite{ro1}, \cite{ro2}, \cite {lm},
           \cite{dt}, \cite{ce}, \cite{cr}.

         Words can be investigated  using
          ideas and techniques from low-dimensional topology.
          The relevance of topology is suggested by
   the connection to
  curves and also by the  phenomenon of
   linking of letters in   words. A prototypical example is
      provided by the words   $abab$ and $aabb$. The letters $a,b$ are obviously
      linked in the first word and unlinked in the second one. A similar
      linking phenomenon for  geometric objects, for instance  knotted circles in
       Euclidean 3-space,  is  studied in knot theory.

A   study of words, based on a transposition of topological ideas,
was started by the author in \cite{tu3}. We begin by fixing an
alphabet (a set) $\alpha$ with involution $\tau:\alpha\to \alpha$.
The concept of generic
  curves, which may  have only  double self-intersections,
   leads to a notion of   nanowords over    $\alpha$.
    Every letter appearing in a nanoword occurs in it exactly twice.
     Using an   analogy with homological intersection numbers of curves on a surface,
     we
      associate  with any  nanoword over $\alpha$ a certain pairing  called
       $\alpha$-pairing.
 The concept of deformation  of curves on a surface
  can be also transposed to the setting of words. One can
   view a deformation of a curve as a sequence of
  local transformations or moves
  following certain simple models. Similar   homotopy moves  can be defined
  for  nanowords over  $\alpha$; they generate
   an equivalence relation of homotopy.

In this paper we introduce further transformations on nanowords
called surgeries. In topology, surgery is an operation on
manifolds consisting in cutting out a certain submanifold (with
boundary) and gluing
 at its place another  manifold with the same boundary.  Various  operations
  of this kind can be considered  for words. We take here the following
    approach:  a surgery  on a  nanoword   deletes a symmetric subnanoword,
    i.e., a subnanoword isomorphic to its opposite.
    More generally, a surgery on a nanoword may delete
    a symmetric subnanophrase.
     Symmetry  plays here  the   role of  the Poincar\'e duality for
     manifolds;   symmetric nanophrases  are moral analogues  of manifolds.

 Surgeries and homotopy moves generate    an  equivalence relation
on the class of nanowords called cobordism. The main aim of the
theory of  cobordisms of  words is to
 classify nanowords up to cobordism or, equivalently,  to  compute the set of   cobordism classes of nanowords over $\alpha$. This set, denoted
  $\N_c =\N_c (\alpha, \tau)$,  is a group  with respect
  to concatenation of nanowords. We use  $\alpha$-pairings and other  homotopy  invariants of nanowords introduced in \cite{tu3} to construct group homomorphisms from $\N_c $ to simpler groups. We prove   that   $\mathcal N_c $ is infinitely generated   provided   $\tau\neq \id$ and  is   non-abelian provided  $\tau$ has at least 3 orbits.

 For   non-cobordant  nanowords,  it is interesting to  measure  how far they are from being cobordant. In other terms we are interested in finding natural  metrics on  $\N_c$. The group structure on $\N_c$ allows us to derive such
 metrics  from norms on  $\N_c$.
We define  two $\ZZ$-valued norms on   $\N_c$:
 the length norm and
the bridge norm. The length norm counts the minimal length of a
nanoword in the given cobordism class. This corresponds to the
topological notion of the  minimal  number of crossings
 of a curve. The bridge norm reflects the idea of a surface of  minimal
  genus  spanned in a 3-manifold $M$ by a loop  in $\partial M$. To define the bridge norm we introduce
so-called bridge moves on nanowords generalizing  surgery.  The
metric on  $\N_c$ induced by the bridge norm has a nice feature of
being invariant under both  left and right translations.   We give
an estimate from below  for this metric involving   a numerical
invariant  of $\alpha$-pairings,  the so-called genus.

  Adding cyclic permutations of nanowords to
the list of moves,  we obtain a notion of weak cobordism and also
define a weak bridge pseudo-metric on   $\N_c$. We estimate this
pseudo-metric from below using the genus of $\alpha$-pairings.

A canonical  procedure, called  desingularization,  transforms any
word $w$ in the alphabet $\alpha$ into a nanoword over
 $\alpha$, see \cite{tu3}.
 The latter  determines an
 element of  $\N_c$, called the cobordism class of $w$.
   This allows one to apply the invariants, metrics, etc. introduced
   in this paper
    to usual words in the  alphabet $\alpha$.

The main results of this paper are  a  construction of a
homomorphism from
 $\N_c$ to a group of cobordism classes of  $\alpha$-pairings
 (Theorem \ref{lemplo}) and the estimates of the bridge
 metric and the weak bridge pseudo-metric
 via the genus  (Corollary \ref{th:dfrovv53}, Theorem \ref{thcircularv53}).

The organization of the paper is as follows. The definitions of
nanowords and nanophrases are given in Sect.\ 2. In Sect.\ 3 and 4
we introduce cobordisms of nanowords and discuss a simple
cobordism invariant. Sect.\ 5 and 6  are concerned with the
general theory of $\alpha$-pairings. In Sect.\ 7 and 8 we discuss
the $\alpha$-pairings of nanowords and the associated cobordism
invariants of  nanowords.  Sect.\   9 -- 11 are devoted to  the
bridge metric, the circular shift on nanowords,   and the weak
bridge pseudo-metric. In Sect.\ 12 -- 14 we discuss connections between
words and bridge moves on the one hand and loops  on surfaces and surfaces in 3-manifolds on the other hand. We use these connections to  prove two  lemmas
from Sect.\  7 and 9.   Note that
 Sect.\  1 -- 11  are  written in a purely algebraic language
 while Sect.\ 12 -- 14 use  elementary   topology.

This work  is a sequel to  \cite{tu1} --  \cite{tu5} but
 a knowledge of these papers is not required.

 Throughout the paper, the symbol $\alpha$ denotes a  fixed set endowed with an  involution $\tau:\alpha\to \alpha$.

    \section{Nanowords and nanophrases}\label{aza}

  In this section we recall   the basics of the theory of nanowords, see  \cite{tu3}.

\subsection{Words and nanowords}\label{word} For a positive integer $n$, set $\widehat  n= \{1,2,...,n\}$. A {\it word  of length} $n$ in the alphabet  $\alpha$ is  a  mapping $w:\widehat  n\to   \alpha$.  Such a
word
$w$  is encoded by the  sequence $w(1)\,  w(2) \cdots  w(n)$.
  Writing the letters of   $w$ in the opposite order  we obtain the {\it opposite word} $w^-= w(n)\, w(n-1) \cdots w(1)$  in the same alphabet.

  An {\it   $\alpha$-alphabet} is a  set $\A$ endowed with  a  mapping  $ \A\to \alpha$
called {\it projection}. The image  of   $A\in
\A$   under this mapping  is denoted $\vert A\vert$. An {\it isomorphism}  of   $\alpha$-alphabets $\A_1$,  $\A_2$ is a bijection  $f:\A_1\to
\A_2$ such that    $\vert A\vert=\vert f(A)\vert$ for all $A\in \A_1$.

A   {\it nanoword  of length} $n$ over $\alpha$ is  a pair (an $\alpha$-alphabet $\A$,   a mapping $w:\widehat  n\to   \A$ such that each element of $\A$ is the image of precisely two elements of  $\widehat  n$).       Clearly,    $n=2\card (\A)$.   By definition, there is a   unique  {\it empty nanoword} $\emptyset$  of length 0.

 We say that  nanowords    ($\A$,
$w$) and  ($\A'$,
$w'$)  over $\alpha$  are {\it   isomorphic} and write $w\approx w'$  if    there is an  isomorphism   of   $\alpha$-alphabets    $f:\A\to \A'$
  such that  $w'  =f  w$.

The {\it concatenation  product}    of two  nanowords  ($\A_1$,
$w_1$) and  ($\A_2$, $w_2$)  is defined as follows. Replacing if
necessary  ($\A_1$, $w_1$) with an isomorphic nanoword we can
assume that $\A_1\cap \A_2=\emptyset$.  Then the product  of $w_1$
and $w_2$  is the  nanoword  $(\A_1 \cup \A_2, w_1 w_2)$ where
$w_1w_2$ is obtained from $w_1, w_2$  by concatenation.  The
nanoword  $(\A_1 \cup \A_2, w_1 w_2)$ is well defined up to
isomorphism. Multiplication of  nanowords    is associative and
has a unit $\emptyset$ (the  empty nanoword).

          A nanoword $w:\widehat  n \to \A$   is {\it symmetric}
           if it is isomorphic to the opposite nanoword
            $ w^-:\widehat  n \to \A$, i.e., if there is a bijection
             $\iota: \A \to \A$ commuting with the projection
             to $\alpha$ and such that $\iota w=w^-$.   The latter means that
                $\iota ( w(i))= w(n+1-i)$ for all $i \in \widehat  n$.
               Clearly,  $\iota$ is uniquely
                determined by $w$ and   $\iota^2=\id$.
                     For example, the nanoword $ABBA$ with arbitrary $\vert A\vert, \vert B\vert \in \alpha$  is
                       symmetric with $\iota=\id$. The nanoword $ABAB$  is symmetric  if and only if
          $\vert A\vert=\vert B\vert$.

\subsection{Homotopy}\label{unkn21}  There are  three basic  transformations of nanowords called {\it homotopy moves}.
 The first  of them  transforms a nanoword        $(\A,   x A
A  y)$ with $A\in \A $    into the nanoword $(\A-\{A\},  xy)$. The
second homotopy move  transforms a nanoword     $(\A ,   xA  B  y
BAz  )$   where $A,B\in  \A$ with $\vert B\vert =\tau ( \vert
A\vert)$    into $(\A-\{A,B\},   xyz)$.   The third move
transforms a nanoword     $(\A, xAByACzBCt)$   where  $A,B,C \in
\A$ are distinct letters with  $ \vert A\vert =\vert B\vert=\vert
C\vert  $    into  $ (\A, xBAyCAzCBt)$.

The  homotopy moves  and isomorphisms generate an equivalence
relation of {\it homotopy} in the class of nanowords. Nanowords
homotopic   to $\emptyset$ are {\it contractible}.     For
example, for $a, b  \in \alpha$ consider  the nanoword   $
w_{a,b}= ABAB$  with $ \vert  A\vert= a, \vert  B\vert=b$. A
homotopy classification of such nanowords is given in  \cite{tu3},
Theorem 8.4.1: $w_{a,b}$  is contractible if and only if  $a=\tau
(b)$, two non-contractible nanowords  $w_{a,b}, w_{a',b'}$ are
homotopic if and only if $a=a'$ and $b=b'$.

The third homotopy move has a  more general version  (see
\cite{tu4}),  but we shall not consider it here.

 \subsection{Nanophrases}\label{nanop}
 A sequence of words $w_1,..., w_k$ in an $\alpha$-alphabet $\A$
 is a {\it nanophrase of length $k$} (over $\alpha$)  if every letter of $\A$
 appears in $w_1,..., w_k$ exactly twice or, in other terms, if
  the   concatenation $w_1w_2\cdots w_k$ is a nanoword.
  We   denote such a nanophrase by $(\A, (w_1\,\vert \, \cdots \, \vert \, w_k))$ or
  shorter by $(w_1\,\vert \, \cdots \, \vert \, w_k)$. For a
   nanophrase $\nabla= (\A, (w_1\,\vert \, \cdots \, \vert \, w_k))$,
  define a function $\varepsilon_\nabla:\A\to \{0,1\}$
 by $\varepsilon_\nabla(A)=0$ if $A\in \A$ occurs twice in the
 same word of   $\nabla$ and $\varepsilon_\nabla(A)=1$ if $A\in \A$ occurs
  in different  words of   $\nabla$. Nanophrases of length 1 are just
   nanowords.

  A nanophrase $\nabla=(w_1\,\vert \, \cdots \, \vert \, w_k)$ is {\it symmetric}
   if there is a bijection  $\iota: \A \to \A$   such that $\iota w_r=w_r^-$
   for $r=1,...,k$ and
  $\vert \iota (A) \vert =\tau^{\varepsilon_\nabla (A)} (\vert A\vert)$ for all
   $A\in \A$. In the sequel we  often write $A^\iota$ for $\iota(A)$.
  The involution $\iota$
   transforms the $i$-th letter of $w_r$ into the $ (n_r+1-i)$-th letter of $w_r$
    for all $i,r$ where $n_r$ is the length of $w_r$.
    Hence $\iota=\iota_\nabla$ is determined by $\nabla$ uniquely and $\iota^2=\id$.
     For nanowords ($k=1$),  this notion of symmetry  coincides with the one in Sect.\ \ref{word}.

A nanophrase $(w_1\,\vert \, \cdots \, \vert \, w_k)$ is {\it
even} if all
 the words $w_1,..., w_k$ have even length. Note that
 the sum  of the  lengths of $w_1,..., w_k$  is always even. All
   nanophrases of length 1  are even.

For example, the nanophrase  $(AB \,\vert \, BA)$ is  even. It is
symmetric if and only if
  $\vert A \vert =\tau  (\vert B\vert)$. The   nanophrase $(A\,\vert \, A)$
  is not even. It is
symmetric if and only if
  $\vert A \vert =\tau  (\vert A\vert)$.

\subsection{Remark}\label{unremk1}  If the involution $\tau:\alpha \to \alpha$ is
  fixed-point-free, then any   symmetric
  nanophrase $\nabla$ over $\alpha$ is even. Indeed,
  if   $\nabla$ contains a word $w $ of odd length, then its central
  letter,  $A$,   satisfies $A^\iota=A$ for $\iota=\iota_\nabla$. If the second entry of $A$ in $\nabla$
  also
   occurs   in $w$, then $A$   occurs  in $w$ at
   least 3 times which contradicts the definition of a nanophrase.
   If the second entry of $A$ occurs in another word of $\nabla$, then
   $\varepsilon_\nabla (A)=1$ and $ \tau (\vert A \vert)=\vert A^\iota
   \vert =\vert A\vert$
   which contradicts the assumption on $\tau$.

  \section{Surgery and cobordism}

\subsection{Surgery}\label{unkngfr21} A nanophrase
$\nabla=(\B, (v_1\,\vert \, \cdots \, \vert \, v_k))$
  is a {\it factor} of a nanoword $(\A, w)$ if $\B\subset \A$  and $$w=x_1
v_1x_2v_2   \cdots  x_{k} v_k x_{k+1}$$ where $x_1,x_2,...,
x_{k+1}$ are words in the $\alpha$-alphabet $\C=\A-\B$. It is
understood  that the projections  $\B\to \alpha$ and $\C\to
\alpha$ are the restrictions of the projection $\A\to \alpha$.
Deleting $v_1,...,v_k$ from $w$, we obtain a nanoword   $(\C, x_1
x_2 \cdots x_{k+1} )$. When $\nabla$ is even and symmetric,    the
transformation $(\A, w)\mapsto (\C, x_1 x_2 \cdots x_{k+1} )$ is
called {\it surgery}.
 Thus,  surgery deletes an even  symmetric factor from a nanoword. The inverse
  transformation      inserts an even symmetric factor.

For example,  the first homotopy move deleting a factor $AA$  is a
surgery since the nanoword $AA$ is symmetric. The second homotopy
move deleting a factor $(AB\,\vert \,BA)$ with $\vert A \vert
=\tau  (\vert B\vert)$  is also a  surgery since the nanophrase
$(AB\,\vert \,BA)$ is even and symmetric. The third homotopy move
is neither a surgery nor an inverse to a surgery.

Here are more examples of even symmetric factors:   $(AB\,\vert
\,AB)$ with $\vert A \vert =\tau (\vert B\vert)$; $(AB\,\vert
\,CACDBD)$ with $\vert A \vert =\tau (\vert B\vert)$, $\vert C
\vert =\vert D\vert$; $(AB\,\vert \,CDEF \,\vert \, DACFBE)$ with
$\vert A \vert =\tau (\vert B\vert)$, $\vert C \vert =\tau (\vert
F\vert)$, $\vert D \vert =\tau (\vert E\vert)$.

\subsection{Group $\mathcal N_c$}\label{fi2grouppsks9}   We say that    nanowords $v,w$  (over $\alpha$)  are {\it cobordant}  and write $v \sim_c w$ if $v $ can be transformed into $w$ by a finite sequence of  moves from    the following list:

(TR) isomorphisms,   homotopy moves,   surgeries, and  the inverse
moves.

\begin{lemma}\label{l333698125mmg1} (i) Cobordism   is an equivalence relation on the class of nanowords.   Homotopic nanowords    are   cobordant.

(ii) If $v \sim_c w$, then $v^- \sim_c w^-$.

(iii) If $v_1\sim_c v_2$ and $w_1 \sim_c w_2$, then $v_1 w_1 \sim_c v_2 w_2$.
           \end{lemma}
                     \begin{proof} Claim (i) follows from the definitions. Consider a sequence of nanowords $v=v_1,v_2,..., v_n=w$ such that $v_{i+1}$ is obtained from $v_i$ by one of the moves (TR) for all $i$.   Then  $v^-_{i+1}$ is obtained from $v^-_i$ by one of the moves (TR) for all $i$.
  Therefore $v^-\sim_c w^-$. To prove (iii), consider a sequence of moves (TR) transforming $v_1$ into $v_2$ (resp.\ $w_1$ into $w_2$). Effecting these moves first on $v_1$ and then on $v_2$, we can transform $v_1v_2$ into $w_1 w_2$. Hence  $v_1 w_1 \sim_c v_2 w_2$.
   \end{proof}

Nanowords  cobordant to $\emptyset$ are said to be {\it slice}.
A    symmetric nanoword is   slice:  being its own symmetric factor it can be deleted to give  $\emptyset$.   Contractible  nanowords are  slice.   A nanoword opposite to a slice nanoword is slice. The concatenation of   two slice nanowords is slice.

  The cobordism classes of nanowords form a group
  $\mathcal N_c=\mathcal N_c(\alpha, \tau)$  with
   multiplication induced by concatenation of nanowords.
   The inverse to a nanoword $w$ in $\mathcal N_c$ is  $  w^-$ since $w  w^-$ is   symmetric and therefore slice.

  \subsection{The length norm}\label{fi2normks9}
A
$\ZZ$-valued   norm on a group $G$ is a mapping $f:G\to \{0,1,2,\ldots\}$ such that   $ f^{-1} (0)= 1$,
$f(g)=f(g^{-1})$ for all $g\in G$,   and $f(g g')\leq f(g )+f(g')$  for any $g ,g'\in G$. Such a norm  $f$ determines a  metric $\rho_f$  on $G$  by $\rho_f(g ,g') = f(g^{-1} g')$. This metric is left-invariant: $\rho_f (hg, hg')=\rho_f (g,g')$ for all $g,g', h\in G$.
We say that $f$ is {\it conjugation invariant} if  $f(h^{-1}gh)=f(g)$ for all $g,h\in G$. It is clear that if $f$ is conjugation invariant, then  $\rho_f$ is right-invariant
in the sense that $\rho_f (gh,  g'h)=\rho_f (g,g')$ for all $g,g', h\in G$.

The length of nanowords determines a $\ZZ$-valued   norm $\vert\vert \cdot \vert \vert_l $ on $\mathcal N_c$ called the {\it length norm}. Its value $\vert\vert w\vert \vert_l $ on a cobordism class of a
nanoword $w$ is half of the minimal length of a nanoword   cobordant  to $w$.  The  axioms of a norm are straightforward. In particular,   $\vert \vert w \vert \vert_l =0$ if and only if $w$
is slice. Since all nanowords of length 2 are
contractible, the length norm
does not take    value 1.   Generally speaking, the length norm is not conjugation invariant.
The associated left-invariant metric on  $\mathcal N_c$ is denoted $\rho_l$ and called the {\it length metric}.

\subsection{Push-forwards and pull-backs}\label{fi222pullbacks9} Given another
set with involution $(\overline \alpha, \overline \tau)$ and an
equivariant mapping
 $f:\overline \alpha\to \alpha$, the induced
  {\it push-forward}   transforms a  nanoword  $(\A,w)$ over $\overline \alpha$ in
   the same nanoword $(\A,w)$ with projection $\A\to \overline \alpha$ replaced by
     its composition with $f$. The push-forward  is compatible with   cobordism and  induces a group homomorphism   $f_*:\mathcal
      N_c (\overline \alpha, \overline \tau) \to \mathcal N_c(\alpha, \tau)$.
 Clearly,  $\vert \vert f_*(x) \vert \vert_l\leq \vert \vert x \vert \vert_l$
 for
  any $x\in \mathcal N_c(\overline \alpha, \overline \tau)$.

For  a $\tau$-invariant  subset  $\beta$ of $\alpha$,  the {\it
pull-back} to $\beta$ transforms any  nanoword $(\A, w)$ over
$\alpha$  in  the  nanoword over $(\beta, \tau\vert_\beta)$
obtained by deleting from both $\A$ and $w$ all letters $A\in \A$
with  $\vert A\vert \in \alpha- \beta$. This transformation is
compatible with   cobordism and  induces a group homomorphism
$\varphi_\beta:\mathcal N_c(\alpha, \tau) \to \mathcal N_c(\beta,
\tau\vert_\beta)$.  Clearly,  $\vert \vert \varphi_\beta (x) \vert
\vert_l\leq \vert \vert x \vert \vert_l$ for any $x\in \mathcal
N_c(\alpha, \tau)$. Composing the push-forward $i_*:\mathcal N_c
(\beta, \tau\vert_\beta) \to \mathcal N_c(\alpha, \tau)$ induced
by the inclusion $i:\beta \hookrightarrow  \alpha$ with
$\varphi_\beta$ we obtain the identity. Therefore $i_*$ is
injective and $\varphi_\beta$ is surjective.

 \subsection{Examples}\label{exa}
1.  For    $a,b\in \alpha$, consider the nanoword $w_{a,b}=ABAB$ with $\vert A\vert =a, \vert B \vert =b$.   If $a= b$, then
 $w_{a,b}$ is  symmetric and therefore slice. If $a= \tau (b)$, then deleting the factor $(AB\,\vert \, AB)$ we obtain $\emptyset$ so that
 $w_{a,b}$ is slice (in fact $w_{a,b}$ is  contractible for $a= \tau (b)$, see  \cite{tu3}, Lemma 3.2.2).  If  $a,b$ belong to   different   orbits of $\tau$, then $w_{a,b}$  is not slice, see
 Sect.\ \ref{groupga}. Obviously, 
$\vert\vert w\vert \vert_l\leq 2 $ and since
 $\vert\vert w\vert \vert_l \neq 0,1$, we have   $\vert\vert w\vert \vert_l= 2 $.

 2.  Pick $a,c \in \alpha$ and consider the nanoword $w= ABA CDCDB$ with   $\vert  A\vert=\vert  B\vert  =a , \vert  C\vert=\vert  D\vert =c$. Deleting the symmetric nanoword  $CDCD$ from $w$, we obtain a symmetric nanoword $ ABA B  $. Therefore $w$  is  slice.
  Note that  $w$  is not symmetric. If $a, c$ belong to different orbits of $\tau$ and $a\neq \tau (a)$, then   the pull-back of $w$ to the orbit of $a$ yields a non-contractible nanoword $ABAB$. Therefore in this case $w$ is not contractible.

 3.  The
  nanoword $  ABA CDCDB$ with   $\vert  A\vert=\tau( \vert  B\vert)$ and
  $  \vert  C\vert=\vert  D\vert  $ is slice since  $CDCD$ is symmetric
  and  $ABAB$ is contractible. The  nanoword $AB
CACDBD$ with $\vert A \vert =\tau (\vert B\vert)$, $\vert C \vert
=\vert D\vert$ is slice since the deletion of the even symmetric
factor  $(AB\,\vert \,CACDBD)$   gives $\emptyset$.

\section{Homomorphism $\gamma$ }\label{ccc56822}

We construct a group homomorphism    from $ \mathcal N_c =\mathcal N_c(\alpha, \tau)$ to a    free product of cyclic groups. This allows us to show that, generally speaking,   $ {\mathcal {N}}_c$ is non-abelian.

\subsection{Group $\Pi$ and homomorphism $\gamma$}\label{groupga}   Let $\Pi$ be the group with generators $\{z_a\}_{a\in \alpha}$ and defining relations $z_a z_{\tau
(a)}=1$ for  all  $a\in \alpha$.
For  a nanoword $ (\A, w:\widehat  n \to \A)$ over $\alpha$, set $\gamma (w)=\gamma_1 \cdots  \gamma_n \in \Pi$
  where   $\gamma_i
  =z_{\vert w(i) \vert}  $ if  $i$ numerates  the first entry of  $w(i)$ in   $w$ (that  is  if $w(i)\neq w(j)$ for $j<i$) and   $\gamma_i    =(z_{\vert w(i) \vert})^{-1}   $  if  $i$ numerates  the second  entry of  $w(i)$ in   $w$.   For example, for  $w= ABAB $ with $\vert A\vert =a\in \alpha, \vert B\vert =b\in \alpha$, we have  $\gamma(w)= z_a z_b z_a^{-1} z_b^{-1} $.

  \begin{lemma}\label{25edcolbrad}
   The element $\gamma(w)\in \Pi$ is invariant under the moves (TR) on $w$.
   The formula $w\mapsto \gamma (w)$ defines a group homomorphism
   $\gamma :   {\mathcal {N}}_c\to \Pi$ and  $\gamma({\mathcal {N}}_c)= [\Pi, \Pi]$. \end{lemma}
    \begin{proof}   It is
easy to check that $ \gamma (w) $    is invariant under
isomorphisms and homotopy moves on $w$. Let us check the
invariance under surgery. It suffices to show that for any even
symmetric factor $\nabla=  (v_1\,\vert \, \cdots \, \vert \, v_k)$
of $w$ and any $r\in \{1,..., k\}$, we have $\gamma (v_r)=1$. Fix
$r$ and set $v=v_r$. Let $n\geq 0$ be the length of $v$. By
definition, $\gamma (v)=\gamma_1 \cdots \gamma_n \in \Pi$ with
$\gamma_i $ defined by the $i$-th letter of $v$ as above. We claim
that $\gamma_i=(\gamma_{n+1-i})^{-1}$ for all $i$. This and the
assumption that $n$ is even would imply that $\gamma (v)=1$.

Pick $i\in \{1, \ldots, n\}$. Consider first the case where the letter $v(i)$ occurs in $v $
twice. Then     $\gamma_i
  =z_{\vert v(i) \vert}  $ if  $i$ numerates  the first entry of  $v(i)$ in   $v$  and
   $\gamma_i    =(z_{\vert v(i) \vert})^{-1}   $ otherwise.
Observe that  if $i$ numerates the first (resp.\ the second) entry
of $v(i)$, then by the symmetry of $\nabla$, the index $n+1-i$
numerates the second (resp.\ the first) entry of the letter
$v(n+1-i)=\iota_\nabla  (v(i))$ in $v$. Also $\varepsilon_\nabla
(v(i))=0$ and by the definition of a symmetric nanophrase, $\vert
v(i)\vert = \vert \iota_\nabla  (v(i)) \vert=\vert v(n+1-i)\vert$.
Hence,   
 $\gamma_i=(\gamma_{n+1-i})^{-1}$. If $v(i)$ occurs in $v=v_r$ only once, then by
  the symmetry,  the same is true for  $v(n+1-i)$. In particular, $v(i)\neq v(n+1-i)$.
   The symmetry implies also that
  the other  entries of these two letters in $\nabla$
   occur   in the same word $v_{r'}$ where $r'\neq r$.  Then $\gamma_i
  =(z_{\vert v(i) \vert} )^\delta $ and $\gamma_{n+1-i}
  =(z_{\vert v(n+1-i) \vert}  )^\delta $ where $\delta =  1$ if $r'>r$ and $\delta=-1$ if $r'<r$.
   Observe that  $\varepsilon_\nabla (v(i))=1$ and so 
   $\vert v(i)\vert =\tau (\vert v(n+1-i)\vert)$.
   Hence   $\gamma_i=(\gamma_{n+1-i})^{-1}$.

  The second claim of the lemma follows from the definitions. The equality  $\gamma({\mathcal {N}}_c)= [\Pi, \Pi]$ follows from \cite{tu3}, Lemma 4.1.1.
  \end{proof}

 The group $\Pi$ is a free product of the cyclic subgroups generated by
 $\{z_a\}$ and  numerated by the orbits of the involution $\tau$.   More precisely,
  $\Pi$  is a free product of   $m$ infinite cyclic
groups and   $l$ cyclic groups of  order 2 where $m$ is the number
of free  orbits of $\tau$ and $l$ is the number of fixed points of
$\tau$.
  The  commutator subgroup   $[\Pi, \Pi]$ is a free group of infinite  rank  if  $m
  \geq 2$  or $m=1$ and $l\geq 1$.  If $m=1$ and $l=0$, then $\Pi=\ZZ$ and
   $[\Pi, \Pi]=0$. If $m=0$, then    $[\Pi, \Pi]$  is a free group of rank
   $2^{l-1} (l-2)+1$.  One can  see  it  by realizing $\Pi$ as the
fundamental group of the connected sum $X$ of $l$ copies of $RP^3$
 and observing that the maximal abelian covering  of $X$ has the same
  fundamental group as  a
connected graph with $2^{l-1} l$ vertices and $2^l (l-1)$ edges.
These computations and Lemma \ref{25edcolbrad} give  the following
information on   the group ${\mathcal {N}}_c$.

\begin{theor} \label{dfroq859s}  If $\tau$ has
at least two orbits, then       ${\mathcal {N}}_c$ is  infinite.
If $\tau$ has at least two orbits and $\tau \neq \id$, then
${\mathcal {N}}_c$ is  infinitely generated. If $\tau$ has at
least three orbits or $\tau$ has two  orbits and $\tau \neq \id$,
then ${\mathcal {N}}_c$ is non-abelian.
 \end{theor}

The free product structure on  $\Pi$  allows us to detect easily whether two given  elements  of $\Pi$ are equal  or not. As an application, consider a nanoword $w=A_1 A_2 \cdots A_n$ such that $\vert A_i \vert, \vert  A_{i+1}\vert \in \alpha $ do not lie in the same  orbit of $\tau$ for $i=1, ...,n-1$. Then there are  no cancellations in the expansion  $\gamma (w)=\gamma_1 \cdots  \gamma_n \in \Pi$. This implies that such $w$ is non-slice and moreover $\vert \vert w\vert \vert_l=n/2$. For instance,
consider the nanoword $w=w_{a,b}=ABAB$ where $\vert A\vert =a\in \alpha, \vert B\vert =b\in \alpha$. By  Example \ref{exa}.1, if $a, b$ lie in the same orbit of $\tau$, then $w$ is slice. If $a,b$ do not lie in the same orbit of $\tau$, then by the criterion above,  $w$ is non-slice and  $\vert \vert w\vert \vert_l=2$.

\begin{theor} \label{dwababafroq859s}  Two non-slice nanowords $w=w_{a,b}$ and $w'=w_{a', b'}$ with $a,b,a',b'\in \alpha$ are cobordant if and only if $a=a'$ and $b=b'$.
 \end{theor}
\begin{proof} If $w\sim_c w'$, then \begin{equation}\label{iko}z_a z_b z_a^{-1} z_{b}^{-1}=\gamma(w)=\gamma(w')= z_{a'} z_{b'} z_{a'}^{-1} z_{b'}^{-1}.\end{equation}  The non-sliceness of $w$ (resp.\ $w'$)  implies that $a,b$ (resp.\ $a', b'$)  belong to different orbits of $\tau$. Therefore there are no cancellations in the expansions for $\gamma (w), \gamma (w')$ above. Formula (\ref{iko}) implies then that $z_a= z_{a'}$ and $z_b= z_{b'}$. Therefore $a=a'$ and $b=b'$.
 \end{proof}

\subsection{Homomorphism $\widetilde  \gamma$}\label{refinga} The homomorphism $\gamma$ admits the  following  refined version.
Let $\widetilde  \Pi$ be the group with generators $\{\widetilde 
z_a\}_{a\in \alpha}$ and defining relations
  $ \widetilde  z_a \widetilde  z_{\tau(a)} \widetilde  z_b=\widetilde  z_b \widetilde  z_a \widetilde  z_{\tau
(a)} $ for   all $a,b\in \alpha$.      The formula $\widetilde 
z_a\mapsto  z_a$ defines a  projection $\widetilde  \Pi \to \Pi$ which
makes  $\widetilde  \Pi$ into a central extension of $\Pi$. Replacing
$z $ with $\widetilde  z $ in the definition of $\gamma$, we obtain  a
lift of $\gamma$ to a group homomorphism $\widetilde  \gamma: {\mathcal
{N}}_c\to \widetilde  \Pi$.

The homomorphisms $\gamma$  and $\widetilde  \gamma$ are not injective.
 For   example, as we shall see
in Sect.\ \ref{remlexamps799}.1,  the nanoword $  w= ABCBAC $ with
$\vert A\vert=\vert B\vert =\vert C\vert \neq \tau (\vert A\vert)
$ is non-slice    but obviously    $\widetilde  \gamma (w)=1$.

    \section{   $\alpha$-pairings and their cobordism}\label{fi:g444599912}

   We now turn to the main theme of this paper: a study of cobordisms of
   nanowords via a study of the linking properties of the letters. In this and the next sections we introduce a purely
   algebraic theory of $\alpha$-pairings; it will be applied to nanowords in later sections.

     Fix   an associative (possibly,  non-commutative)    ring   $R$ and a left $R$-module~$\pi$. The
      module $\pi$ will be the  target of all $\alpha$-pairings.

\subsection{$\alpha$-pairings}\label{fi222799}    An  {\it $\alpha$-pairing}   is  a set $S$ endowed with a distinguished element
$s\in S$ and
mappings  $S-\{s\}\to \alpha$ and   $e:S\times
 S
\to \pi$.   The conditions on $S $ can be rephrased by saying that $S$ is a disjoint union of an $\alpha$-alphabet $S^\circ=S-\{s\}$ and a distinguished element $s$. The image of   $A\in S^\circ$ under the
projection  to $\alpha $
 is denoted $\vert A\vert$. The pairing $e:S\times
 S
\to \pi$ uniquely extends to a  bilinear form  $RS \times RS \to \pi$ where $RS$
 is the free $R$-module with basis $S$. This form  is   denoted by $\widetilde  e$
 or, if it cannot lead to a confusion,  simply by   $e$.  Every   $A\in S $ determines a basis vector in $RS $ denoted by the same symbol $A$.

An {\it isomorphism} of
$\alpha$-pairings  $(S_1, s_1,e_1), (S_2, s_2,e_2)$  is a bijection $S_1\to S_2$ transforming $s_1, e_1$ into $s_2, e_2$, respectively, and inducing an isomorphism of
$\alpha$-alphabets $S_1^\circ\to S_2^\circ$. Isomorphism of $\alpha$-pairings is denoted   $\approx$.

For each $\alpha$-pairing $p=(S, s,e)$, we  have     the {\it opposite}
$\alpha$-pairing    $  p^-=(  S,  s ,    e^-) $ where   $ e^-(  A,   B)= -e(A,B)$ for $A,B\in S $.

 \subsection{ Hyperbolic  $\alpha$-pairings}\label{fi:gdr5999}
 Consider  an  $\alpha$-pairing
$p=(S,s,e)$.   A vector $x\in RS$ is  {\it short} (with respect to
$p$)  if   $x \in S^\circ  \subset S\subset RS$  or $x =  A+B$ for
distinct   $A,B \in S^\circ$ with $\vert A\vert = \vert B\vert $
or $x =  A-B$ for   distinct  $A,B \in S^\circ$ with  $\vert
A\vert =\tau(\vert B\vert)$.   Note that if   $\vert A\vert =
\vert B\vert=\tau(\vert B\vert)  $ then both $A+B$ and $A-B$ are
short.

  A  {\it filling} of
$p$ is  a finite family  of  vectors  $  \{\lambda_i\in RS \}_i$  such that one of the  $\lambda_i$'s  is equal to $s$, all the other $\lambda_i$ are short,  and  every element of  $S^\circ$ occurs in exactly one of   $  \lambda_i $ with non-zero coefficient (this
 coefficient is then
$\pm 1)$.   For example, the family  $\{A\}_{A\in S } $ is a
filling
 of
$p$. It is called the {\it tautological} filling.

   A filling  $ \{\lambda_i\}_i$ of
$p$ is {\it annihilating} if      $\widetilde 
e(\lambda_i,\lambda_j)=0$ for all $i,j$.  The $\alpha$-pairing $p$
is {\it hyperbolic} if it has an annihilating filling. Since the
number of fillings of  $p$ is finite, one can detect in a finite
number of steps whether  $p$ is hyperbolic or not.

 If
an $\alpha$-pairing   is  hyperbolic, then the opposite
$\alpha$-pairing and  all  isomorphic $\alpha$-pairings    are
hyperbolic.

\subsection{Summation of $\alpha$-pairings}\label{summ799}
The {\it   sum} $p_1\oplus p_2$ of    $\alpha$-pairings $p_1=(S_1,
s_1, e_1)$ and $p_2=(S_2, s_2, e_2)$ is the $\alpha$-pairing  $  (
S=S_1^\circ\amalg  S_2^\circ \amalg \{s\} , s, e) $ where
$e:S\times
 S
\to \pi$  is defined as follows. Consider the bilinear form
$\widetilde  e_i:RS_i\times RS_i \to \pi$ extending  $e_i$ for $i=1,2$.
The direct sum  $\widetilde  e_1 \oplus\widetilde   e_2$ is a bilinear form
on $RS_1\oplus RS_2$. Consider the $R$-linear embedding
$f:RS\hookrightarrow RS_1\oplus RS_2$ which extends the embeddings
$S_i^\circ \hookrightarrow S_i \subset RS_i$ with $i=1,2$ and
sends $s\in S$ to $s_1\oplus s_2$. For  $x,y \in S$, set
$$e(x,y)=(\widetilde  e_1 \oplus\widetilde   e_2) (f(x), f(y))\in \pi.$$
 The values of $e$ can be computed explicitly:
 $e(S_1^\circ, S_2^\circ)=e(S_2^\circ, S_1^\circ)=0$;
   $e\vert_{  S_i^\circ}=e_i\vert_{  S_i^\circ} $;
    $e(s,s)=e_1(s_1,s_1)+e_2 (s_2, s_2)$;
$e(A,s)=e_i(A,s_i), e(s, A)= e_i (s_i, A)$ for $i=1,2$ and   $A\in
S_i^\circ$.  It is clear  that  the bilinear extension $\widetilde  e
:RS\times RS\to \pi$ of $e$ is obtained by pushing back $\widetilde 
e_1 \oplus\widetilde  e_2$ along $f$. Observe that   the projection $S^\circ =S_1^\circ\amalg  S_2^\circ \to \alpha$ is the disjoint union of the given projections
$S_1^\circ \to \alpha$ and $S_2^\circ \to \alpha$.
  In these constructions  we assume
$S_1, S_2$ to be disjoint; if it is not the case, replace $p_1$ by
an isomorphic $\alpha$-pairing and proceed as above. 

We shall routinely describe fillings  of $p_1\oplus p_2=(S,s,e)$
in terms of their images  under the embedding $f:RS
\hookrightarrow RS_1\oplus RS_2$. By abuse of the language, the
image of a filling of $p_1\oplus p_2$ under $f$ will sometimes be
called a filling of $p_1\oplus p_2$. A finite family of vectors
$\lambda= \{\lambda_i \}_i \subset RS_1 \oplus RS_2$ is the image
of a filling of $p_1 \oplus p_2$ if and only if it satisfies the
following conditions: $\lambda$ consists of $s_1+s_2$ and  
vectors of the form $A\in S_1^\circ \cup S_2^\circ $ or $ A\pm \, B$, where $A, B$ are distinct elements of $S_1^\circ \cup
S_2^\circ$ with $  \vert A\vert \in \{\vert B\vert, \tau
(\vert B\vert)\}$, and  the sign in front of $B$ is necessarily $+$ if $\vert
A\vert =\vert B\vert \neq \tau (\vert B\vert)$ and  is necessarily $-$ if $\vert
A\vert = \tau (\vert B\vert)\neq \vert B\vert $; every element of
$S_1^\circ \cup S_2^\circ $ occurs in exactly one    $ \lambda_i $
with  non-zero coefficient (equal then to $\pm 1$).
 It is clear that    $\lambda $ corresponds to an
annihilating filling of $p_1\oplus p_2$  if and only if $(\widetilde 
e_1 \oplus \widetilde  e_2)(\lambda_i, \lambda_j)=0$ for all $i,j$.

The sum  $p_1\oplus p_2$ is well-defined  up to isomorphism and
$p_1\oplus p_2 \approx p_2\oplus p_1$. The sum of hyperbolic
$\alpha$-pairings is hyperbolic.

\subsection{ Cobordism of  $\alpha$-pairings}\label{fi:999hjf9} We say that  two  $\alpha$-pairings
$p_1 , p_2$ are {\it   cobordant} and write $p_1\simeq_c  p_2$ if
the  $\alpha$-pairing  $p_1\oplus  p_2^-$ is  hyperbolic.

\begin{lemma}\label{l:gg1}  Cobordism   is an equivalence relation on the class of $\alpha$-pairings.  Isomorphic  $\alpha$-pairings are cobordant.        \end{lemma}
                     \begin{proof}  Consider an isomorphism  $\varphi:S_1\to S_2$ of   $\alpha$-pairings  $p_1=(S_1, s_1,e_1), p_2=(S_2, s_2,e_2) $.   The set of vectors
 $\{A+\varphi (A)\}_{A\in S_1^{\circ}}\cup \{s\}$  is  an annihilating filling of   $p_1\oplus  p_2^-$.     Therefore this pairing
  is   hyperbolic and    $p_1\simeq_c  p_2$. In particular,  $p_1\simeq_c  p_1$.

If $p_1\oplus  p_2^-$  is  hyperbolic, then so is its opposite
  $p_1^-\oplus  p_2\approx p_2 \oplus  p_1^-$.  This implies the  symmetry of   cobordism.

Let us prove the transitivity.   Let $p_1,p_2,p_3$ be
$\alpha$-pairings such that $p_1 \simeq_c  p_2 \simeq_c  p_3$. We
verify that $p_1\simeq_c  p_3$.     Let $p_k=(S_k,s_k,e_k)$ for
$i= 1,2,3$ and   $p'_2=(S'_2,s'_2,e'_2)$ be a  copy of $p_2$
 where
$S'_2=\{A' \,\vert \, A\in S_2\}$, $s'_2=(s_2)'$,
and  $e'_2(A',B')=e_2(A,B)$ for all $A,B\in S_2$.   Replacing the pairings by isomorphic ones, we can assume that the sets $ S_1,S_2,S'_2, S_3$ are   disjoint.    Let  $\Lambda_1 , \Lambda_2, \Lambda'_2, \Lambda_3 $ be  free $R$-modules
with bases
$S_1, S_2, S'_2,  S_3$, respectively. Set $\Lambda=  \Lambda_1\oplus  \Lambda_2\oplus \Lambda'_2\oplus
\Lambda_3 $.
There is  a  unique   bilinear form $e=e_1\oplus e_2^-\oplus  e'_2 \oplus e_3^-:\Lambda\times \Lambda \to \pi$ such that the
sets
$S_1, S_2, S'_2, S_3\subset \Lambda$ are  orthogonal with respect to $e$ and the restrictions of $e$ to these  sets  are equal
 to $e_1, e_2^-, e'_2,e_3^-$, respectively.

 Let $\Phi$   be the submodule  of
$ \Lambda_2\oplus \Lambda'_2$ generated by    the vectors
$\{A+A'\}_{A\in S_2}$.  Set $L= \Lambda_1\oplus  \Phi \oplus \Lambda_3\subset \Lambda$.     The
projection $q:
L  \to \Lambda_1\oplus    \Lambda_3$  along $\Phi$ transforms $e$ into $e_1\oplus e_3^-$.
Indeed, for any  $A_1,B_1\in S_1, A , B\in S_2, A_3,B_3\in S_3$,
$$e(A_1+A  +A' +A_3, \,\,B_1+ B + B'+B_3)$$
$$=
e_1(A_1 ,B_1) - e_2 (A ,  B)+  e'_2  (A', B')  -e_3(A_3,B_3)=
e_1(A_1 ,B_1)     -e_3(A_3,B_3).$$

Pick    a  filling  $\lambda= \{\lambda_i\}_i\subset \Lambda_1
\oplus \Lambda_2$ of $p_1\oplus  p_2^- $ and   a   filling $ \mu=
\{\mu_j\}_j\subset \Lambda'_2\oplus \Lambda_3$    of $p'_2\oplus
p_3^- $. Consider the $R$-modules $V_\lambda\subset
\Lambda_1\oplus \Lambda_2$ and $V_\mu\subset \Lambda'_2\oplus
\Lambda_3$ generated respectively by $ \{\lambda_i\}_i$ and    $
\{\mu_j\}_j$. Below we construct a  finite  set $\psi\subset
(V_\lambda+V_\mu)\cap L$ such that  $q(\psi) \subset
\Lambda_1\oplus    \Lambda_3$ is a filling of $p_1 \oplus  p_3^-
$.  Choosing $\lambda, \mu$ to be annihilating fillings, we obtain
that $e(V_\lambda , V_\lambda )=  e( V_\mu,  V_\mu)= 0$ and
therefore $e(V_\lambda+ V_\mu, V_\lambda+ V_\mu)= 0$. Since
$q: L  \to \Lambda_1\oplus    \Lambda_3$    transforms $e$ into
$e_1\oplus e_3^-$, the filling $q (\psi) $ of $p_1\oplus  p_3^- $
is annihilating.  Hence $p_1\simeq_c  p_3$.

To define $\psi$, we first  derive
from the    filling $ \lambda$   a 1-dimensional manifold
$\Gamma_{\lambda}$. If $\lambda_i= A\pm B$ with distinct  $A, B
\in S_1^\circ \cup S_2^\circ$, then $\lambda_i $ yields   a
component of $\Gamma_{\lambda}$ homeomorphic to $[0,1]$ and
connecting   $A$ with $B$. (By the definition of a filling, $\vert
A\vert \in \{\vert B\vert ,\tau (\vert B \vert)\}$). If
$\lambda_i=A\in S_1^\circ \cup S_2^\circ$, then $\lambda_i $
yields  a component of $\Gamma_{\lambda}$ homeomorphic to  $[0,
\infty)$ where $0$ is identified with  $A$. The vector   $
s_1+s_2\in \lambda $   does  not contribute to $\Gamma_{\lambda}$.  The
definition of a filling implies that $\partial \Gamma_{\lambda}=
S^\circ_1\cup {S^\circ_2} $.  The filling $\mu$ similarly gives
rise to a 1-dimensional manifold $\Gamma_{\mu}$ with $\partial
\Gamma_{\mu}= (S'_2)^\circ \cup {S^\circ_3}  $.  We can assume
   the manifolds $\Gamma_{\lambda}$ and $\Gamma_{\mu}$ to be 
disjoint. Gluing them along $S_2^\circ\approx (S'_2)^\circ$,  we
obtain     a 1-dimensional manifold  $\Gamma=\Gamma_{\lambda} \cup
\Gamma_{\mu}$     with $\partial \Gamma= S^\circ_1 \cup
{S^\circ_3} $.

We  associate with each component  $K$ of $\Gamma$ with  $\partial K\neq \emptyset$ a vector $\psi_K\in \Lambda$. Suppose first that $K$ is compact and let
 $A, B \in  S^\circ_1 \cup    {S^\circ_3}$ be its endpoints.   The 1-manifold  $K$   is glued from several components of
$\Gamma_{\lambda}\amalg  \Gamma_{\mu}$ associated with certain  vectors  $\lambda_i ,  \mu_j$ (the components of $\Gamma_{\lambda}$  are intercalated in $K$ with the components  of~$\Gamma_{\mu}$). We define $\psi_K$ as an algebraic  sum
of these vectors  $\sum_i \varepsilon_i \lambda_i +\sum_j \eta_j \mu_j$, where the signs $\varepsilon_i, \eta_j=\pm 1$ are defined from the following two conditions:
 $\psi_K\in L$ and
$q(\psi_K) = A\pm B $. The  signs $\varepsilon_i, \eta_j$ are determined by induction
moving along $K$  from $A$ to $B$.  For example, if $K$ is a union
of a component of  $\Gamma_{\lambda}$ connecting  $A\in S_1^\circ$
to  $C\in S_2^\circ$ and corresponding to $\lambda_{i}=A+C $  with
a component of  $\Gamma_{\mu}$ connecting  $C'$ to $B\in
S_3^\circ$  and corresponding to $\mu_{j}=-C' \pm B $, then
$\psi_K=\lambda_{i}-\mu_{j}$. If    $\mu_{j}=C'\pm B $,  then
$\psi_K=\lambda_{i}+\mu_{j}$. In both examples  $\vert A\vert =
\vert C\vert =\vert C'\vert \in \{ \vert B\vert, \tau (\vert
B\vert)\}$.

  An easy inductive argument shows that $\vert A\vert
\in \{\vert B\vert, \tau(\vert B \vert)\}$ for any compact
component $K$ of $\Gamma$ with endpoints $A,B$. We claim that
$q(\psi_K)= A\pm   B$ is short, i.e., that the sign $\pm $
satisfies the requirements in the definition of a short vector.
  If $\tau (\vert B\vert)=\vert B\vert$, then there is nothing to prove  since both $+$ and $-$ satisfy these  requirements. If  $\tau (\vert B\vert)\neq \vert B\vert$, then this claim is obtained by a  count of  minuses in the sequence of vectors $\lambda_i,  \mu_j$  corresponding to the  components of
$\Gamma_{\lambda}\amalg  \Gamma_{\mu}$ forming $K$. Note that the vectors  $\psi_K$ determined as above  by moving along $K$  from $A$ to $B$ and from $B$ to $A$  may differ;  we take any of them.
If  $K$ has only  one endpoint $A$, then $\psi_K$ is similarly defined as an algebraic sum of the vectors associated with the  components of
$\Gamma_{\lambda}\amalg  \Gamma_{\mu}$ forming $K$, where the signs are determined inductively from two  conditions:  $\psi_K \in L$ and $ q(\psi_K)  = A $.
It follows from the definitions that in all cases   $\psi_K\in (V_\lambda+V_\mu)\cap L$.

Set
$\psi_0=s_1+s_2+s'_2+ s_3
\in
  (V_\lambda+V_\mu)\cap L$ and
 set $\psi=\{\psi_0\}\cup \{\psi_K\}_K$, where $K$ runs over the  components of $\Gamma$ with non-void boundary.
  All  vectors in the family $q(\psi) $ besides  $q(\psi_0)=s_1+s_3$  are short and all  elements  of  $S_1^\circ \cup S_3^\circ  $ occur  in exactly one of these  vectors   with non-zero coefficient.     This means that $q(\psi)$ is a
filling of
$p_1\oplus  p_3^-
$ as  required.
 \end{proof}

 \subsection{The group ${\mathcal{P}}$}\label{ermlp99}
The cobordism classes of $\alpha$-pairings form an abelian group with respect to summation  $\oplus$. This group is denoted ${\mathcal{P}}={\mathcal{P}}(\alpha, \tau, \pi)$. The neutral element of  ${\mathcal{P}}$ is the  class of  the
 {\it  trivial $\alpha$-pairing}  $ ( S=\{s\}, s, e(s, s)=0)$. The opposite  in ${\mathcal{P}}$  to the class of an $\alpha$-pairing  $p$  is the class of  $p^-$.

 An $\alpha$-pairing $p=(S,s,e)$ is {\it skew-symmetric} if $e(A,A)=0$ and $e(A,B)=-e(B,A)$ for all $A,B\in S$. In particular, we must have $e(s,s)=0$ and $e(s,B)=-e(B,s)$ for all $B\in S^\circ$.
It is clear that the sum of  skew-symmetric  $\alpha$-pairings is  skew-symmetric and the  $\alpha$-pairing opposite to a skew-symmetric one is itself skew-symmetric. Therefore the cobordism classes of skew-symmetric $\alpha$-pairings form a subgroup of ${\mathcal{P}}$. It is  denoted  ${\mathcal{P}}_{sk}={\mathcal{P}}_{sk}(\alpha, \tau, \pi)$.

 \subsection{Normal $\alpha$-pairings}\label{ermgiofprnisdnf99} Although we
 shall not need  it in the sequel, we briefly discuss   so-called normal $\alpha$-pairings.  An $\alpha$-pairing $p=(S,s,e)$ is {\it normal} if  $e(s,s)=0$.
The cobordism classes of normal  $\alpha$-pairings form a subgroup
of ${\mathcal{P}}$  denoted  ${\mathcal{P}}_{n}$.
Clearly, ${\mathcal{P}}_{sk}\subset {\mathcal{P}}_{n}$.
The following  lemma computes  ${\mathcal{P}}$
from ${\mathcal{P}}_{n}$. Denote by $\underline R$ the underlying
 additive group of $R$. For $r\in R$, denote by $i(r)$  the $\alpha$-pairing $(S=\{s\},s, e(s,s)=r)$.

 \begin{lemma}\label{cooplussst} The formula $r\mapsto i(r)$ defines an injective group homomorphism $i: \underline R\to \mathcal P$ and
  $\mathcal P= {\mathcal{P}}_{n}\oplus i(\underline R) $.\end{lemma}
\begin{proof} The additivity of  $i$  follows from the definitions.
For an $\alpha$-pairing  $p=(S,s,e)$, set $r_p=e(s,s)\in R$ and consider
the $\alpha$-pairing $p'=(S,s,e')$ where $e'(s,s)=0$ and
 $e'(A,B)=e(A,B)$ for all pairs $A,B\in S$ distinct from
 the pair $(s,s)$. It follows from the definitions  that $p\approx p'\oplus i(r_p )$. Therefore  $\mathcal P={\mathcal{P}}_{n}+i (\underline R)$. Observe that if two $\alpha$-pairings $p_1=(S_1,s_1, e_1), p_2=(S_2, s_2, e_2)$ are cobordant, then the vector $s_1 +s_2$ belongs to an annihilating filling of  $p_1\oplus  p_2^-$ and therefore
$$0= (e_1 \oplus e_2^-) (s_1+s_2, s_1 +s_2)=e_1(s_1,s_1) - e_2(s_2,s_2)=
r_{p_1}-r_{p_2}.$$
Hence the formula $p\mapsto r_p$ defines a group homomorphism
$r:\mathcal P\to \underline R$. Clearly, $r\circ i=\id$ and $r({\mathcal{P}}_{n})=0$. Therefore
$i$ is an injection and $ {\mathcal{P}}_{n}\cap i(\underline R)=0$.
\end{proof}

\section{Cobordism invariants of $\alpha$-pairings}

We give two constructions  of cobordism invariants of $\alpha$-pairings.
Fix as above  a left module $\pi$  over a ring   $R$.

\subsection{The polynomial $u$}\label{eruuuuusp99} Let $I$ be
 the free $R$-module with basis $\{\delta_g\}_{g\in \pi-\{0\}}$.
Let $J$ be   the  submodule of $I$  generated by the vectors   $\{\delta_g+\delta_{-g}\}_{g\in \pi-\{0\}}$.
 For an  $\alpha$-pairing $p=(S,s,e)$  and any   $a\in \alpha$,
 set $$[a]_p= \sum_{ A \in S^\circ, \vert A\vert =a, { e(A,s)} \neq 0}
  \delta_{ e(A,s)} \,\,
\in I.$$ The {\it  $u$-polynomial} $u^p$ of $p$ is the function on
$ \alpha$
  defined  as follows: for $a\in \alpha$ with $\tau (a)\neq a$,  $$
  u^p(a)=  [a]_p-[\tau(a)]_p \,(\modu J)\in I/J$$  and  for $a\in \alpha$ with $\tau
  (a)=a$,
  $$
  u^p(a)=  [a]_p \,(\modu J+2I)\in I/(J+2I).$$  Clearly, $u^p(a)=-u^p(\tau(a))$ for all $a\in \alpha$. The
function $u^p$ was introduced in  \cite{tu3} without factorization
by $J$. This factorization is needed here  to ensure the following
theorem.

   \begin{theor}\label{coppprvplo}   $u^p $ is  an additive  cobordism invariant
    of  $p$.     \end{theor}
\begin{proof}   It follows from the definitions that $u^p $ is   additive and $u^{-p}=-u^p$ for all $p$. It remains to show that if $p$ is hyperbolic, then $u^p(a)=0$ for all $a\in \alpha$.  Pick      an annihilating filling
$\{\lambda_i\in RS\}_i$ of  $p$ and let $\lambda_0$ be the vector
of this filling equal to $s$.        If $\lambda_i= A\in S^\circ$,
 then the equalities
  $e(A,s)=e (\lambda_i, \lambda_0)=0$ imply  that   $A$  contributes 0
   to   $[b]_p$ for all  $b\in \alpha$.
Therefore $A$ contributes $0$ to  $u^p(a)$.

Consider a vector of this filling  $\lambda_i= A\pm B$ with $A,
B\in S^\circ$. Recall that $\vert A\vert \in \{\vert B \vert, \tau
(\vert B \vert) \}$. The condition $e(\lambda_i,\lambda_0)=0$
implies that $ e(A,s) =\mp \,  e(B,s)$. If  $e(A,s)=0$, then
$e(B,s)=0$ so that both $A$ and $  B$ contribute   $0$ to $[b]_p$
for all $b\in \alpha$ and hence contribute  $0$ to $u^p(a)$. If
$\vert A\vert \notin \{a,\tau (a)\}$, then $\vert B\vert \notin
\{a,\tau (a)\}$ so that   both $A$ and $B$ contribute  $0$ to
$[a]_p$, $[\tau (a)]_p$ and $u^p(a)$. Suppose from now on that  $e
(A,s)\neq 0$ and $\vert A\vert \in \{a,\tau (a)\}$. If $a=\tau
(a)$, then  the equality  $ e(A,s) =\mp \, e(B,s)$ implies that
$\delta_{ e(A,s)}+\delta_{ e(B,s)}\in J+2I$ and therefore the pair
$A,B$ contributes $0$ to $u^p(a)=[a]_p\,(\modu J+2I)$. Suppose
that $a\neq \tau (a)$. If $\vert A\vert =\vert B\vert =a$, then
$\lambda_i= A + B$,  $e(A,s) =-  e(B,s)$,  and $A, B$ contribute
$\delta_{e(A,s)}+\delta_{e (B,s)}\in J$ to $[a]_p$ and $0$ to
$[\tau(a)]_p$. Hence  $A, B$ contribute   $0$ to  $u^p(a)$.    If
$\vert A\vert =a, \vert B\vert =\tau (a)$, then   $\lambda_i= A -
B$,  $e(A,s) =   e(B,s)$  and $A, B$ contribute $\delta_{e(A,s)} $
to $[a]_p$ and $\delta_{e(B,s)}$ to $[\tau(a)]_p$. Hence  $A, B$
contribute   $0$ to  $u^p(a)$. The cases where $\vert A\vert =\tau
(a)$ and  $\vert B\vert =\tau (a)$ or $ \vert B\vert =a$ are
similar.
 Since every letter of $ S^\circ$ appears in exactly
 one $\lambda_i$, summing up the contributions of all letters to $u^p(a)$ we obtain   $u^p(a)=0$.
     \end{proof}

\subsection{Genus of $\alpha$-pairings}\label{u125452}
Let $F$ be a
 commutative $R$-algebra without zero-divisors.
    Fix  an $R$-module homomorphism $\varphi: \pi\to F$.  For an
  $\alpha$-pairing
$p=(S,s,e)$, we define   a non-negative half-integer $\sigma_\varphi(p)$
as follows.
 Consider the  bilinear pairing $\varphi \, e=\varphi\, \widetilde  e:RS \times
 RS \to F$.  For a filling $\lambda=\{\lambda_i\}_i$ of
$p$, the matrix $(\varphi  e(\lambda_i, \lambda_j  ))_{i,j}$ is a
square matrix over   $F$.  Let $$\sigma_\varphi
(\lambda)= (1/2) \, \rk (\varphi  e(\lambda_i, \lambda_j  ))_{i,j} \in \frac{1}{2}\ZZ$$ be half of its rank.  The rank $\rk$ for   matrices and  bilinear forms over $F$  is defined by extending $F$ to  its quotient field and using the standard definitions for the latter.      Set
$$\sigma_\varphi(p)=\min_\lambda \sigma_\varphi (\lambda), $$
where $\lambda$ runs over all fillings of~$p$. The number $\sigma_\varphi(p)$ is called the {\it $\varphi$-genus} of~$p$. It is obvious that this number is invariant under  isomorphisms
 of  $\alpha$-pairings and  $\sigma_\varphi(p^-)=\sigma_\varphi(p)$.
   If  $ p $ is hyperbolic, then $\sigma_\varphi(p)=0$. If $p$ is
   skew-symmetric, then the matrix
   $(\varphi  e(\lambda_i, \lambda_j  ))_{i,j}$ is skew-symmetric, so that
    $\sigma_\varphi (\lambda) \in \ZZ$ for all $\lambda$ and  $\sigma_\varphi(p) \in \ZZ$.

 \begin{lemma}\label{th:verdic3}  For any  $\alpha$-pairings  $   p_1, p_2, p_3$,
\begin{equation}\label{trian}    \sigma_\varphi (p_1\oplus    p_2^- )
+ \sigma_\varphi (p_2 \oplus  p_3^-)\geq \sigma_\varphi(p_1\oplus
p_3^-).\end{equation}
\end{lemma}
\begin{proof}    We use  notation introduced in the proof of Lemma \ref{l:gg1}. Pick    a  filling  $\lambda= \{\lambda_i\}_i\subset \Lambda_1 \oplus \Lambda_2$ of
$p_1\oplus  p_2^- $ such that $\sigma_\varphi (p_1 \oplus   p_2^- )=\sigma_\varphi (\lambda)$. Pick    a   filling $ \mu=
\{\mu_j\}_j\subset
\Lambda'_2\oplus \Lambda_3$    of
$p'_2 \oplus  p_3^- $ such that $\sigma_\varphi (p'_2 \oplus  p_3^-)=\sigma_\varphi (\mu)$. Consider the $R$-modules $V_\lambda\subset \Lambda_1\oplus \Lambda_2$ and $V_\mu\subset \Lambda'_2\oplus \Lambda_3$ generated respectively by $ \{\lambda_i\}_i$
and    $ \{\mu_j\}_j$.  Recall  the
projection $q:
L=  \Lambda_1\oplus \Phi\oplus   \Lambda_3 \to \Lambda_1\oplus    \Lambda_3$   transforming  $e=e_1\oplus e_2^-\oplus  e'_2 \oplus e_3^-$ into $e_1\oplus e_3^-$.
The proof of Lemma \ref{l:gg1} yields a  finite  set $\psi\subset (V_\lambda+V_\mu)\cap L$ such that  $q(\psi) \subset
\Lambda_1\oplus    \Lambda_3$ is a filling of
$p_1\oplus  p_3^- $.  Denote by $V$ the $R$-submodule of $\Lambda_1\oplus    \Lambda_3$ generated by $q(\psi)$.
Clearly,
$$\sigma_\varphi(p_1\oplus  p_3^-)\leq  \sigma_\varphi (q(\psi))= (1/2)  \rk ((\varphi   e_1\oplus \varphi e_3^-) \vert_V)= (1/2)  \rk (\varphi e \vert_{q^{-1}(V)}).$$
Observe that $q^{-1}(V)\subset  (V_\lambda +V_\mu)\cap L + \Ker q$
and that $\Ker q=\Phi$ lies in both the left and the right
annihilators of $e\vert_L$. Therefore
$$(1/2)  \rk (\varphi e \vert_{q^{-1}(V)})  \leq
(1/2)
\rk (\varphi  e
\vert_{(V_\lambda+V_\mu)\cap L} )
 \leq   (1/2) \rk (\varphi  e
\vert_{V_\lambda +V_\mu})$$
$$=  (1/2)  \rk (\varphi  e \vert_{V_\lambda}) +
(1/2) \rk (\varphi  e \vert_{V_\mu})=\sigma_\varphi (\lambda) +
\sigma_\varphi (\mu)$$ $$=\sigma_\varphi (p_1\oplus p_2^- ) +
\sigma_\varphi (p'_2 \oplus  p_3^-)=\sigma_\varphi (p_1\oplus p_2^- ) +
\sigma_\varphi (p_2 \oplus  p_3^-) .$$
\end{proof}

\begin{theor} \label{859s}    The $\varphi$-genus   of an  $\alpha$-pairing  is a cobordism invariant.
 \end{theor}
\begin{proof}
If   $p_1\simeq_c  p_2$, then   $p_1\oplus  p^-_2$ is hyperbolic
and
  $\sigma_\varphi(p_1 \oplus   p^-_2)=0 $.   Applying the previous lemma to the triple  $(p_1, p_2, p_3)$ where  $p_3=(\{s\}, s, e=0)$ is  the trivial $\alpha$-pairing,  we obtain
$\sigma_\varphi(p_2)\geq \sigma_\varphi(p_1)$.      By symmetry,
$\sigma_\varphi(p_1)= \sigma_\varphi( p_2)$.
 \end{proof}

  \subsection{Remark}\label{fi:grema9568}
   For any $\alpha$-pairings $p_1=(S_1, s_1, e_1),p_2=(S_2, s_2, e_2)$, \begin{equation}\label{w3}
 \sigma_\varphi ( p_{1})+  \sigma_\varphi ( p_{2})\geq  \sigma_\varphi (p_1 \oplus  p_2).\end{equation}
 This can be deduced  from  (\ref{trian}) by choosing there $p_2$
to be the trivial $\alpha$-pairing. A direct proof of (\ref{w3})
goes by taking the  union of a filling  $\lambda^{(1)}$  of $p_1$
with a
 filling  $\lambda^{(2)}$ of $p_2$ and replacing in this union the elements $s_1, s_2$ by
 $s_1+s_2$. This gives   a filling  $\lambda$ of $p_1\oplus p_2$ such that
 $\sigma_\varphi(\lambda^{(1)})+
 \sigma_\varphi(\lambda^{(2)})\geq \sigma_\varphi(\lambda)$. Hence  (\ref{w3}).
In general,  (\ref{w3}) is not an equality. This is
   clear already
    from the fact that the $\varphi$-genus takes only non-negative values and
    annihilates $p\oplus p^-$ for any $\alpha$-pairing~$p$.

\section{Homomorphism ${\mathcal{N}}_c\to {\mathcal{P}}_{sk}$}

\subsection{The group $\pi$}\label{fipipipiib99} From now on,
unless explicitly stated to the contrary,  $\pi=\pi(\alpha, \tau)$
is the abelian group with generators $\{a\}_{a\in \alpha}$ and
defining relations $ a + \tau (a)=0 $ for all ${a\in \alpha}$.
This group is the abelianization of the group $\Pi$ considered in
Sect.\ \ref{ccc56822}. Clearly, $\pi$
  is a direct sum   of copies of $\ZZ$ numerated by the  free  orbits of $\tau$
  and
 copies of $\ZZ/2\ZZ$ numerated by the  fixed points  of $\tau$.
 The group $\pi$, considered as a module over  $R=\ZZ$,
  will be   the target  of
  $\alpha$-pairings. Note that  in
\cite{tu3}   the group operation in $\pi$  is written  
multiplicatively rather than
   additively as here.

\subsection{$\alpha$-pairings of nanowords}\label{fi2ampsortnib99} By \cite{tu3}, each nanoword
$(\A,w:\widehat  n \to \A)$ gives rise to a skew-symmetric
$\alpha$-pairing $p(w)=(S=\A\cup \{s\}, s, e_w:S\times S\to \pi)$
with target $\pi=\pi(\alpha, \tau)$. Here the projection $S^\circ
=\A \to \pi$ is determined  by the structure of an
$\alpha$-alphabet in $\A$.  Recall   the definition of $e_w$.
First, for  any $A, B\in \A$, we define an integer $n_w(A,B) $ to
be $+1$ if   $w=\cdots A \cdots  B\cdots  A \cdots  B  \cdots $,
to be $-1$ if  $w=\cdots  B\cdots  A \cdots  B \cdots  A  \cdots $
and to be $0$ in all other cases.  Given $A\in \A$, denote by
$i_A$ (resp.\ $j_A$) the minimal (resp.\ the maximal) element of
the 2-element set $w^{-1} (A)\subset \widehat  n$. For $A,B\in \A$, set
$$ A\circ_w B= \sum_{D\in \A, i_A<i_D<j_A \, {\rm {and}} \,
i_B<j_D<j_B} \vert  D \vert \in \pi.
 $$  Set  $$e_w(A, B)=2 (A\circ_w B -  B\circ_w A)+ n_w(A,B)\,( \vert A\vert+ \vert B\vert) \in \pi, $$
$$ e_w(A,s)=-e_w(s,A)=
\sum_{ D \in \A }  {n_w (A,D)} \, \vert D\vert \in \pi $$ and
$e_w(s,s)=0$. All  these expressions are sums over   $D\in \A$ of
terms $m_D \vert D\vert$ with $m_D\in \ZZ$. The term $m_D \vert
D\vert$ is   the {\it contribution} of   $D$. For instance, the
contributions of $A,B$ to $e_w(A, B)$ are $ n_w(A,B)\,  \vert
A\vert$ and $n_w(A,B)\, \vert B\vert$, respectively. The
contribution of $A$ to $e_w(A,s)$ is 0. It is clear that a letter
$D\neq A, B$ may contribute non-trivially to $e_w(A,B) $ only if
(i) each entry of $D$ appears between the two entries of $A$ or
between the two entries of $B$ (or both) and (ii) $D$ occurs at
least once between the   entries of $A$ and at least once between
the   entries of $B$.

  It is useful  to have more direct formulas for $e_w(A, B)$. For   words $x,y$ in
the alphabet $\A$, set  $\langle x, y\rangle=\sum_{D } \vert
D\vert\in \pi$  where $D$ runs over the letters in $\A$ occurring
exactly  once in $x$ and exactly once in $y$. If $w=\cdots A x A y
B z B \cdots $, where $x,y,z $ are words in the alphabet $\A$,
then \begin{equation}\label{mn1}e_w(A, B) =2 \langle x, z\rangle
\end{equation}  If $w=\cdots  A x B y B z A \cdots$, then
\begin{equation}\label{mn}e_w(A, B)= 2 \langle x, y\rangle - 2
\langle y, z\rangle .\end{equation}  If $w=\cdots A x B y A z B
\cdots $, then \begin{equation}\label{x56}e_w(A, B)= 2   \langle
x, y\rangle  + 2
  \langle x, z\rangle + 2
\langle y, z\rangle +\vert A\vert +\vert B\vert.\end{equation} If
$B$ occurs in $w$ before both entries of $A$, then   $e_w(A,
B)=-e_w(B, A)$ can be computed applying   these   formulas to the
pair $(B,A)$.
 Similarly, if $w=xAyAz$,
then $e_w (A,s)= \langle y, z\rangle- \langle x, y\rangle$.

It is
easy to describe the behavior of the  pairings associated with
nanowords under pushing forward. Consider  a set with involution
$(\overline \alpha, \overline \tau)$ and an equivariant mapping
 $f:\overline \alpha\to \alpha$.   This mapping induces an additive homomorphism
 $\overline \pi=\pi(\overline \alpha,\overline \tau)\to  \pi(\alpha,\tau)=\pi$  denoted   $f_{\#}$.
   Consider  a
nanoword $\overline w$ over  $\overline \alpha$ and let $w$ be the
nanoword over $\alpha$ obtained from $\overline w$ by pushing
forward along $f$. Then the $\alpha$-pairing $(S,s,e_{w}:S\times S
\to \pi)$ associated with $w$ is
 obtained from the $\overline \alpha$-pairing
 $(S,s,e_{\overline w}:S\times S \to \overline \pi)$ associated with $\overline w$ by
 composing $e_{\overline
 w}$ with $f_{\#}$ (the set $S$ is preserved). Thus, $e_w= f_{\#} \circ
 e_{\overline w}$.

As an exercise, the reader may check that $p(w^-) = (p(w))^-$ and
$p(w_1 w_2)= p(w_1)\oplus p(w_2)$ for any nanowords $w, w_1, w_2$ over $\alpha$.

\subsection{Homomorphism
 $p:{\mathcal{N}}_c \to {\mathcal{P}_{sk}}$.}\label{fihomomlb99}
 Let ${\mathcal{N}}_c={\mathcal{N}}_c (\alpha, \tau)$ and ${\mathcal{P}_{sk}}={\mathcal{P}}_{sk}(\alpha, \tau, \pi)$ be the groups introduced in Sections \ref{fi2grouppsks9} and \ref{ermlp99}, respectively, 
 where $\pi$ is the group defined in Section  \ref{fipipipiib99}.

\begin{theor}\label{lemplo} The formula $w\mapsto p(w)$
 defines a  group
 homomorphism
 $p:{\mathcal{N}}_c \to
 {\mathcal{P}_{sk}} $. \end{theor}

 The homomorphism  $p:{\mathcal{N}}_c\to {\mathcal{P}_{sk}}$  is in general not injective.
 This is clear already from the fact that  ${\mathcal{N}}_c$ may be non-commutative while $ {\mathcal{P}_{sk}}$
 is commutative.

We now reduce Theorem \ref{lemplo} to a   lemma which will be   proven in Sect.\
\ref{fi:mxlxlxlxlx40} using topological techniques.  

\subsection{Proof  of Theorem \ref{lemplo} (modulo a lemma)}\label{fsymnano9}
The multiplicativity of $p:{\mathcal{N}}_c\to {\mathcal{P}_{sk}}$
follows from the definitions. We need to check only that $p$ is
well-defined, i.e., that cobordant nanowords give rise to
cobordant $\alpha$-pairings. A direct comparison shows that  two
nanowords related by the third homotopy move have isomorphic
$\alpha$-pairings (cf.\ \cite{tu3}, proof of Lemma 7.6.1). Since
the first and second homotopy moves are special instances of
surgery, it remains to show that nanowords related by a surgery
have cobordant $\alpha$-pairings.

We begin by fixing notation. Consider a nanoword $(\A, w)$ and its
even symmetric factor $\nabla=(\B, (v_1\,\vert \, \cdots \, \vert
\, v_k))$. Thus, $\B\subset \A$  and $w=x_1 v_1x_2v_2   \cdots
x_{k} v_k x_{k+1}$ where $x_1,x_2,..., x_{k+1}$ are words in the
$\alpha$-alphabet $\C=\A-\B$. Deleting $\nabla$ we obtain the
nanoword $(\C, x=x_1 x_2 \cdots x_{k+1})$. Let
$\iota=\iota_\nabla:\B \to \B$ and
$\varepsilon=\varepsilon_\nabla:\B \to \{0,1\}$ be the involution
and the mapping associated with $\nabla$ in Sect.\ \ref{nanop}.

 We must prove that the $\alpha$-pairings $p(w) $ and
$p(x) $ are cobordant.
 Replacing each letter $C\in \C$ by its copy $C'$,  we obtain  a nanoword
 $ (\C'=\{C'\}_{C\in \C}, x')$  isomorphic to  $(\C,x)$.
 It suffices to  verify that the  $\alpha$-pairing  $ p(w)\oplus  p( x')^-$
  has an annihilating filling.
 Let $p=p(w)=(S,s,e_w)$ and $p'=p(x')= (S',s', e_{x'})$ where
 $S=\A \cup \{s\}$ and $S'=\C'\cup \{s'\}$.
 With  every   $C\in \C$ we associate the vector
 $\lambda_C=C+ C'\in  \ZZ S\oplus \ZZ  S'$. With  every    $B\in \B$ we associate
 the vector
 $\lambda_B\in \ZZ S \subset  \ZZ S\oplus \ZZ  S'$ equal to $B$ if $B=B^\iota$ and equal to
 $B+ (-1)^{\varepsilon (B)} B^\iota$  if $B\neq B^\iota$.
 (Note that $\lambda_{B^\iota}= (-1)^{\varepsilon (B)} \lambda_B$.)
 To proceed,   pick a set $\B_+\subset \B$ meeting every orbit of
 $\iota:\B\to \B$ in one  element.
The  set of vectors $\{\lambda_C\}_{ C\in \C} \cup\{\lambda_B\}_{ B\in
\B_+}\cup \{   s+s'\}$ is  a filling of $ p\oplus (p')^-$.   We
claim that this filling is annihilating. Indeed, for $C\in \C$,
$$ (e_w\oplus e_{x'}^-) ( \lambda_{C}, s+s')= e_w(C, s) +
(e_{x'})^- (C', s')= e_w(C, s) -e_x(C,s)=0.$$ The latter equality
follows from two facts:  for   $D\in \B$ with $D=D^\iota$, we have
$n_w(D,C)=0$ so that   $D$ contributes 0 to $e_w(C, s)$; for $D\in
\B$ with  $D\neq D^\iota$, we have either
$n_w(D,C)=n_w(D^\iota,C)=0$ or $n_w(D,C)=n_w(D^\iota,C)=\pm 1$. In
the latter case $\varepsilon (D)=1$ and   $\vert D^\iota
\vert=\tau(\vert D \vert)$. In all  cases, the pair $D, D^\iota$
contributes $0$ to $e_w(C, s)$. Hence the total contribution of
all $D\in \B$ to $e_w(C,s)$ is $0$. Therefore $e_w(C, s)
=e_x(C,s)$.

 Similarly,  for  $C_1, C_2 \in \C$,
 $$ (e_w\oplus e_{x'}^-) (\lambda_{C_1},
 \lambda_{C_2})=e_w(C_1,C_2)-e_x(C_1,C_2)=0.$$ Indeed, if
 $D\in \B$ contributes non-trivially to $e_w(C_1,C_2)$ then
    $\varepsilon (D)=1$,   $\vert D^\iota
  \vert =\tau(\vert D \vert)$ and the sum of the contributions of $D,
  D^\iota$ to $e_w(C_1,C_2)$ is $0$. Here we use the evenness of
  $\nabla$ which implies that
  $\varepsilon (D)=1\Rightarrow D\neq D^\iota$.

 It remains  to prove that for any $B, B_1, B_2\in \B, C\in \C$,
$$ (e_w\oplus e_{x'}^-) (\lambda_{B_1}, \lambda_{B_2})= (e_w\oplus e_{x'}^-) (\lambda_{B}, \lambda_C)=
 (e_w\oplus e_{x'}^-) (\lambda_{B}, s+s')=0.$$
Since $\lambda_{B}, \lambda_{B_1}, \lambda_{B_2} \in \ZZ \A$,
these formulas are equivalent to
\begin{equation}\label{plo} e_w (\lambda_{B_1},
\lambda_{B_2})=e_w (\lambda_{B}, C)=
 e_w (\lambda_{B}, s)=0.\end{equation}

In the rest of the proof, we denote by $\alpha_0$ the 2-letter alphabet $\{+,-\}$ with involution $\tau_0$ permuting $+$ and $-$.

 \begin{lemma}\label{l:edlocfc12}
Formula \ref{plo} holds for     $\alpha=\alpha_0$ and  $
\tau=\tau_0$.  \end{lemma}

Lemma \ref{l:edlocfc12} will be proven in Sect.\ \ref{Words and loops} using topological techniques.

\begin{lemma}\label{l:edlocfzzzzpc12}
  Let
$e\in \pi =\pi(\alpha,\tau)$ be one of the expressions $e_w
(\lambda_{B_1}, \lambda_{B_2})$,  $e_w (\lambda_{B}, C)$,
 $e_w (\lambda_{B}, s)$  in Formula \ref{plo}, 
 where $\alpha$ is an arbitrary alphabet with involution $\tau$.
For any additive homomorphism $\varphi:\pi \to \ZZ$ sending all
the generators  $a\in \alpha$ to $\{-1,+1\}\subset \ZZ$, we have
$\varphi(e)=0$.
\end{lemma}
                     \begin{proof}
                     The homomorphism $\varphi$ induces a mapping
                      $f:\alpha\to \alpha_0=\{-1,+1\}$
                     sending each $a\in \alpha$ to $\varphi(a) \in \alpha_0$.
                      The additivity of $\varphi$ implies that $f$ is
                      equivariant with respect to the involutions $\tau:\alpha \to
                      \alpha$ and $\tau_0:\alpha_0 \to \alpha_0$.
                     Let $w_0$ be the nanoword over
                     $\alpha_0$ obtained    by pushing
                      $w$ forward along $f$. The
                     $\alpha_0$-pairing $(S,s,e_{w_0})$ of $w_0$ is obtained from
                     the $\alpha$-pairing $(S,s,e_w)$ of $w$ by
                     composing $e_w:S\times S \to \pi$ with
                     $\varphi$.   Thus
$\varphi (e_w (\lambda_{B_1}, \lambda_{B_2}))=e_{w_0}
(\lambda_{B_1}, \lambda_{B_2})$,
 $\varphi(e_w (\lambda_{B}, C))=e_{w_0}
(\lambda_{B}, C)$,  and $ \varphi(
 e_w (\lambda_{B}, s))=e_{w_0}
(\lambda_{B}, s)$.  Lemma \ref{l:edlocfc12} implies that the
right-hand sides of these formulas are equal to  $0$.
\end{proof}

\begin{lemma}\label{2304040edlocfzzzzpc12}
 Formula \ref{plo} holds for any alphabet
 $\alpha$ with fixed-point-free involution $\tau$.
\end{lemma}
                     \begin{proof} If $\tau$ is fixed-point-free,
                     then   $\pi=\pi(\alpha,\tau)$ is a
                     free abelian group with basis $\beta$
                     where $\beta$ is any subset of
                     $\alpha$ meeting every orbit of $\tau$ in one element.
                     Let
$e\in \pi $ be one of the expressions $e_w (\lambda_{B_1},
\lambda_{B_2})$,  $e_w (\lambda_{B}, C)$,
 $e_w (\lambda_{B}, s)$  in Formula \ref{plo}. We
 expand $e=\sum_{x\in \beta} k_x \, x$ with  $k_x\in \ZZ$. We claim that $k_x=0$ for
 all $x$. Indeed, pick any $x \in \beta$.
 Consider the  additive homomorphisms $\varphi_+, \varphi_-:\pi\to \ZZ$ defined by  $\varphi_\pm (y) =1$ for all
  $y\in \beta -\{ x\}$ and $\varphi_\pm (x)=\pm 1$. By
the previous lemma, $\varphi_\pm (e) =0$. Therefore $2
k_{x}=\varphi_+ (e)-\varphi_-(e)=0$. Hence $k_{x}=0$.
\end{proof}

We can  now
prove Formula \ref{plo} in its full generality. We begin by
associating with the nanoword $(\A,w)$ another nanoword as
follows. For $i=0,1$, set $\B_i=\{B\in \B\,\vert \,
\varepsilon_\nabla (B)=i\}$.   It follows from the definition of
the  involution $\iota=\iota_\nabla$ on $\B$ that
$\iota(\B_i)=\B_i$.    Set $X=\B_0/\iota$ and let $X'=\{x'\, \vert
\, x\in X\}$ be a copy of $X$. Similarly, let $\C'=\{C'\, \vert \,
C\in \C\}$ be a copy of $\C=\A-\B$. Set $\overline \alpha=X \cup
\B_1 \cup \C \cup X' \cup \C' $ where it is understood that the
five sets on the right-hand side are disjoint. There is a unique
involution $\overline \tau $ on $ \overline \alpha$ such that
$\overline \tau (x)=x'$ for $x\in X$, $\overline \tau (C)=C'$ for
$C\in \C$ and $\overline \tau (B)=B^\iota$ for $B\in \B_1$. The
projection $\B_0\to X$ and the identity on $\B_1 \cup \C$ form
  a mapping $ p:\A =\B_0\cup \B_1 \cup \C \to \overline \alpha$.
This mapping  makes $\A$ into an $\overline \alpha$-alphabet. The
pair $(\A,w)$ becomes thus a nanoword over $\overline \alpha$. We
denote this nanoword by $\overline w$.

The nanoword $\overline w$ and the original nanoword $(\A, w)$
over $\alpha$  coincide as words in the alphabet $ \A$ and differ
in the choice of the ground alphabet. The nanoword  $w$
 is a push-forward of $\overline w$ as follows.
Define a mapping $f:\overline \alpha \to \alpha$ by $f(x)=\vert
B\vert, f(x')=\tau (\vert B\vert)$ for any $x\in X, B\in
p^{-1}(x)\subset \B_0$ and $ f(B)=  \vert B\vert , f(C)= \vert
C\vert,  f(C')= \tau (\vert C\vert)$ for   $B\in \B_1, C\in \C$.
The mapping $f$ is well-defined on $X\subset \overline \alpha $
because $\vert B^\iota\vert =\vert B\vert$ for $B\in \B_0$. The
definition of $\overline \tau$  and the equality $\vert
B^\iota\vert =\tau (\vert B\vert)$ for $B\in \B_1$ imply that $f$
is equivariant with respect to the involutions $\overline
\tau:\overline \alpha \to \overline \alpha$ and $\tau:\alpha\to
\alpha$. Since the composition of $f$ with $p:\A \to \overline
\alpha$ is the given projection $\A\to \alpha, A\mapsto \vert
A\vert$, the nanoword $w$ is the push-forward of $\overline w$
along $f$.

 The
 phrase $(\B,\nabla)$ with projection $p \vert_{\B}:\B \to \overline \alpha$ is
 a nanophase over $\overline \alpha$ denoted $\overline \nabla$. Clearly,
 $\overline \nabla$ a factor of $\overline w$ since this
property does not involve the   ground alphabet.  The nanophrase
$\overline \nabla$ is   even because so is $\nabla$. Obviously,
$\iota_{\overline \nabla}  = \iota_{ \nabla} =\iota:\B\to \B$ and
$\varepsilon_{\overline \nabla}  = \varepsilon_{ \nabla} :\B \to
\{0,1\}$. Our   definition of $\overline \alpha$ ensures   that
  $\overline \nabla$ is symmetric. Indeed, for any
$B\in \B$, if $\varepsilon_{ \nabla} (B)=0$, then $B\in \B_0$ and
$p(B^\iota)=p(B)\in X$ by the definition of $X,p$. If
$\varepsilon_{ \nabla} (B)=1$, then $B\in \B_1$ and
$p(B^\iota)=B^\iota =\overline \tau (B)= \overline \tau ( p(B))$.

Observe that the involution $\overline \tau$ is fixed-point-free.
This follows from the fact that  $\iota=\iota_\nabla$ has no fixed
points in $\B_1$  which in its turn follows from the evenness of
$\nabla$. By Lemma \ref{2304040edlocfzzzzpc12}, $$e_{\overline w}
(\lambda_{B_1}, \lambda_{B_2})= e_{\overline w} (\lambda_{B}, C) =
 e_{\overline w} (\lambda_{B}, s)=0.$$  Since $w$ is a push-forward of
$\overline w$,  the latter formula implies Formula \ref{plo}.

\subsection{Examples and remarks}\label{remlexamps799} 1. Pick three (possibly coinciding) elements    $a,b, c\in \alpha$.
Consider the nanoword $w =ABCBAC$ with $\vert A\vert =a, \vert B
\vert =b, \vert C\vert =c$. If $a= \tau (b)$, then
 $w $ is contractible and therefore slice. We use Theorem \ref{lemplo} to
 verify   that $w$  is not slice for
$a\neq \tau(b)$. It suffices to verify that the $\alpha$-pairing $
p(w)= (S=\{s,A,B,C\},s,e_w)$ is not hyperbolic.  A direct
computation from definitions shows that $e_w$ is given by the
matrix $$  \left [
\begin{array}{ccccc}
0& -c    & -c & a+b  \\
      c& 0    & 0 & a+2b+c  \\
c& 0    & 0 & b+c   \\
 -a-b& -a-2b-c    & -b-c & 0  \\
\end{array} \right ] $$
where the rows and columns correspond to $s, A, B, C$,
respectively. An easy check shows that   $p(w)$ has no
annihilating fillings. Indeed, the tautological filling formed by
the  vectors $s, A, B, C$ is non-annihilating since $c\neq 0$ in
$\pi$. The vectors $A-B$ and $ B-A$ cannot belong to a filling since $a
\neq \tau (b)$. A vector  $ \lambda$ of the form $A \pm C, B\pm C,
C - A, C-B $ cannot belong to an annihilating filling since
$e_w(\lambda, s)$ is an algebraic sum of $a,b,c$ and is non-zero
in $\pi$.   It remains to consider the family of vectors $\{s,
A+B, C\}$ which is a filling if     $a=b$. We have   $e_w(A+B, s)
=2c$ and $ e_w(A+B, C)=a+2c+3b=4a+2c$. If this filling is
annihilating, then $2c=4a+2c=0$ in $\pi$ and therefore $4a=0$.
This may happen only when $a=\tau(a)$ which is excluded by   $a=b,
a\neq \tau (b)$. Thus, for $a\neq \tau (b)$, the pairing $p(w)$ is
not hyperbolic and the nanoword $w  $ is not slice. Note that if
$a,b,c$ belong to one orbit of $\tau$, then $\gamma (w)=1$. This
shows that the
  $\alpha$-pairings may provide more information than the
homomorphism $\gamma$.

2. Consider the nanoword $w=ABCADCBD$ with $\vert A\vert =\tau
(\vert B\vert)$, $ \vert C \vert =\vert D\vert$. Direct
computations  show  that $\gamma (w)=1$ and $p(w)$ is hyperbolic
with annihilating filling $\{s, A-B, C+D\}$. (A general
construction producing such examples will be discussed  in Sect.\
\ref{ushifts21}). The author does not know whether $w$  is slice
except in the case where $\vert A\vert =\vert D\vert$ (then $w$ is
symmetric and therefore slice).

3.  Given   a nanoword $(\A , w)$ over $ \alpha $ and a family
$H=\{H_a \}_{a\in \alpha}$ of subgroups of $\pi$ such that
$H_a=H_{\tau(a)}$ for all $a$, we define    the {\it $H$-covering}
of $w$
 to be the nanoword  $(\A^H,w^H)$ over
$\alpha$ obtained  by deleting   from both $\A$ and $w$ all
letters $A$ such that $e_w(A,s) \notin H_{\vert A\vert}$ (cf.\
\cite{tu3}). One  may check that the
$H$-coverings of cobordant nanowords are cobordant (we shall not
use it). Moreover, the formula $w\mapsto w^H$ defines  a group
endomorphism of $ \mathcal N_c(\alpha, \tau)$.

\section{Cobordism invariants of nanowords}

    \subsection{The $u$-polynomial}\label{fipolynomial99}
   The   {\it  $u$-polynomial} of a nanoword $w$ is defined by
    $u^w= u^{p(w)}$ where $p(w)$ is the $\alpha$-pairing   associated with $w$ and
      $u^{p(w)}$ is its
     $u$-polynomial.
     By Theorems \ref{coppprvplo}  and   \ref{lemplo},
$u^w$  is a cobordism invariant of $w$.

    Consider in more detail the case where  $\tau $ is fixed-point-free.  Let $k$ be the number of orbits of $\tau$ and $t_1,t_2,...,t_k\in  \alpha$ be representatives of the  orbits  so that each orbit contains exactly one $t_i$. It is convenient to switch to the multiplicative notation for the group operation in $\pi=\pi(\alpha, \tau)$. Thus, any
    $g\in  \pi$     expands uniquely
as a monomial  $g=\prod_{i=1}^k  t_i^{m_i}    $
with   $m_1,..., m_k \in \ZZ$. Identifying $\delta_g\in I$
 with this monomial, we  identify the  $\ZZ$-module  $I$   from   Sect.\ \ref{eruuuuusp99}  with the additive group of    Laurent polynomials over  $\ZZ$ in the commuting variables
$t_1, ...,t_k$ with zero free term.  The $\ZZ$-submodule $J\subset
I$ consists of those Laurent polynomials  with zero free term
which are invariant under the  inversion of the
variables $t_1\mapsto t_1^{-1},...,t_k\mapsto t_k^{-1} $.   The
quotient   $I/J$ is an infinitely generated  free abelian group
with basis    $ \prod_{i=1}^k  t_i^{m_i}   (\modu J) $ where the
tuple  $(m_1,..., m_k )$ runs over    $k$-tuples of integers such
that at least one of its entries is non-zero and the first
non-zero entry   is positive.  The {\it degree} of such a basis
monomial is the number $\sum_{i=1}^k \vert m_i\vert \geq 1$. For
$x\in I/J$, we define its  {\it degree}  $\deg (x) $ to be the
maximal degree of a basis monomial appearing in $x$ with non-zero
coefficient. The number $\deg(x)$ does not depend on the choice of
$t_1,..., t_k$. It can be used to estimate the length norm of  nanowords. It follows from the definitions that for any nanoword $w$ and any  $a\in \alpha$,
\begin{equation}\label{pop111}\vert \vert w\vert \vert_l \geq \deg (u^w(a)) +1.
\end{equation}

\begin{theor} \label{8fgpm59s}  If $\tau \neq \id$, then the group $ \mathcal{N}_c= \mathcal{N}_c(\alpha, \tau)$ is infinitely generated.
\end{theor}
\begin{proof}  Let $\beta \subset  \alpha$ be a free orbit of $\tau$.  Since the pull-back homomorphism
$\mathcal N_c(\alpha, \tau) \to \mathcal N_c(\beta, \tau\vert_\beta)$ is surjective, it suffices to prove  the theorem in the case where $\alpha$ consists of two  elements permuted by $\tau$.
We give a more general argument working for all fixed-point-free $\tau$.

Pick   representatives $t_1,t_2,...,t_k\in  \alpha$ of the  orbits of $\tau$ as above.  For   $j=1,...,k$,   define an additive homomorphism $\partial_j:I \to \ZZ$ by
$\partial_j (\prod_{i=1}^k  t_i^{m_i} )=m_j$. Clearly, $ \partial_j (J)=0$. The  induced  group homomorphism $I/J \to \ZZ$ is also denoted $\partial_j$.

Given a nanoword $w$, the function $u^{w}:\alpha\to I/J$ satisfies
$u^{w} (a)= - u^{w}(\tau (a))$ and   is therefore determined by
its values on $t_1,..., t_k$. Theorem 6.3.1 of \cite{tu3} implies
that a sequence  $f_1,...,f_k\in I/J$ is realizable as a sequence
 $u^{w} (t_1), ..., u^{w}(t_k)$ for a  nanoword $w$ if and only if
 $    \partial_j (f_i)+ \partial_i  (f_j)=0$ for all $i, j$. This gives $k(k+1)/2$ conditions
 so that realizable sequences form a subgroup of $(I/J)^k$ of corank $\leq k(k+1)/2$.
 Since $I/J$ is an infinitely generated  free abelian group,
the  image of the  homomorphism $ \mathcal{N}_c \to (I/J)^k, w \mapsto (u^{w}(t_1),...,u^w(t_k))$ is  infinitely generated. Therefore  $ \mathcal{N}_c$ is infinitely generated. \end{proof}

\subsection{The genus}\label{fcoinnan68}
 Any    commutative  domain  $F$ is a $\ZZ$-algebra in the
 usual way.  Given   an additive homomorphism
$\varphi: \pi \to F$, the $\varphi$-{\it genus} of a  nanoword $w$
is defined by $\sigma_\varphi (w)=\sigma_\varphi (p (w)) \in \ZZ$.
 By
Theorems \ref{859s} and \ref{lemplo},  this
 number  is a cobordism invariant of $w$. For any nanowords
  $w,w_1,w_2$ we have $\sigma_\varphi
(w)=\sigma_\varphi (w^-)$ and $\sigma_\varphi (w_1w_2) \leq
\sigma_\varphi (w_1) + \sigma_\varphi (w_2)$. The $\varphi$-genus
can be used to estimate  the length  norm from below: for any
non-slice $w$,
\begin{equation}\label{pop} \vert \vert w\vert
 \vert_l \geq \sigma_\varphi (w )/2 +1,
\end{equation}
see  Sect.\ \ref{brieeeeeword}. If $w$ is slice, then $\vert \vert
w \vert \vert_l = \sigma_\varphi (w )=0$.

\subsection{Example}\label{reeeexcs799}   Consider  the nanoword
$w =ABCBAC$ from Example \ref{remlexamps799}.1.  Suppose    that
$a=\vert A\vert,b=\vert B\vert, c=\vert C\vert \in \alpha$ are not
fixed points of $\tau$. We show how to use the $u$-polynomial
$u^w$ to compute $ \vert \vert w\vert \vert_l$. Pushing back, if
necessary,  to the union of the orbits of $a,b,c$ we can assume
that $\tau$ is fixed-point-free.
 We have   $[c]_{p(w)} =\delta_{-a-b}+r \delta_c $ where
  $r\in \ZZ$ is zero unless $c=a$ and/or $c=b$.
Since $a\neq \tau (b)$, the    monomials $\delta_{-a-b}, \delta_c$
have degrees $2$ and $1$, respectively. Therefore  $\deg u^w(c)
=\deg u^{p(w)} (c)=2$ and by Formula \ref{pop111}, we have $ \vert
\vert w\vert \vert_l \geq 3$. Since $w$ is a nanoword of length 6,
we have  $ \vert \vert w\vert \vert_l=3$.

Assume additionally that $a,b,c$ belong to different orbits of
$\tau$. Then the $\alpha$-pairing $p(w)$ has only one filling
$\lambda$, the tautological one. A direct computation shows that
for any homomorphism  $\varphi$ from $\pi$ to a commutative domain   such that
$\varphi(a+b)\neq 0$ and $\varphi (c)\neq 0$, we have
$\sigma_\varphi (p(w))=\sigma_\varphi (\lambda)=2$. Formula
\ref{pop} gives in this case $ \vert \vert w\vert \vert_l \geq 2$
which is weaker   than Formula \ref{pop111}.

  \section{Bridges and the  bridge  norm}\label{avgpm1za}

\subsection{Bridges}\label{brinword}
A {\it quasi-bridge} in a nanoword  $(\A,w)$ is a pair consisting of a
factor
  $\nabla=(\B, (v_1\,\vert \, \cdots \, \vert \, v_k))$  of $w$ with $k\geq 1$ and an involutive
    permutation $\kappa:\widehat  k \to \widehat  k$ of the set $\widehat  k=   \{1,2,..., k\}$
      satisfying the
     following two conditions:

(a)  the length of $v_r$ is even for any $r\in \widehat  k$
such that $\kappa (r)=r$;
 
 (b)    there is a  mapping  $\iota:\B \to \B$ such that   $\iota v_r= v_{\kappa(r)}^-$ for all $r\in \widehat  k$.

Consider  a quasi-bridge $(\nabla, \kappa)$. Let $n_r=n_{\kappa(r)}$ be  the length of $v_r$ for $r=1, \ldots, k$. Any entry
of a letter $B\in \B$   in $w$  appears in some $v_r$, say,
on the $j$-th position where $  1\leq j\leq n_r$. By
  (b), the letter $B^\iota=\iota(B)$ appears   in
$v_{\kappa(r)}$ on the $(n_r+1-j)$-th position. The latter entry of
$B^\iota$ in $w$ is said to be {\it symmetric} to the original
entry of $B$ in $w$. Thus, the mapping $\iota=\iota_{\nabla, \kappa}$   
is uniquely determined by  $(\nabla, \kappa)$ and  $\iota^2=\id$.

For $B\in \B$, set $\varepsilon_{\nabla, \kappa}(B)=1$ if the entry symmetric to the  leftmost  entry of $B$  in $w$ is the leftmost  entry of $B^\iota$ in $w$. Otherwise, set 
$\varepsilon_{\nabla, \kappa}(B)=0$. Clearly,
$\varepsilon_{\nabla, \kappa}(B)=\varepsilon_{\nabla,
\kappa}(B^\iota)$.  The quasi-bridge $(\nabla, \kappa)$ is  a  {\it bridge} if for all $B\in \B$,
$$\vert B^\iota \vert =\tau^{\varepsilon_{\nabla, \kappa}(B)}
(\vert B\vert).$$

 Given  a bridge $(\nabla=(\B, (v_1\,\vert \, \cdots \, \vert \, v_k)), \kappa)$  in a nanoword  $(\A,w)$,
we can delete all letters of the set  $\B\subset \A $ from $\A$
and $w$. For $w=x_1 v_1x_2v_2   \cdots  x_{k} v_k x_{k+1}$, the
deletion yields the nanoword $(\C=\A-\B, x_1 x_2... x_{k+1})$.
This transformation $(\A,w) \mapsto (\C, x_1 x_2... x_{k+1})$ is
the   {\it bridge move} determined by $(\nabla, \kappa)$. A bridge
move   is always associated with a specific  bridge. Thus,
two  bridges $(\nabla, \kappa)$ and $(\nabla', \kappa')$
in $w$ determine the same move if and only if $\nabla=\nabla'$ and
$\kappa=\kappa'$.

For a bridge move $m$ determined by  a bridge $(\nabla,
\kappa)$, the  free (2-element) orbits of the involution
$\kappa:\widehat  k \to \widehat  k$, where $k$ is the length of $\nabla$,
are called {\it arches} of $m$. The number of arches of $m$ is
denoted $\ar(m)$. Obviously, $\ar(m)\leq [k/2]$.  For the inverse
move $m^{-1}$, set $\ar(m^{-1})=\ar (m)$.

 For $\kappa=\id$,  a pair ($\nabla$,
 $\kappa$) as above is  a bridge if and only if $\nabla$ is an even
   symmetric factor. This follows from the fact that in this case
   $\varepsilon_{\nabla, \kappa}=\varepsilon_{\nabla}$ is the
   function introduced in Sect.\ \ref{nanop}.
 Therefore     bridges generalize  even symmetric
  factors. The latter are precisely  the bridges   with 0 arches.
  Surgeries are precisely the bridge moves determined by  bridges with
0 arches.

 \subsection{Examples.}\label{brinexaword} 1. Given a nanoword $(\A, w)$ and a
 letter $A\in \A$ we can split $w$ uniquely as $x_1Ax_2Ax_3$ where $x_1, x_2, x_3$
  are words in the alphabet $\A-\{A\}$.
  The factor $(A\,\vert \,A)$ of $w$ endowed with transposition
   $\kappa= (12)$ and the identity mapping $\iota: \{A\}\to \{A\}$ is
   a bridge
    with one  arch
    (here  $\varepsilon_{\nabla, \kappa}(A)=0$). Deleting this bridge,
    we obtain the nanoword $(\A-\{A\}, x_1 x_2 x_3)$.

2.  Consider  a nanoword $(\A, w =x_1ABx_2 BA x_3)$ where $A,B\in
\A$. The factor $(AB\,\vert \,BA)$ of $w$ endowed with    $\kappa=
(12)$ is  a bridge. Its deletion gives  the nanoword
$(\A-\{A,B\}, x_1 x_2 x_3)$.

3.  Consider  a nanoword $  ABCBCDAD $
 with  $\vert A \vert, \vert B \vert= \vert C\vert, \vert D \vert \in \alpha$.
  The factor $(A\,\vert \,  BCBC  \,\vert \,A)$ of $w$  with
    $\kappa= (13)$
 is    a bridge.  Its deletion gives  $DD$.

4.  Consider  a nanoword $  ADEBC DCAEB$  with
 $\vert A \vert= \vert B \vert, \vert C\vert, \vert D \vert,
  \vert E \vert\in \alpha$. The factor $(A\,\vert \,  BC \,\vert \,
    CA  \,\vert \,B)$ of $w$  with    $\kappa= (14)(23)$ is     a bridge.
     Its deletion gives  $DEDE$.

\subsection{The bridge norm}\label{brinmetricword}  The list (TR) of transformations on nanowords considered  in Sect.\ \ref{fi2grouppsks9}  can   be extended to the following wider list:

(TR+) isomorphisms,   homotopy moves,   bridge moves, and  the inverse
moves.

Given two nanowords $v,w$,  {\it a metamorphosis} $m:v \to w$ is a
finite sequence   $m=(m_1,...,m_n)$ of   moves from    the list
(TR+)  transforming $v$ into $w$.   The inverse metamorphosis
$m^{-1}=(m_n^{-1},..., m_1^{-1})$ transforms $w$ into $v$. Set $\ar
(m)=\ar(m_1)+\cdots + \ar(m_n)$ where $\ar(m_i)$ is the number of
arches of $m_i$ if $m_i$ is a bridge move or its inverse and
$\ar(m_i)=0$ for all other moves.  Clearly,  $\ar(m^{-1})=\ar(m)$
and
  $\ar(m)=0$ if  and only if all the bridge moves in   $m$ are surgeries
  or inverses of surgeries.

For any   nanoword    $w$, there is a metamorphosis   $w\to \emptyset$. For instance,  one can  consecutively delete the   letters of $w$ as in Example \ref{brinexaword}.1.  Set
$$\vert \vert w \vert \vert_{br}=\min_m \ar (m) \geq 0$$
where $m$ runs over all metamorphoses   $w\to \emptyset$.

\begin{lemma}\label{th:v54937892153}  The function  $w\mapsto \vert \vert w \vert \vert_{br}$  induces a conjugation  invariant $\ZZ$-valued   norm on $\N_c=
\N_c(\alpha, \tau)$.
\end{lemma}
\begin{proof} If $w\sim_c v$, then there is a metamorphosis $m:w\to v$ with $\ar (m)=0$.
Composing $m$ with a metamorphosis $M:v\to \emptyset$  we obtain a metamorphosis $M':w\to \emptyset$ with $\ar(M')=\ar (M)$. Therefore $\vert \vert w \vert \vert_{br} \leq \vert \vert v \vert \vert_{br}$.
By symmetry, $\vert \vert w \vert \vert_{br}=\vert \vert v \vert \vert_{br}$.
Therefore  the formula
  $w\mapsto \vert \vert w \vert \vert_{br}$  defines    a  function
   $\N_c \to \ZZ$. That the latter   satisfies all axioms of a $\ZZ$-valued norm
 directly follows from the definitions.  To show that it is  invariant under conjugation, it is enough to show that $\vert \vert v w v^- \vert \vert_{br}\leq  \vert \vert w \vert \vert_{br}$ for all nanowords $v,w$.
Any metamorphosis $m:w\to \emptyset$, extends by the identity on $v, v^-$ to a metamorphosis $m':vwv^- \to v\emptyset v^-=vv^-$ with $\ar(m')=\ar(m)$.
The  symmetric  nanoword $vv^-$  can be transformed into $\emptyset$ by a single surgery.
Therefore $\vert \vert  v w v^-  \vert \vert_{br}\leq  \vert \vert w \vert \vert_{br}$.
\end{proof}

 The   $\ZZ$-valued   norm on $\N_c$ provided by this lemma is denoted $\vert \vert  \cdot  \vert \vert_{br}$   and  called the {\it bridge norm}.
A   consecutive deletion of all but one letters of a nanoword $w$
yields a metamorphosis of $w$ into a contractible nanoword of type
$AA$. Therefore if $w$ is non-slice, then
\begin{equation}\label{zaz}\vert \vert w \vert \vert_l \geq
\vert \vert w \vert \vert_{br} + 1.\end{equation}

\subsection{The bridge metric}\label{brieeeeeword}  The bridge norm induces  a  left- and right-invariant metric $\rho_{br}$  on  $\mathcal N_c$ by $\rho_{br} (w,v)= \vert \vert wv^- \vert \vert_{br}$ for any nanowords $w,v$, cf.\  Sect.\ \ref{fi2normks9}.  Formula (\ref{zaz})
 implies  that
$\rho_l (w,v) \geq \rho_{br} (w,v)+1$.

\begin{lemma}\label{th:v54937ccccc153}  For any   nanowords $v,w$,
$$\rho_{br} (w,v)=\min_m \ar (m)  $$
where $m$ runs over all metamorphoses   $w\to v$.
\end{lemma}
\begin{proof} Denote the right-hand side by $\eta  $. Given a metamorphosis $m:w\to v$ we can extend it by the identity on $v^-$ to  a metamorphosis $m':wv^-\to vv^-$. Composing the latter with the surgery $vv^-\to \emptyset$, we obtain a  metamorphosis $m'':wv^-\to  \emptyset$
 with $\ar (m'')=\ar (m')=\ar(m)$.
Therefore $  \vert \vert wv^- \vert \vert_{br} \leq \eta  $.
Conversely, given a metamorphosis $M:wv^-\to \emptyset$,   extend
it by the identity on $v$ to  a metamorphosis $M':wv^-v\to v $.
Composing $M'$  with the inverse surgery $w \to wv^- v $, we
obtain a  metamorphosis $M'':w \to  v$ with $\ar (M'')=\ar
(M')=\ar(M)$. Therefore $\vert \vert wv^- \vert \vert_{br}  \geq
\eta  $. Hence $\rho_{br} (w,v)= \vert \vert wv^- \vert \vert_{br}
= \eta  $.
\end{proof}

We can  estimate the bridge norm and the bridge metric  via the
genus. Recall the  group $\pi=\pi(\alpha, \tau)$ from Sect.\ \ref{fipipipiib99}.

  \begin{theor}\label{tcdre53} Let $\varphi: \pi \to   \ZZ$ be an additive homomorphism  such that
  $\varphi(a)\in \{+1, -1 \} $ for all $a\in \alpha$. For  
  any nanoword $w$ over $\alpha$,
  $$      \vert \vert w\vert \vert_{br}\geq  \sigma_\varphi (p_w )/2, $$
where $p_w=p(w)$ is the $\alpha$-pairing associated with $w$.
\end{theor}

This theorem  yields a computable {\it a priori} estimate from below for the
total number of arches in any metamorphosis $w\to \emptyset$.  I do not know whether the assumption
$\varphi(a)\in \{+1, -1 \} $ for all $a\in \alpha$ is really
necessary here.

  Theorem \ref{tcdre53} and inequality
(\ref{zaz}) directly imply inequality (\ref{pop}). Note also the
following corollary.

 \begin{corol}\label{th:dfrovv53}  For   any nanowords $w, v$
and any    $\varphi $  as in Theorem \ref{tcdre53},
 $$\rho_{br} (w,v)\geq   \sigma_\varphi (p_w \oplus  p_v^-)/2.$$
\end{corol}

Indeed,
$$ \rho_{br} (w,v)=\vert \vert wv^-\vert \vert_{br}\geq \sigma_\varphi (p_{wv^-})/2= \sigma_\varphi (p_w \oplus  p_v^-)/2.$$

We now deduce Theorem \ref{tcdre53}  from the following lemma
whose proof, postponed to Sect.\ \ref{fdvmvmvmvmmsospmlensionf40},
uses topological techniques.

\begin{lemma}\label{tphi}
For any bridge move $m:w\to x$ and any $\varphi $ as in Theorem
\ref{tcdre53},
  \begin{equation}\label{bert}  \ar(m)\geq  \sigma_\varphi (p_w \oplus  p_x^-)/2.\end{equation}
\end{lemma}

\subsection{Proof  of Theorem \ref{tcdre53}}
We claim  that the inequality (\ref{bert}) holds for any move
$m:w\to x$ from the list (TR+). If $m$ is a bridge move, then this
is Lemma \ref{tphi}. If $m$ is an isomorphism or
  a homotopy move, then
$p_w, p_x$ are  cobordant and    $p_w\oplus  p_x^-$ is hyperbolic.
Then    $\ar(m)=0= \sigma_\varphi (p_w\oplus  p_x^-)$. If
(\ref{bert}) holds for $m: w \to x$, then it holds for
$m^{-1}:x\to w$ since
$$ \ar (m^{-1})=\ar(m)\geq  \sigma_\varphi (p_w \oplus
p_x^-)/2=\sigma_\varphi (p_x^- \oplus p_w)/2=\sigma_\varphi (p_x
\oplus p_w^-) /2.$$

Consider a metamorphosis $m:w\to x$ that splits as a composition
of two metamorphoses $m':w\to v$
 and $m'':v \to x$ where $v$ is  a nanoword. If   $m'$ and
 $m''$ satisfy (\ref{bert}), then so does $m$ since,
 by  (\ref{trian}),
$$\ar (m) =\ar(m') + \ar (m'') \geq  \sigma_\varphi
(p_w\oplus  p_v^-)/2+\sigma_\varphi (p_v\oplus  p_x^-)/2 \geq
\sigma_\varphi (p_w\oplus  p_x^-)/2.$$

We conclude that   (\ref{bert}) holds for all metamorphoses
$m:w\to x$.
 For  $x=\emptyset$,  this gives   $  \ar(m)\geq
\sigma_\varphi (p_w )/2$. Taking the minimum over all  $m:w\to
\emptyset$, we obtain the claim of the theorem.

 \section{Circular shifts and a weak bridge metric}\label{circubrimetr1za}

    \subsection{Shifts}\label{ushifts21}
 The {\it (circular) shift} of a
nanoword $(\A, w)$   is the nanoword $(\widehat  \A, \widehat  w )$ obtained
by moving the first letter $A=w(1) $ of $w$ to the end and
applying $\tau$ to $\vert A\vert \in \alpha$. More precisely,
$\widehat  \A =(\A-\{A\}) \cup \{\overline A\}$ where $\overline A$ is
a \lq\lq new" letter not belonging to $\A$. The   projection  $
\widehat  \A\to \alpha$ extends the   given  projection    $\A-\{A\}\to
\alpha$ by $\vert \overline A\vert=  \tau(\vert  A\vert)$.  The
word $\widehat  w$ in the alphabet $\widehat  \A$ is defined by $\widehat  w  = x
\overline A y \overline A$ for $w=AxAy$.

The $n$-th power of the shift transforms   a nanoword   of length
$n$    into itself. Hence  the inverse to the shift is a power of
the shift.

Two nanowords are {\it weakly cobordant} if they can be related by
a finite sequence of homotopy moves, surgeries, circular shifts
and inverse moves.
   For
example, for  $a,b \in \alpha $,  the shift transforms $w_{a,b}$
into $w_{b, \tau (a)}$. Therefore $w_{a,b}$ and  $w_{b, \tau (a)}$
are weakly cobordant.  If $a,b$  belong to different orbits of
$\tau$, then these two nanowords are not cobordant.

A simple invariant of weak cobordism   is provided by the
conjugacy class of $\gamma $: if nanowords $w, v$ are weakly
cobordant, then $\gamma (w), \gamma (v)\in \Pi$ are conjugate  in
$\Pi$.  This follows from Lemma \ref{25edcolbrad} and the identity
  $\gamma (\widehat  w)= z_{\vert w(1) \vert}^{-1} \gamma (w) z_{\vert w(1) \vert}$.
  In particular, if $\gamma (w)=1$, then $\gamma(v)=1$ for all nanowords $v$
  weakly cobordant to   $w$.

 We will see below    that  
    the genera and the $u$-polynomial of nanowords are
 weak cobordism invariants. Here
 we note the following result.

\begin{lemma}\label{ltytutu12} If    the cobordism class of a nanoword $w$ lies in
    $\Ker (p:\N_c \to {\mathcal{P}_{sk}})$, then all nanowords weakly cobordant to
     $w$ have the same property.
     \end{lemma}
                     \begin{proof}
                     It suffices to verify that if     $p_w$ is hyperbolic, then so is $p_{\widehat 
w}$ where $\widehat  w$ is obtained from $w$ by the shift. Let
$p_w=(S,s,e)$ and  $A=w(1)\in S-\{s\}$.  Then $p_{\widehat  w}= (\widehat 
S, s, \widehat  e)$ where $\widehat  S =(S-\{A\}) \cup \{\overline A\}$. A
direct computation shows that  $\widehat  e:\widehat  S \times \widehat  S \to
\pi$ is
     the unique skew-symmetric pairing
      such that $\widehat  e       \vert_{\widehat  S-\{\overline A\}}
=e\vert_{S-\{A\}}$ and $\widehat  e           (\overline A,B)=
e(2s-A,B)$ for  all $B\in S-\{A\} $. In particular, $e_A
(\overline A,s)=-e(A,s)  $. Observe now that any filling
$\lambda=\{\lambda_i\}_i$ of $p_{w}$ yields a filling $\widehat 
\lambda$ of $p_{\widehat  w}$ by changing the unique vector $\lambda_{i_0}$
in which the letter $A=w(1)$ occurs: if $\lambda_{i_0}=A$, then it is
replaced with $\overline A$; if   $ \lambda_{i_0}= \pm A+B$, then  $ \lambda_{i_0}$
is replaced with $\mp\, \overline A +B$; if $\lambda_{i_0}=A-B$, then
  $ \lambda_{i_0}$ is replaced with $ \overline A+B$. It is easy to see  that   if
$\lambda$ is an annihilating filling of $p_w$, then $\widehat  \lambda$
is an annihilating filling of $p_{\widehat  w}$. Therefore if $p_w$ is
hyperbolic, then so is $p_{\widehat  w}$.
\end{proof}

\subsection{The weak bridge pseudo-metric}\label{weakbrinmetricword}
The list (TR+) of moves on nanowords considered  in Sect.\
\ref{brinmetricword}  can   be extended to the following wider
list:

(TR++) isomorphisms,   homotopy moves,   bridge moves, circular
shifts, and  the inverse moves.

For nanowords $v,w$,  {\it a circular metamorphosis} $m:w \to v$
is a    finite sequence   $m=(m_1,...,m_n)$ of   moves from the list
(TR++)  transforming $w$ into $v$.  Set $\ar
(m)=\ar(m_1)+\cdots + \ar(m_n)$ where $\ar(m_i)$ is the number of
arches of $m_i$ if $m_i$ is a bridge move or its inverse and
$\ar(m_i)=0$ for all other moves. Set  
$$\rho_{wbr} (w,v)=\min_m \ar (m) \geq 0$$
where $m$ runs over all circular metamorphoses   $w\to v$. The
resulting function $\rho_{wbr}$ on $\N_c\times \N_c$  is a
pseudo-metric, i.e.,  it is symmetric, non-negative, satisfies the triangle inequality, and $\rho_{wbr} (w,w)=0$ for all $w$.  Lemma \ref{th:v54937ccccc153} implies that
$\rho_{br} (w,v) \geq \rho_{wbr} (w,v)$. The following theorem
yields an  estimate of $\rho_{wbr}$ from below  via the genus.

 \begin{theor}\label{thcircularv53}  For  any   nanowords $w, v$ and
 any
  $\varphi $    as in Theorem \ref{tcdre53},
 $$\rho_{wbr} (w,v)\geq   (\sigma_\varphi (p_w \oplus  p_v^-)-1)/2.$$
\end{theor}

In the next section, we   deduce Theorem \ref{thcircularv53} from
Lemma \ref{tphi}.

    \section{Weak cobordism  of   $\alpha$-pairings}\label{fi:133339912}

      Fix   a  ring    $R$ and a left $R$-module $\pi$.
      In this section we    study  algebraic properties of
         $\alpha$-pairings and apply them to  nanowords.

 \subsection{ Hyperbolic  $\alpha$-pairings}\label{2fi:gdr5999}
 The theory of fillings and hyperbolic $\alpha$-pairings
 extends  to tuples of $\alpha$-pairings (with values in $\pi$) as follows.
 Consider  a tuple of  $\alpha$-pairings
$p_1=(S_1,s_1,e_1)$,..., $p_r=(S_r,s_r,e_r)$ with $r\geq 1$.
Replacing these $\alpha$-pairings  by isomorphic  ones, we can
assume that the sets $S_1,...,S_r$ are disjoint.  Set
$S=\cup_{t=1}^r S_t$ and $S^\circ=\cup_{t=1}^r
S^\circ_t=S-\{s_1,...,s_r\}$.  Let $\Lambda=R S$ be the free
$R$-module   with basis $S$. Let $\Lambda_s$ be the submodule of
$\Lambda$ generated by the basis vectors  $s_1,...,s_r$.   A
vector $x\in \Lambda$ is  {\it weakly short} if   $x =  A \,
(\modu \Lambda_s)$ for   $A\in S^\circ $  or $x =  A+B\, (\modu
\Lambda_s)$ for   distinct   $A,B \in S^\circ$ with $\vert A\vert
= \vert B\vert $ or $x =  A-B\, (\modu \Lambda_s)$ for   distinct
$A,B \in S^\circ$ with  $\vert A\vert =\tau(\vert B\vert)$.
Removing the expression $(\modu \Lambda_s)$ in these formulas we
obtain a notion of a short vector.

  A  {\it weak  filling}  of the tuple
$p_1,...,p_r$ is  a finite family $  \{\lambda_i\}_i$ of  vectors
in $\Lambda$ such that one of $\lambda_i$ is equal to
$s_1+s_2+...+s_r$, all the other $\lambda_i$ are weakly short,
and  every element of  $S^\circ$ occurs in exactly one of   $
\lambda_i $ with non-zero coefficient (this
 coefficient is then
$\pm 1)$. The basis  vectors $ s_1,...,s_r $  may appear in
several $\lambda_i$ with non-zero coefficients.  For example, the
families $\{A\}_{A\in S^\circ} \cup \{s_1+s_2+...+s_r\}$ and
$\{A+s_1\}_{A\in S^\circ} \cup \{s_1+s_2+...+s_r\}$ are  weak
fillings
 of
$p_1,...,p_r$.

The pairings  $\{e_t:S_t\times S_t\to \pi\}_{t=1}^r$ induce a
bilinear form $e=\oplus_t e_t:\Lambda\times \Lambda\to \pi$ such
that $e  \vert_{S_t}=e_t $  and $e(S_t,S_{t'})=0$ for $t\neq t'$.
A weak filling  $ \{\lambda_i\}_i$ of $p_1,...,p_r$ is {\it
annihilating} if $e(\lambda_i,\lambda_j)=0$ for all $i,j$.  The
tuple $p_1,...,p_r$ is {\it hyperbolic} if it has an annihilating
weak filling.  The hyperbolicity is preserved under permutations
of $p_1,...,p_r$.

  For $r=1$, the notion of
a weak filling   is   wider than the one of a filling, cf.  Sect.\
\ref{fi:gdr5999}. Any weak filling $  \{\lambda_i\}_i$ of $p_1$
can be transformed into a filling of $p_1$   by adding appropriate
multiples of $s_1$ to all $\lambda_i \neq s_1$. Therefore, for
$r=1$, the notions of hyperbolicity introduced in this section and
in Sect.\ \ref{fi:gdr5999} are equivalent.

For $r=2$, the construction of Sect.\ \ref{summ799} shows  that
each filling of the $\alpha$-pairing $p_1\oplus p_2$ yields a weak
filling of the pair $(p_1, p_2)$. If the former is annihilating,
then so is the latter.
  We   conclude that if   $p_1\oplus
p_2$ is hyperbolic, then so is the pair $(p_1, p_2)$. The converse
may be not true.

\subsection{Weak cobordism}\label{fi:999hjf92} We say that    $\alpha$-pairings
$p , q$ are {\it   weakly cobordant} and write $p\simeq_{wc} q$ if
the pair $(p, q^-)$ is  hyperbolic. By the remarks above, if the
$\alpha$-pairing $p\oplus q^-$ is hyperbolic, then so is
  the pair $(p,q^-)$. Therefore cobordant
$\alpha$-pairings are weakly cobordant.

\begin{lemma}\label{l:gfgg12}  Weak cobordism of $\alpha$-pairings  is an equivalence relation.      \end{lemma}
                     \begin{proof}
It is clear that if a tuple of $\alpha$-pairings $p_1,...,p_r$ is
hyperbolic, then so is the tuple of opposite $\alpha$-pairings
$p^-_1,...,p^-_r$.  Thus, if a  pair  $(p, q^-)$ is  hyperbolic,
then so is the pair  $(p^-, q)$.  This implies the  symmetry of
the weak  cobordism.

The transitivity  of the weak cobordism is proven similarly to the
transitivity of cobordism in Lemma \ref{l:gg1} and we indicate
only the necessary changes. As $\lambda=\{\lambda_i\}_i$ (resp.\ $
\mu=\{\mu_j\}_{j }$), we take any weak filling of the pair $(p_1,
p_2^-)$ (resp.\ $(p_2, p_3^-)$). Before constructing $ \psi$, we
modify $ \lambda $ as follows. Let $\lambda_0$ be the vector of
  $\lambda$ equal to $s_1+s_2$. Adding appropriate
multiples of $\lambda_0 $ to the other $\lambda_i$, we can assume that
the basis vector $s_2 \in S_2 $ appears in all
$\{\lambda_i\}_{i\neq 0}$  with  coefficient $0$. This  does not
change the $R$-module $V_\lambda$ generated by $\{\lambda_i\}_i$.
Similarly, there is a vector $\mu_0$ of   $ \mu$  equal to
$s'_2+s_3$, and we can assume that $s'_2 \in S'_2$ appears in all
$\{\mu_j\}_{j\neq 0}$ with coefficient $0$. In the rest of the
proof instead of
    $q(\psi_K) = A\pm B $ and $q(\psi_K) = A $ it  should be
respectively $q(\psi_K) = A\pm B  \, (\modu R s_1 +R s_3)$ and
$q(\psi_K) = A   \, (\modu R s_1 +R s_3)$. Instead of
$\lambda_{i}=A+C $ and $\mu_{j}= - C' + B $ it should be
respectively $\lambda_{i}=A+C  \, (\modu R s_1)$ and  $\mu_{j}= -
C' + B  \, (\modu  R s_3)$, etc. The word \lq\lq short" should be
replaced with \lq\lq weakly short".
  \end{proof}

 \subsection{Invariants}\label{erfg9}     We can generalize the  genus of $\alpha$-pairings to   tuples
  as follows.
  Let   $F$ and $\varphi:\pi\to F$ be   as in Sect.\ \ref{u125452}.
    For a   tuple of  $\alpha$-pairings $p_1=(S_1,s_1,e_1),...,
p_r=(S_r,s_r,e_r)$, set $S=\cup_{t=1}^r S_t$  and let $e=\oplus_t
e_t:RS \times RS \to \pi$    as in Sect.\ \ref{2fi:gdr5999}.
  For a weak filling $\lambda=\{\lambda_i\}_i$ of the tuple
$p_1,...,p_r$, the matrix $(\varphi  e(\lambda_i, \lambda_j
))_{i,j}$ is a   square matrix over $F$.  Let $\sigma_\varphi
(\lambda)\in \frac{1}{2} \ZZ$ be half of its rank and
$$\sigma_\varphi(p_1,...,p_r)=\min_\lambda \sigma_\varphi (\lambda)\geq 0$$
where $\lambda$ runs over all weak fillings of $(p_1,...,p_r)$.
The half-integer $\sigma_\varphi(p_1,...,p_r)$ is called the {\it
$\varphi$-genus} of the tuple $p_1,...,p_r$.   It is obvious that
the $\varphi$-genus is preserved when $p_1,...,p_r$ are permuted
and $\sigma_\varphi(p_1^-,...,p_r^-)=\sigma_\varphi(p_1,...,p_r)$.
If $p_r=(\{s_r\},s_r, e_r=0)$ is the trivial $\alpha$-pairing, then
$\sigma_\varphi(p_1,...,p_r)=\sigma_\varphi(p_1,...,p_{r-1})$
 because then the vector $s_r\in RS$ lies in the annihilator of
$e$.  If the tuple $ p_1,...,p_r $ is hyperbolic, then
$\sigma_\varphi(p_1,...,p_r)=0$. If
  $ p_1,...,p_r $ are skew-symmetric,  then
$\sigma_\varphi(p_1,...,p_r)\in \ZZ$.

\begin{lemma}\label{thccvcvccvvccic3}
For any $\alpha$-pairings $p_1, p_2, p_3$,
\begin{equation}\label{ntrian}
\sigma_\varphi (p_1,   p_2^- ) + \sigma_\varphi (p_2 , p_3^-) \geq
\sigma_\varphi(p_1,  p_3^-) .\end{equation}
\end{lemma}
\begin{proof}     Pick    a  weak filling  $\lambda $
of $(p_1, p_2^-) $ such that $\sigma_\varphi (p_1,   p_2^-
)=\sigma_\varphi (\lambda)$. Pick    a   weak filling $ \mu $ of
$(p'_2, p_3^-) $ (where $p'_2$ is a copy of $p_2$) such that
$\sigma_\varphi (p_2 , p_3^-)=\sigma_\varphi (\mu)$. We modify
$\lambda$ and $\mu$ as in the proof of Lemma \ref{l:gfgg12}. This
modification preserves the $R$-modules $V_\lambda, V_\mu$
generated by these families of vectors and   therefore preserves $
\sigma_\varphi (\lambda)$ and  $ \sigma_\varphi (\mu)$. The rest
of the argument goes   as in the proof of Lemma \ref{th:verdic3}.
\end{proof}

\begin{theor} \label{1859s}    The $\varphi$-genus   of
  $\alpha$-pairings  is a   weak cobordism invariant.
 \end{theor}
\begin{proof}
 We need to prove that  $p_1 \simeq_{wc} p_2 \Rightarrow
 \sigma_\varphi(p_1)=\sigma_\varphi(p_2)$.
    The hyperbolicity of   the pair $p_1,
p_2^-$ implies that   $\sigma_\varphi(p_1,  p_2^-)=0 $. Applying
Lemma \ref{thccvcvccvvccic3} to the triple $p_1, p_2, p_3$ where
$p_3=(\{s\},s , e=0)$   is a trivial $\alpha$-pairing, we obtain the inequality 
$\sigma_\varphi(p_2)\geq \sigma_\varphi(p_1)$. By symmetry,
$\sigma_\varphi(p_1)= \sigma_\varphi(p_2)$.
 \end{proof}

\begin{lemma} \label{yyuyuys} For any   $\alpha$-pairings $p_1,p_2$,
 $$\sigma_\varphi (p_1 \oplus
 p_2)\geq \sigma_\varphi (p_1 ,  p_2) \geq \sigma_\varphi (p_1 \oplus
 p_2)-1.$$
 \end{lemma}
\begin{proof} Let $p_1=(S_1,s_1,e_1)$ and $p_2=(S_2,s_2,e_2)$. By  Sect.\ \ref{summ799},
every filling $\lambda$ of $p_1\oplus p_2$ yields a weak filling
of the pair $(p_1,p_2)$. Therefore $\sigma_\varphi(\lambda) \geq
\sigma_\varphi (p_1 ,  p_2)$. Taking minimum over all fillings
$\lambda$ of  $p_1\oplus p_2$, we obtain  $\sigma_\varphi (p_1 \oplus
 p_2)\geq \sigma_\varphi (p_1 ,  p_2)$. Conversely, any weak
 filling $\mu$ of the pair
$(p_1,p_2)$ gives rise to a filling $\mu'$ of $p_1\oplus p_2$ by
adding appropriate multiples of $s_1, s_2$ to all vectors of $\mu$
distinct from $s_1+s_2$. Let $V, V'$ be the submodules of $
RS_1\oplus RS_2$ generated respectively by $\mu, \mu'$. Clearly
$V'\subset V+Rs_1+Rs_2=V+Rs_1$ (since $s_1+s_2\in   V$). Therefore
the rank of the pairing $\varphi \circ (e_1\oplus e_2)$ restricted
to $V'$ does not exceed the rank of this pairing restricted to $V$
plus 2. For the half-ranks, we have  $\sigma_\varphi (\mu) \geq
\sigma_\varphi (\mu')-1 \geq \sigma_\varphi (p_1 \oplus
 p_2) -1$.
 Taking minimum over all $\mu$, we obtain
 $\sigma_\varphi (p_1 ,  p_2) \geq \sigma_\varphi (p_1 \oplus
 p_2) -1$. \end{proof}

 \subsection{Applications to nanowords}\label{aptona}
  \begin{lemma}\label{ledr:gnbnbnbg1} If  nanowords $w, v$ are weakly cobordant,
  then  $p(w)\simeq_{wc}  p(v)$.
           \end{lemma}
                     \begin{proof} We begin by defining
for any integer $m$, a transformation of skew-symmetric $\alpha$-pairings  called
{\it $m$-shift}. Consider   a skew-symmetric    $\alpha$-pairing $p=(S,s,e)$.
Pick $A\in S^\circ$ and replace it with a \lq\lq new" element
$\overline A$ such that $\vert \overline A\vert=\tau (\vert
A\vert)$. Endow the resulting  set $S_A =(S-\{A\}) \cup
\{\overline A\}$
  with  the unique skew-symmetric pairing  $e_A             :S_A             \times S_A             \to \pi$ such that $e_A             \vert_{S_A-\{\overline A\}}
=e\vert_{S-\{A\}}$ and $e_A             (\overline A,B)=
e(ms-A,B)$ for   $B\in S-\{A\} $. In particular, $e_A (\overline
A,s)=-e(A,s)  $.   We say that the $\alpha$-pairing $p_A=(S_A
,s,e_A             )$ is obtained from  $p$ by the $m$-shift   at
$A$. We claim that $p$ and $p_A$ are weakly cobordant. Consider a
copy $p'= (S'=\{B'\}_{ B\in S},s',e')$ of $p$  and  the weak
filling
 of  the pair  $(p_A, (p')^-)$ formed by the vectors  $\lambda_0=s +  s', \lambda_A= (\overline A- ms)  -A' $  and   $\{\lambda_B=B + B'\}_{B\in S-\{A\}}$.  This weak filling is
 annihilating. In particular,
 $$(e_A\oplus (e')^-) (\lambda_A, \lambda_0)=   e_A  (\overline A-ms , s )+
  (e')^-  ( - A',  s)= -e(A,s)+ e'(A',  s')=0,
  $$
 $$(e_A\oplus (e')^-) (\lambda_A, \lambda_B)=   e_A  (\overline A-m s , B )  +(e')^-  ( - A',  B')=- e(A,B)+  e'(A',  B')=0 .$$
  Thus
the pair $(p_A, (p')^-)$ is   hyperbolic so that $p_A\simeq_{wc}
p'\approx   p$.

To prove the lemma, we need only to show that   $p(\widehat  w
)\simeq_{wc} p(w)$, where $\widehat  w$ is obtained from $w$ by the
shift. The proof of Lemma \ref{ltytutu12} shows that   $p(\widehat  w)$
is obtained from $p(w)$ by the 2-shift at $A=w(1)$. Hence $p(\widehat 
w )\simeq_{wc} p(w)$.
\end{proof}

\begin{theor} \label{85956565656s}   For any additive homomorphism $\varphi$ from $\pi=\pi(\alpha, \tau)$
to a commutative domain, the $\varphi$-genus   of
 nanowords  is a
 weak cobordism invariant.
 \end{theor}

 This theorem follows from  Theorem \ref{1859s} and Lemma    \ref{ledr:gnbnbnbg1}.

 \subsection{Proof of Theorem \ref{thcircularv53}}\label{apprtjhohona}
  We claim that for  any nanowords $w, v$,
$$
\rho_{wbr} (w,v)\geq    \sigma_\varphi (p_w,  p_v^-) /2.$$
 By Lemma \ref{yyuyuys}, this will imply
 the theorem. By the definition of  $\rho_{wbr}$,
 it suffices to prove that  for any circular metamorphosis
$m:w\to v$,
  \begin{equation}\label{opl}\ar (m)\geq    \sigma_\varphi (p_w,  p_v^-)
  /2.
  \end{equation} If $m$ is a bridge move, then this inequality directly follows from
    Lemma
\ref{tphi} and the left inequality in Lemma \ref{yyuyuys}. If $m$
is an isomorphism or
  a homotopy move or a circular shift, then
$p_w, p_v$ are   weakly cobordant so that
    $\ar(m)=0= \sigma_\varphi (p_w,  p_v^-)$.
    If (\ref{opl}) holds for $m: w \to v$, then it holds for
$m^{-1}:v\to w$ since $ \ar (m^{-1})=\ar(m)$ and $ \sigma_\varphi
(p_w , p_v^-) =\sigma_\varphi (p_v^- , p_w) =\sigma_\varphi (p_v,
p_w^-)  $.   Finally, if a metamorphosis $m:w\to v$   splits as a
composition of  two  metamorphoses $m':w\to x$
 and $m'':x \to v$  satisfying (\ref{opl}), then
   Lemma \ref{thccvcvccvvccic3} ensures that
 $m$ also satisfies (\ref{opl}):
$$\ar (m) =\ar(m') + \ar (m'') \geq  \sigma_\varphi
(p_w,  p_x^-)/2+\sigma_\varphi (p_x,  p_v^-)/2 \geq \sigma_\varphi
(p_w,  p_v^-)/2.$$ We conclude that   (\ref{opl})  holds for all
circular metamorphoses $m:w\to v$.

 \subsection{Remarks}\label{fiertilorema9568} 1. The
 results obtained above for the $\varphi$-genera of   pairs  extend to tuples as follows.
 Pick an arbitrary   tuple $p_1,..., p_r$ of $\alpha$-pairings with $r\geq 1$.  Lemma
 \ref{thccvcvccvvccic3} generalizes to the following claim:
   for any $1\leq k \leq l \leq r$,
   $$   \sigma_\varphi (p_1,..., p_l)
+ \sigma_\varphi ( p_{k}^-,..., p_l^-,  p_{l+1},..., p_r)+l-k\geq
\sigma_\varphi(p_1,..., p_{k-1}, p_{l+1},...,p_r).$$ Setting here
$k=l=r-1$, we can deduce that $ \sigma_\varphi (p_1,..., p_r)$
depends only on the weak cobordism classes of  $p_1, ..., p_r$.
Lemma
 \ref{yyuyuys} generalizes to
$$\sigma_\varphi(p_1 \oplus  ... \oplus  p_r) \geq
  \sigma_\varphi(p_1,   ...,  p_r) \geq
   \sigma_\varphi(p_1 \oplus  ... \oplus  p_r) +1-r
   .$$
The arguments of Lemma \ref{l:gfgg12} extend to show that
   if the tuples $p_1, ..., p_k$ and
$p_k^-,p_{k+1},..., p_r$  are hyperbolic, then so is the tuple
$p_1,..., p_{k-1}$, $ p_{k+1},..., p_r$.
 Taking $k=1, r=2$, we obtain  that an $\alpha$-pairing  weakly cobordant
  to a hyperbolic $\alpha$-pairing is itself hyperbolic.

2. It is easy to show  that   the $u$-polynomial  
 of skew-symmetric (more generally, normal) $\alpha$-pairings is    invariant under weak cobordism.
Therefore the $u$-polynomial of nanowords  is   invariant under weak cobordism.

 \section{Words and loops}\label{Words and loops}

In this section we study nanowords over the 2-letter alphabet
$\alpha_0=\{+,-\}$ with involution $\tau_0$ permuting $+$ and $-$.
These nanowords are shown to be   disguised forms of generic loops
on surfaces. As an application, we   prove Lemma \ref{l:edlocfc12}.

\subsection{Loops}\label{v4040irss} By a {\it loop} $f:S^1\to \Sigma$, we   mean
  a generic  immersion of an oriented circle $S^1$ into an
  oriented connected surface $\Sigma$. A
loop may have  only a finite number of self-intersections which
are all double and transversal.  We shall sometimes use the term
\lq\lq loop"   for the set 
$f(S^1)$. A loop is {\it pointed} if it is endowed with a base
point (the origin) which is not a self-intersection. A loop
$f:S^1\to \Sigma$ is {\it spinal} if $\Sigma $ is a compact connected oriented surface  that deformation retracts on the set $f(S^1)$.
  Two
pointed spinal loops are {\it homeomorphic} if there is a an
orientation preserving homeomorphism of the ambient surfaces
mapping the first loop onto the second one keeping the origin and
the orientation of the loop.  

We    associate  with any
    pointed   loop $f$ a  nanoword  over  $\alpha_0$. To this end,
    label  the self-intersections of $f$  by (distinct) letters
    $A_1,...,A_m$ where $m$ is the number of self-intersections.
   Starting at the origin of $f$ and following along $f$
    in the positive direction  we write
   down the labels of all self-intersections until the   return to the
   origin. Since every self-intersection is traversed twice, this gives
   a   word $ w $   in the alphabet $\A=\{A_1,...,A_m\}$ such that
   every   $A_i$ appears in $w$ twice. The word $w$, called the Gauss word
    of $f$, was
   first constructed by  Gauss \cite{ga}. We define a
   projection $\A\to \alpha_0$ as follows.
   For
   $i=1,...,m$, we may speak about the first and second branches of $f$
   appearing at the
   first and second
   passages of $f$ through the self-intersection labelled by $A_i$.
    Let $ t_i^1$ (resp.\ $t_i^2$) be a positively oriented  tangent
   vector of  the first (resp.\ second) branch of $f$
    at this self-intersection. Set $\vert A_i\vert =+ $ if the
   pair $(t_i^1, t_i^2)$ is positively oriented and $\vert A_i\vert =-  $
    otherwise. This makes $(\A, w)$
     into a nanoword over $\alpha_0$ of length $2m$. It  is well defined up to
   isomorphism and is called the {\it underlying
      nanoword} of $f$.   Obviously,    homeomorphic loops
   have isomorphic underlying nanowords.

   \begin{theor}\label{edrdcppmolbrad}
    The map assigning to a pointed loop its underlying nanoword establishes a
     bijective correspondence between the set of
     homeomorphism classes of pointed spinal  loops and
     the set of isomorphism classes
    of nanowords over $\alpha_0=\{+,-\}$.
 \end{theor}

    \begin{proof}     Given a
nanoword  $(\A, w:\widehat  n \to \A)$   over $\alpha_0$ we define a
pointed spinal loop as follows. Let $S^1=\RR \cup \{\infty\}$ be
the circle
   obtained by the compactification of the line $\RR$ with  right-handed orientation. Since
   every letter of $\A$ appears in $w$ twice, the family
   $\{w^{-1}(A)\}_{A\in \A}$ is a partition  of the set $  \widehat  n
   \subset \RR\subset S^1$ into pairs.
     Identifying the elements of $w^{-1}(A)$ for every $A\in \A$,
  we transform $S^1$ into   a graph (i.e., a 1-dimensional
  CW-complex)
$\Gamma=\Gamma_w$. This graph has $n$   edges, which we endow with
orientation  induced by the one in $S^1$,  and $n/2$ four-valent
vertices $\{V_A\}_{A\in \A}$ where $V_A$ is the image of
$w^{-1}(A)$ under the projection $S^1\to \Gamma$. Next, we thicken
$\Gamma$ to a surface $\Sigma=\Sigma_w$.
 If $n=0$ (so that $w=\emptyset$), then
$\Gamma=S^1\subset \Sigma= S^1\times [-1, +1]$. Assume    
that $n\geq 2$.    A neighborhood of a vertex $V_A\in \Gamma$
embeds into a copy $d_A$ of the standard  unit 2-disk   $
\{(p,q)\in \RR^2\,\vert \, p^2+q^2\leq 1\}$ as follows.   Suppose
that $w^{-1}(A)=\{i,j\}$ with  $1 \leq i<j\leq n$.  Note that any
point $x\in S^1$ splits its small neighborhood in $S^1$ into two
oriented arcs,    incoming and  
outgoing with respect to $x$. A neighborhood  of $V_A$ in $\Gamma$
consists of four  arcs which can be identified with   incoming and
outgoing arcs of $i, j$ on $S^1$. We embed this neighborhood into
$d_A$ so that $V_A$ goes to the origin $(0,0)$ and the incoming
(resp.\ outgoing) arcs of $i,j$ go to the intervals $[-1,0]\times
0$,  $0 \times [-1,0]$ (resp.\ $[0,1]\times 0$, $0 \times [0,1]$),
respectively. We endow $d_A$ with counterclockwise orientation if
$\vert A \vert =+$ and with clockwise orientation if $\vert A
\vert =-$. In this way the vertices of $\Gamma$ are thickened
 to disjoint oriented  copies of the unit 2-disk.
An  edge  of $\Gamma$ leads from a vertex, $V_A$, to a vertex,
$V_B$, (possibly $A=B$). Its thickening is the union of $d_A,d_B$
and   a ribbon $R_{A,B}$ connecting these 2-disks. The ribbon
$R_{A,B}$  is a copy of the rectangle $   [0,1]\times
[-1/10,+1/10]$ endowed with counterclockwise orientation. The copies in
$R_{A,B}$ of the intervals $ 0\times  [-1/10,+1/10]$, $1\times
[-1/10,+1/10]$, $   [0,1]\times 0$ are called the left side, the
right side, and the core of $R_{A,B}$,  respectively. It is
understood that  $R_{A,B}$ meets $\Gamma$ along its core and meets
$d_A\cup d_B$ along its  sides. More precisely, the ribbon
$R_{A,B}$ is glued to the disk $d_A$ (resp.\ $d_B$) along a
length-preserving  embedding of its left (resp.\ right) side into
the boundary of the disk such that the orientations of this disk
and $R_{A,B}$ are compatible.
 Thickening in this way
all  the vertices and edges of $\Gamma$, we embed $\Gamma$ into a
compact   connected  oriented surface  $ \Sigma$.      Composing the
projection $S^1\to \Gamma$ with the inclusion $\Gamma
\hookrightarrow \Sigma $, we obtain a
spinal loop $f:S^1\to \Sigma$   with origin $f(0)$ for $0\in
\RR\subset S^1$. It is straightforward to see that the underlying
nanoword of $f$ is isomorphic to $w$. Applying this construction
to the underlying nanoword  of a pointed  spinal  loop, we obtain a
homeomorphic pointed loop. This proves the claim of the theorem.
\end{proof}

\begin{corol}\label{edrcorrrolorolbrad}
There is a bijective correspondence between  the set of
     homeomorphism classes of non-pointed spinal loops and the set of isomorphism
classes
    of nanowords over $\alpha_0 $ considered up to shifts.
 \end{corol}

 It suffices to observe  that when the base point of a loop is pushed
along the loop across a self-intersection, the corresponding
nanoword over $\alpha_0$ changes via the circular shift determined
by $\tau_0$.

              \subsection{Homological computations}\label{fi:g32}
We  analyze in more detail the relationships between a nanoword
$(\A, w:\widehat  n \to \A )$ over $\alpha_0$ and the corresponding
pointed spinal loop  $f:S^1\to \Sigma=\Sigma_w$  constructed in
Theorem \ref{edrdcppmolbrad}.     The orientation of $\Sigma$
determines a skew-symmetric intersection pairing $b:
H_1(\Sigma)\times H_1(\Sigma) \to \ZZ$, where
$H_1(\Sigma)=H_1(\Sigma;\ZZ)$.   By abuse of notation, the homological
intersection number of two loops $x,y$ in $\Sigma$ will be denoted
$b(x,y)$. Thus, $b(x,y)=b([x], [y])$, where $[x], [y]\in
H_1(\Sigma)$ are the homology classes of $x,y$, respectively.
To compute $b(x,y)$, one deforms $x,y$ on $\Sigma$ so
that they have only a finite number of intersections which are all
transversal and distinct from the self-crossings of $x,y$. Then
$b(x,y)$ is equal to the number of intersections where $x$
crosses $y$ from left to right  minus the number of intersections
where $x$ crosses $y$ from right to left.

   For a letter $A\in \A$,
 we define a loop on $ \Sigma $ as follows.
 Let $w^{-1}(A)= \{i,j\}$ with $1\leq i<j \leq n$.
 Since $f(i)=f(j)$,
 the map $f$ transforms the  interval $  [i,j]\subset \RR \subset S^1$, oriented
from $i$ to $j$, into a loop  on $\Sigma$ with origin
$V_A=f(i)=f(j)$. This loop is denoted $f_A$.

Recall the abelian group $\pi=\pi(\alpha_0,\tau_0)$ generated by
the elements of $\alpha_0$ subject to the relations $a+\tau_0
(a)=0$. The group homomorphism $\pi\to \ZZ$ sending $+$ to $+1$
and $-$ to $-1$ is an isomorphism, and we   use it to identify
$\pi$ with $\ZZ$. The $\alpha_0$-pairing $(S=\A\cup \{s\},
s, e_w:S\times S\to \pi=\ZZ)$ associated with $w$ is related to $b:
H_1(\Sigma)\times H_1(\Sigma) \to \ZZ$ as follows.

\begin{lemma}\label{ltye4rtutu12}   For any $A\in \A$,
\begin{equation}\label{lp}e_w(A,s)=   b(f_A,f).\end{equation} For any $A,B\in \A$,
\begin{equation}\label{lp2}e_w(A,B)= 2\, b(f_A,f_B).\end{equation}
     \end{lemma}
                     \begin{proof}  We need an additional piece of notation.
 Let as above $A\in \A$ and $w^{-1}(A)= \{i,j\}$ with $1\leq i<j \leq n$.
   Denote by $[j,i]$ the oriented interval in 
$S^1$  going from $j$ to $+\infty=-\infty$ and then from $-\infty$
to $i$. Thus,  $[i,j]\cup [j,i]=S^1$ and $[i,j]\cap [j,i]=\{i,j\}$.
The mapping $f$ transforms $[j,i]$ into a loop, $f^-_A$, on
$\Sigma$ such that $[f_A]+[f^-_A]=[f]$. Drawing a picture of the
loops $f_A, f^-_A$ in
  the disk neighborhood $d_A$ of their common origin $V_A=f(i)=f(j)$, one
observes that a little deformation makes  $f_A, f^-_A$ disjoint in
$d_A$. Outside  $d_A$,  these loops meet transversely at the
points $\{V_D\}$, where $D$ runs over letters in $ \A$ such that
either $w=\cdots A \cdots D\cdots A \cdots D \cdots $ or $w=\cdots
D \cdots A\cdots D \cdots A \cdots$. The intersection sign of
$f_A, f^-_A$ at $V_D$ is $\vert D\vert\in \alpha_0=\{\pm \}$ in the first case and
$-\vert D\vert$ in the second case. Therefore
$$b(f_A, f^-_A)=\sum_{D\in \A} n_w (A, D) \vert D\vert=e_w(A,s).$$
 This implies  Formula
(\ref{lp}):
$$ b(f_A,f)= b(f_A, f_A)+b(f_A,  f^-_A) = b(f_A,  f^-_A)=e_w(A,s).$$

 Let us prove Formula
(\ref{lp2}). If $A=B$, then both sides are equal to 0.  Assume
 that $A\neq B$. Let $w^{-1}(A)=\{i,j\}$ with $i <j $
and $w^{-1}(B)=\{\mu, \nu\}$ with $\mu <\nu $. Note that the
numbers $i,j,\mu, \nu$ are pairwise distinct. By the skew-symmetry
of $e_w$ and $b$, if Formula (\ref{lp2}) holds for   $A,B$, then
it also holds for   $B,A$. Permuting if necessary $A$ and $B$, we
can assume that $i<\mu$. We distinguish three cases depending on the
order of $j,\mu, \nu$.

Case   $i<j<\mu<\nu$. Then $w=x' Ax A y B z B z' $ where $x,y,z,
x',z'$ are words in the alphabet $\A$. Observe that the intervals
$ [i,j] $ and $ [\mu, \nu]  $ are disjoint. Therefore the loops
$f_A, f_B$ meet transversely at the points $\{V_D\}$ where $D$
runs over letters in $ \A$ which appear once between the entries
of $A$ and once
  between the entries of
$B$. The intersection sign of $f_A, f_B$ at $V_D$ is $\vert
D\vert$. Therefore, in the notation of Sect.\
\ref{fi2ampsortnib99},
 $b(f_A, f_B)=\langle x,  z\rangle $. Formula \ref{mn1} implies that
 $e_w(A,B)= 2\, b(f_A, f_B)$.

Case   $i<\mu< \nu<j$. Then $w=x' Ax B y B z A z' $ where $x,y,z,
x',z'$ are words in the alphabet $\A$. Observe that the intervals
$ [j,i] $ and $ [\mu,\nu] $ on $ S^1$ are disjoint. Therefore the
loops $f_A^-, f_B$ meet transversely at the points $\{V_D\}$ where
$D$ runs over letters in $ \A$ which appear once between the
entries of $B$ and once
 before the first entry of $A$ or   after the last entry of
$A$. The intersection sign of $f^-_A, f_B$ at $V_D$ is $\vert
D\vert$ in the first case and $-\vert D\vert$ in the second case.
Therefore
$$b(f^-_A, f_B)=\langle x',y\rangle -\langle y, z'\rangle .$$
As we know,
$$ b(f, f_B)=- b(f_B,f)  =-e_w(B,s)=
\langle x',y\rangle +\langle x,y\rangle -\langle y,z \rangle
-\langle y, z'\rangle.$$ Then
$$b(f_A, f_B) =b(f, f_B)-
b(f^-_A, f_B)= \langle x,y\rangle -\langle y,z \rangle.$$ Now,
Formula \ref{mn} implies that
 $e_w(A,B)= 2\, b(f_A, f_B)$.

 Case   $i<\mu<j<\nu$. Then $w=x' Ax B y A z B z' $ where $x,y,z,
x',z'$ are words in the alphabet $\A$.  This case is more involved
since neither the loops $f_A, f_B$ nor the complementary loops are
transversal. Note that composing the projection $\A\to \alpha_0$
with $\tau_0:\alpha_0\to \alpha_0$ we obtain a new nanoword
$(\A,w')$ over $\alpha_0$ such that $e_{w'}(A,B)=-e_w(A,B)$. The
spinal loop corresponding to $w'$ is obtained from $f$   by
reversing  orientation in the ambient surface;  the associated
intersection form is $-b$.  Therefore, replacing if necessary $w$
by $w'$, we can assume     that $\vert B\vert =+$.
 Choose  coordinates $(p,q)$ in
the disk neighborhood $d_B\subset \Sigma$ of the point
$V_B=f(\mu)=f(\nu)$  so that
 $f_A\cap d_B$ is the line
  $q=0$,  $f_B\cap d_B$ is   the union of
half-lines $p=0, q\leq 0$ and $q=0, p\geq 0$,  and the orientation
on $f_A,f_B$ is right-handed   on the latter half-line. Since
$\vert B\vert =+$,  the coordinates $(p,q)$ determine the
orientation of $\Sigma$. Pushing $f_B$ slightly to its left in
$\Sigma$, we obtain a \lq\lq parallel" loop, $f'_B$, transversal
to $f_A$. We can assume that $f'_B\cap d_B$ is the union of
 half-lines $p=-1, q\leq 1$ and $q=1, p\geq -1$.

 To compute
$b(f_A, f_B)=b([f_A], [f_B])= b(f_A, f'_B)$, we  split the set
$f_A\cap f'_B$ into five disjoint subsets. The first of them
consists of the single intersection of $f_A$ and $f'_B$ in $ d_B$,
given in the coordinates above by $p=-1,q=0$. The intersection
sign of $f_A, f'_B$ at this point is  $+1$. The second subset of
$f_A\cap f'_B$ consists of the intersections of $f_A$ and $f'_B$
in the disk neighborhood $d_A$ of $V_A=f(i)=f(j)$. An inspection
shows that if $\vert A\vert =+1 \in \ZZ$, then $f_A$ and $f'_B$ do
not meet in $d_A$ and if $\vert A\vert =-1\in \ZZ$, then $f_A$ and
$f'_B$ meet transversely in one point in $d_A$ and their 
intersection sign at this point is $-1$. The joint contribution of
the first and second sets to $b(f_A, f'_B)$ is equal to   $( \vert
A\vert+1)/2 = ( \vert A\vert+ \vert B\vert)/2$. The third subset
of $f_A\cap f'_B$ is $f([i,\mu])\cap f'_B$; its points $\{V_D\}$
are numerated by letters $D$   which appear once between the first
entry of $A$ and the first entry of $B$ and once
  between the entries of
$B$. The intersection sign of $f_A, f'_B$ at such $V_D$ is $\vert
D\vert$.  The contribution of these crossings to $b(f_A , f'_B)$
is equal to $\langle x, y\rangle+\langle x, z\rangle$.  The forth
subset of $f_A\cap f'_B$ is numerated by the   crossings of
$f([\mu,j])$ with the part of $f'_B$ obtained by pushing
$f([j,\nu])\subset f_B$ to the left; they are numerated by letters
$D$   which appear once in $y$ and once in $z$. These crossings
contribute $\langle y, z\rangle$ to $b(f_A, f'_B)$. The remaining
subset of $f_A\cap f'_B$ is numerated by the self-crossings of
$f([\mu,j])$: each of them gives rise to two points of $f_A\cap
f'_B$ with opposite intersection signs. This   subset contributes
0 to $b(f_A, f'_B)$. Summing up these contributions we obtain
$$b(f_A, f_B)=  b(f_A, f'_B)= ( \vert A\vert+ \vert
B\vert)/2 +  \langle x, y\rangle+\langle x, z\rangle+\langle y,
z\rangle.$$ Now, Formula \ref{x56}  implies that
 $e_w(A,B)= 2\, b(f_A, f_B)$.
\end{proof}

 \subsection{Proof of Lemma \ref{l:edlocfc12}}\label{fi:mxlxlxlxlx40}
      The idea of the
proof is as follows. Let $f:S^1\to \Sigma$ be the pointed spinal
loop constructed from $w$ in Theorem \ref{edrdcppmolbrad}. Lemma
\ref{ltye4rtutu12} allows us to interpret the expression
$e_w(\lambda_{B_1}, \lambda_{B_2})$    in Formula \ref{plo} as an
intersection number  of certain loops on $\Sigma$ associated with
$B_1,B_2$. We   construct an oriented  3-dimensional manifold $M$,
depending on the nanophrase $\nabla$,  such that $\Sigma \subset
\partial M$ and the loops on $\Sigma$ associated with all $B\in
\B$ are homologically trivial in $M$. This implies that  the
intersection number of two such loops, $e_w(\lambda_{B_1},
\lambda_{B_2})$, is equal to 0. Other equalities in Formula
\ref{plo} are proven similarly. The construction of $M$ needs
a few  preliminaries which we now discuss.

We keep notation introduced in the second paragraph of Sect.\
\ref{fsymnano9} and in the proof of Theorem
\ref{edrdcppmolbrad}.  
  Thus,  each letter $A\in \A$ gives rise to a
self-intersection  of $f$ and to its disk neighborhood $d_A$ which
is a copy of the   unit 2-disk $\{(p,q)\in \RR^2\,\vert\,
p^2+q^2\leq 1\}$. The curve $f$ traverses $d_A$ first time along
$[-1,+1]\times 0$ and second time along $0\times [-1,+1]$, both
times   from $-1$ to $+1$. We call the points
$(-1,0), (1,0), (0,-1) , (0,1)\in \partial d_A$, respectively, the
first input, the first output,  the second input, and the second
output of $d_A$.
 Each
consecutive pair of letters $A,B$ in $w$ gives rise to a ribbon
$R_{A,B}\subset \Sigma$ which is a copy of the rectangle $
\{(p,q)\in \RR^2\,\vert \, p\in [0,1], q\in [-1/10,+1/10]\}$
endowed with counterclockwise orientation. The curve $f$ traverses the
ribbon $R_{A,B}$ once along its core $[0,1]\times 0$  in the
direction from $0$ to $1$. Warning: the notation $R_{A,B}$ may be
misleading since this ribbon depends not only on $A,B$ but on the
exact position of $AB$ in $w$: if the sequence of two consecutive
letters $AB$ occurs in $w$ twice, then it gives rise to two distinct
ribbons. In our arguments it will be always clear which sequence
$AB$ is implied. One more ribbon  $R_{w(n),w(1)}$  in $\Sigma$ is obtained by
thickening the interval $[n,1]\subset S^1$, where   $n$ is the length of $w$. This ribbon connects $d_{w(n)}$
to $d_{w(1)}$.  Each ribbon $R_{A,B}$ meets $d_A,d_B$ along its
sides; otherwise these $n$ ribbons and $n/2$ disks   are disjoint.

Let $\Sigma_1$ be the compact subsurface of $\Sigma$ formed by the disks
$\{d_B\}_{B\in \B}$ and the ribbons $R_{B_1, B_2}$ associated with
pairs of consecutive letters $B_1,B_2$ in $w$ contained in one of
the words $v_1,...,v_k$ forming the nanophrase $\nabla$. (The number of such pairs is equal to
$2\card (\B)-k$. Note that if there are no letters in $w$ between
$v_{i-1}$ and $v_i$, then the pair consisting of the last letter
of $v_{i-1}$ and the first letter of $v_i$ does not contribute to
$\Sigma_1$.)  The orientation of $\Sigma$ induces an
orientation of $\Sigma_1$.

We  define  an orientation reversing involution   $ I:\Sigma_1\to
\Sigma_1 $. We begin by defining it on  $\cup_{B\in \B}\,  d_B$. For
$B\in \B$, let $I_B:d_B\to d_{B^\iota}$  be the homeomorphism
acting as follows:  a point on $d_B$ with coordinates $(p,q)$
   goes to the point on $d_{B^\iota}$
  with coordinates $(-p,-q)$ if
   $\varepsilon  (B)=1$ and with coordinates $(-q,-p)$
   if $\varepsilon  (B)=0$.  Recall that   $d_B$
   is oriented counterclockwise (with respect to the coordinates $p,q$) if
$\vert B\vert =+$ and clockwise if $\vert B\vert =-$. That $I_B$
is orientation
  reversing follows from  the
 assumption that $\vert B\vert =\tau_0^{\varepsilon (B)}(\vert B^\iota \vert)$.
 The equality $\varepsilon (B)=\varepsilon
(B^\iota)$ implies that $I_{B^\iota} I_B=\id$.
 Note that $I_B$ transforms the outputs into the
inputs and {\it vice versa}. More precisely, if $\varepsilon
(B)=1$, then $I_B$ sends   the $i$-th output of $d_B$   to the
$i$-th input of $d_{B^\iota}$ for $i=1,2$. If $\varepsilon (B)=0$,
then $I_B$ sends the $i$-th output of $d_B$   to the $(3-i)$-th
input of $d_{B^\iota}$ for $i=1,2$.

We   define a similar involution on the    ribbons forming
$\Sigma_1$. Consider the ribbon $R_{B_1,B_2}\subset \Sigma_1$
arising from a 2-letter segment $B_1 B_2$ in $v_r$ where $1\leq r
\leq k$. Let $B_2^\iota B_1^\iota$ be the symmetric segment in
$v_r$: if  $B_1, B_2$ appear on the $j$-th and $(j+1)$'st positions
in  $v_r$ and the length of $v_r$ is $n_r$, then the symmetric
segment is formed by the letters appearing on the $(n_r-j)$-th and
$(n_r+1-j)$-th positions in $v_r$. We define a homeomorphism
$I_{B_1,B_2}:R_{B_1,B_2}\to R_{B_2^\iota, B_1^\iota}$ using the
coordinates $(p,q)$ on these ribbons:
 a point on $R_{B_1,B_2}$ with coordinates $(p,q)$
   goes to the point on $R_{B_2^\iota, B_1^\iota}$
  with coordinates $(1-p,q)$. This homeomorphism is orientation
  reversing and exchanges the sides   left $\leftrightarrow $ right  of the
    ribbons. We claim
   that $I_{B_1,B_2}$ coincides with
 $I_{B_1 }:d_{B_1 } \to   d_{B_1^\iota}$ on $R_{B_1,B_2}\cap d_{B_1 }$,
  i.e., on the left side of  $R_{B_1,B_2}$. Since both these
  homeomorphisms are orientation reversing and length preserving,
  it suffices to check
that $I_{B_1 }$  sends the output of $d_{B_1 }$ lying on the left
side of $R_{B_1,B_2}$ (in its metric center) into the input of
$d_{B_1^\iota}$ lying on the right side of $R_{B_2^\iota,
B_1^\iota}$ (again in its metric center). If $\varepsilon (B_1)=0$
and the entry of $B_1$ in question is its $i$-th entry in $v_r$
with $i=1,2$, then the entry of $B_1^\iota$ in question is its
$(3-i)$-th entry in $v_r$. Thus, $R_{B_1,B_2}$ is incident to the
$i$-th  output of $d_{B_1 }$ and $R_{B_2^\iota, B_1^\iota}$ is
incident to the $(3-i)$-th input of $d_{B_1^\iota}$. As observed
above, these output and input are related by $I_{B_1 }$.
Similarly, if $\varepsilon (B_1)=1$ and the entry of $B_1$ in
question is its $i$-th entry in $w$ with $i=1,2$, then the entry
of $B_1^\iota$ in question is also its $i$-th entry in $w$. Thus
$R_{B_1,B_2}$ is incident to the $i$-th  output of $d_{B_1 }$ and
$R_{B_2^\iota, B_1^\iota}$ is incident to the $i$-th input of
$d_{B_1^\iota}$. These output and input are related by $I_{B_1 }$. A
similar argument shows that the homeomorphism $I_{B_1,B_2}$  is
compatible with
 $I_{B_2 }:d_{B_2 } \to   d_{B_2^\iota}$; indeed the latter
sends the input of $d_{B_2 }$ lying on the right side of
$R_{B_1,B_2}$ to the output of $d_{B_2^\iota}$ lying on the left
side of $R_{B_2^\iota, B_1^\iota}$.

We conclude that  the homeomorphisms $ \{I_B\} $ and
$\{I_{B_1,B_2}\}$   extend to  an orientation reversing
homeomorphism $I:\Sigma_1\to \Sigma_1$. Clearly, $I^2=\id$.   We
describe the set of fixed points $Fix (I)$ of $I$. If $B\neq
B^\iota$, then $I\vert_{d_B}=I_B:d_B \to d_{B^\iota}$ has no fixed
points. If $B= B^\iota$, then $\varepsilon (B)=0$ and
$I\vert_{d_B}=I_B:d_B \to d_{B}$ is defined by $(p,q)\mapsto
(-q,-p)$. The set $Fix (I) \cap d_B$ is then the interval
$p+q=0$ connecting  the points of $\partial d_B$ with coordinates
$(-1/\sqrt{2},1/\sqrt{2})$ and $(1/\sqrt{2},-1/\sqrt{2})$. Both
these points lie on $\partial \Sigma_1$. Similarly, the
homeomorphism $I_{B_1,B_2}:R_{B_1,B_2}\to R_{B_2^\iota,
B_1^\iota}$ may have  fixed points if and only if the 2-letter
segment $B_1 B_2$ in $v_r$ lies precisely in the center of $v_r$.
If it is the case, then $I_{B_1,B_2}$ is an involution on $
R_{B_1,B_2}$ given  by $(p,q)\mapsto (1-p,q)$. Its set of fixed
points  is the interval $(1/2) \times [-1/10,1/10]\subset
R_{B_1,B_2}$ with endpoints on $\partial \Sigma_1$. This interval
meets $ f(S^1)$ in one point  with coordinates $(1/2, 0)$. Since
the length of $v_r$   is even for all $r=1,..., k$, the word $v_r$
has a unique central 2-letter sequence which gives rise to a
component of $Fix (I)$.

It is crucial for the sequel that   for all $r=1,...,k$,    the sub-path of $f$
corresponding to $v_r$ lies in $\Sigma_1$ and is folded by $I$ in two at the middle point. More
precisely, if $i,j\in \widehat  n$ numerate the first and the last
letter of $v_r$, then $f([i,j])\subset \Sigma_1$. For all
$u=0,1,..., j-i-1$, we have $I(f(i+u))=f(j-u)$ and $I$ maps the
arc $f([i+u,i+u+1])$ bijectively onto $f([j-u-1,j-u])$ (reversing
orientation). This shows that each path $f([i+u,j-u])$ is folded
in two in the quotient  $\Sigma_1/I$, that is it becomes a loop of type $\delta\delta^{-1}$, 
where $\delta$ is a path in $\Sigma_1/I$ and $\delta^{-1}$ is the inverse
path. Such a loop is contractible in $\Sigma_1/I$.

 Let $M$ be the
topological space obtained from   the cylinder $\Sigma \times
[0,1]$ by the identification   $a\times 1 = I(a)\times 1$ for all
$a\in \Sigma_1$. An   inspection of   neighborhoods of points
shows that $M$ is a  3-manifold. The fact that
 $I$ is orientation-reversing implies that $M$ is orientable.
  We identify $\Sigma$ with   $\Sigma\times
0 \subset \partial M$ and denote  by $b$ the intersection pairing
$H_1(\Sigma) \times H_1(\Sigma) \to \ZZ$.  It is well-known   that
$b(x,y)=0$ for any
   $x,y\in H_1(\Sigma)$ lying in the kernel of the
inclusion homomorphism $H_1(\Sigma) \to H_1(M)$.

We can now prove that $e_w (\lambda_{B_1}, \lambda_{B_2})=0$ for
all $B_1, B_2\in \B$.  Define an additive homomorphism $\rho:\ZZ
\A \to H_1(\Sigma)$ by $\rho (A)=[f_A]$ for $A\in \A$, where $f_A$
is the loop introduced  in Sect.\ \ref{fi:g32}. By Lemma
\ref{ltye4rtutu12}, $e_w (\lambda_{B_1}, \lambda_{B_2})=2\, b
(\rho(\lambda_{B_1}), \rho(\lambda_{B_2}))$. To prove the equality
$e_w (\lambda_{B_1}, \lambda_{B_2})=0$ it suffices to prove the
following claim:

$(\ast)$  for all $B\in \B$, the homology class $\rho
(\lambda_{B})\in H_1(\Sigma)$ is homologically trivial in $M$.

Suppose first that  $B^\iota=B$ so that  $B$ appears twice in the same
word $v_r$ on symmetric spots. The loop $f_B=f_B\times 0$ on
$\Sigma=\Sigma\times 0$ is obviously homotopic to the loop
$f_B\times 1$ in $\Sigma\times [0,1]$.  By the argument  above,
the latter loop lies on  $\Sigma_1\times 1$ and projects to a
contractible loop in $(\Sigma_1\times 1)  /I\subset M$.  Therefore
the loop $f_B$ is  contractible in $M$.
 Hence $\rho (\lambda_{B})=\rho(B)=[f_B]$
is homologically trivial in $M$.

 Suppose now  that $B^\iota\neq B$.  Let $i<j$ (resp.\ $\mu<\nu$) be the
numbers numerating the entries of $B$ (resp.\ of $B^\iota $) in $w$.
Exchanging if necessary $B, B^\iota$ and using that
$\lambda_{B^\iota}=\pm \lambda_B$, we can assume that $i<\mu$.

Consider first the case where $\varepsilon (B)=0$. Then
$\lambda_{B}=B+B^\iota$  and  $B, B^\iota$ appear twice in the
same word $v_r$ with $r=1,...,k$. The definition of $\iota$
implies that either $i<j<\mu<\nu $ or $i<\mu<j<\nu $.  If
$i<j<\mu<\nu $, then the path $ f\vert_{[i,\nu]}$ is
  the  product of the loop $f_B=f\vert_{[i,j]}$, the path  $
f\vert_{[j,\mu]}$, and the loop $f_{B^\iota}=f\vert_{[\mu,\nu]}$.
Therefore $\rho (\lambda_{B})=[f_B]+[f_{B^\iota}]$ is
 the homology class of the loop $f\vert_{[i,\nu]}
(f\vert_{[j,\mu]})^{-1}$ in $\Sigma_1$. Both paths forming the
letter loop project to contractible loops in $\Sigma_1/I$. This
implies  $(\ast)$. If $i<\mu<j<\nu $, then the path $
f\vert_{[i,\nu]}$ is
  the  product of the loop $f_B=f\vert_{[i,j]}$, the path $
(f\vert_{[\mu,j]})^{-1}$, and the loop
$f_{B^\iota}=f\vert_{[\mu,\nu]}$. Therefore $\rho
(\lambda_{B})=[f_B]+[f_{B^\iota}]$ is
  the homology class of the loop $f\vert_{[i,\nu]}
 f\vert_{[\mu,j]} $. As above, this
implies   $(\ast)$.

Consider the case where  $\varepsilon (B)=1$. Then $\lambda_{B}=B-B^\iota$
and   both  $B$ and  $B^\iota$ appear  once in $v_r$ and once in $v_{r'}$
with $r<r'$. We have  either $i<\mu<j<\nu $ or $i<\mu<\nu<j $.  In
the first case $f_B=f\vert_{[i,j]}= f\vert_{[i,\mu]}
f\vert_{[\mu,j]}$ and $f_{B^\iota}=f\vert_{[\mu,\nu]}=
f\vert_{[\mu,j]} f\vert_{[j,\nu]}$. Therefore $\rho
(\lambda_{B})=[f_B]-[f_{B^\iota}]$ is
  the homology class of the loop $f\vert_{[i,\mu]}
(f\vert_{[j,\nu]})^{-1}$. Both paths forming the letter loop
project to contractible loops in $\Sigma_1/I$. This implies
$(\ast)$.  If $i<\mu<\nu<j $, then   $f_B=f\vert_{[i,j]}$ is the
product of the paths $ f\vert_{[i,\mu]}$,
$f_{B^\iota}=f\vert_{[\mu,\nu]}$, and $ f\vert_{[\nu,j]}$.
Therefore $\rho (\lambda_{B})=[f_B]-[f_{B^\iota}]$ is
  the homology class of the loop $f\vert_{[i,\mu]}
 f\vert_{[\nu,j]} $. As above, this implies     $(\ast)$.

To prove the remaining equalities $e_w (\lambda_{B}, C)=e_w
(\lambda_{B}, s)=0$, we need more notation. 
For every $r=1,...,k$,  we define two
points  $F_r,G_r\in f(S^1) \cap \partial \Sigma_1$. Let the first
and the last  letters of $v_r$ be numerated by $i=i(r),j=j(r)\in
\widehat  n$ with $i<j$. Then $F_r$ is the input of $d_{w(i)}$ lying on
the right side of    $R_{w(i-1), w(i)}  $ and $G_r$ is the output
of $d_{w(j)}$ lying on the left side of $ R_{w(j), w(j+1)}$.
  (If $i=1$, then $i-1$ should be replaced  
with $n$, and if $j=n$, then $j+1$ should be replaced with~$1$.)
   The sub-path of $f$
leading from $F_r$ to $G_r$ lies in $\Sigma_1$, and  the sub-path
of $f$ leading from $G_r$ to $F_{r+1}$ lies in  $\Sigma_2=\overline{\Sigma-\Sigma_1} \subset \Sigma$. Clearly,  $\Sigma_2$ is a compact 
(possibly, disconnected) surface. We endow $\Sigma_2$ with the orientation induced by the  one in $\Sigma$.  The set $Y= \Sigma_1\cap \Sigma_2= \partial \Sigma_1 \cap \partial \Sigma_2$ consists of $2k$ disjoint closed  intervals  each meeting $f(S^1)$  transversely in one of the points $\{F_1,G_1,..., F_k,G_k\}$.
  The involution $I$ on
$\Sigma_1$ satisfies $I(F_r)=G_r$  for all $r$  and   sends the interval in $Y$ containing $F_r$ to the interval in $Y$ containing $G_r$. Hence $I(Y)=Y$. It is clear that the involution $I \vert_{Y}$ inverts the  orientation  on
$ Y$ induced from  the one on $\Sigma_2$. Let $\Psi$ be the compact oriented
 surface obtained from $\Sigma_2$ by
identifying each point $y\in Y\subset  \partial \Sigma_2$ with $I(y)\in
Y$. The embedding $ \Sigma_2\times 1 \hookrightarrow
\Sigma\times 1$ induces an embedding $ \Psi \hookrightarrow
\partial M$ whose image is disjoint from $  \Sigma=\Sigma \times 0 \subset \partial M$. The projection $\eta: \Sigma\times [0,1]\to M$ maps $\Sigma_2\times 1$
   onto  $ \Psi $ and maps
$ (\Sigma_1-\partial \Sigma_1 )\times 1$ to $ M-\partial M$.

We   define a loop $g$ on $\Psi\subset \partial M$. It  starts  in
$f(0)$ and goes along $f$ in $\Sigma_2$ until hitting $F_1$, then
it switches to $I(F_1)=G_1$ and goes along $f$ in  $\Sigma_2$
until hitting $F_2$, then it switches to $G_2$, etc., until
finally returning to $f(0)$. The loop $g$ is continuous since the
points $F_r, G_r$ are identified  in $\Psi$ for all $r$.  In a
sense, $g$ is obtained by cutting out from $f$ the $k$ sub-paths
lying on $\Sigma_1$ and corresponding to $v_1,...,v_k$. Since
these sub-paths project to contractible loops in $\Sigma_1/I$, the
loops $ g$ and $\eta (f\times 1)$ are homotopic in
$  M$. The loop $f\times 1$ being
homotopic to $f=f\times 0$ in $\Sigma\times [0,1]$, we can conclude
that $ g$ is homotopic to $ f$ in $M$. Therefore the homology class
$[f]-[ g]\in   H_1(\partial M)$ lies in the kernel of the
inclusion homomorphism $H_1(\partial M)\to H_1(M)$. Since
$\rho(\lambda_{B}) = [f_B]\pm [f_{B^\iota}]$ also lies in this
kernel, its intersection  number with $[f]-[ g]$ is equal to 0.
 On the other hand, this   number is equal to $b(\rho(\lambda_{B}),[f]
)$ since the  loops $f_B, f_{B^\iota}$ do not meet $ g$ (they lie
in disjoint subsurfaces of $\partial M$).  Thus
$b(\rho(\lambda_{B}),[f] )=0$. By Lemma \ref{ltye4rtutu12}, $e_w
(\lambda_{B}, s)=0$.

If $C\in \C$, then the loop $f_C$ on $\Sigma$     intersects
$\partial \Sigma_2 $ in the points $\{F_r,G_r\}$, where $r$ runs
over all indices $1,...,k$ such that the word $v_r$  lies between
the two entries of $C$ in $w$. Cutting out from  $f_C$ the
sub-paths in $\Sigma_1$   corresponding to all such $v_r$, we
obtain a loop $g_C$ in $\Psi\subset \partial M$  homotopic to
$\eta (f_C\times 1)$ in $M$. Since the loop $f_C\times
1$ is homotopic to $f_C=f_C\times 0$, we conclude that $g_C$ is
homotopic to $\eta f_C$  in $M$.  The  same  argument as in   the
previous paragraph  shows that $e_w (\lambda_{B},
C)=b(\rho(\lambda_{B}),[f_C])=0$.
\qed

\subsection{Remarks}\label{remsssrss} 1.   The geometric interpretation of nanowords over $\alpha_0$ may
be extended to  nanowords over an arbitrary alphabet $\alpha$. One
possibility is to consider   equivariant mappings $\alpha\to
\alpha_0$ and the corresponding push-forwards of nanowords. In
this way any nanoword over $\alpha$ determines a family of pointed
spinal loops on surfaces numerated by the equivariant mappings
$\alpha\to \alpha_0$. Another geometric interpretation of
nanowords may be obtained by considering loops     with additional
data in the self-intersections. This data may be an
over/under-crossing information or a label. For more on this, see
\cite{tu4}.

2. Consider a nanoword $(\A, w)$ over $\alpha_0$ and  the tautological filling $\lambda=\{\lambda_i\}_i$ of
the associated $\alpha_0$-pairing $ e_w $. Let $\varphi: \pi(\alpha_0, \tau_0) \to   \ZZ$
 be the   identification isomorphism.
 By Lemma \ref{ltye4rtutu12}, the matrix $(\varphi e_w
 (\lambda_i, \lambda_j))_{i,j}$, considered up to multiplication of rows and
 columns by $2$,
  is  the matrix of homological intersections of the loops
 $\{f_A\}_{A\in \A}, f$ on the surface $\Sigma_w$ associated with $w$. Since the homological
 classes of these loops generate $H_1(\Sigma_w)$, the rank of
 this
 matrix  is equal to  $2 \ar (\Sigma_w)$, where $\ar (\Sigma_w)$ is the genus
 of $\Sigma_w$. Hence $\sigma_\varphi(\lambda)=\ar (\Sigma_w)$. This equality prompted 
  the term \lq\lq genus" for
$\sigma_\varphi$.
 We can conclude that $\sigma_\varphi(w)\leq \ar (\Sigma_w)$.

  \section{Surfaces in 3-manifolds}\label{fdprellleisionf40}

We discuss  properties of surfaces in   3-manifolds   needed in
the next section to prove Lemma \ref{tphi}.

\subsection{Simple surfaces}\label{2prelimitoprss} Let $F$ be a compact subspace
of a  3-manifold $N$.  A point $ a\in F $  is a {\it branch point}
 if   it lies inside a closed 3-ball $D^3\subset  N $
such that  $F\cap D^3$ is the cone over a figure eight loop in
$S^2=\partial D^3$ with cone point $a \in \Int D^3$. Here a figure
eight loop in $S^2$ is a loop with one transversal
self-intersection. The set of branch points of $F$ is denoted
$Br(F)$. Clearly, $Br(F)\subset \Int N=N-\partial N$.    We call
$F$ a  {\it simple
  surface in $N$} if
 any point of $ F-Br(F)$ has  an open neighborhood $V\subset
N$ such that the pair $(V, V\cap F)$ is homeomorphic to either
$(\RR^3, \RR^2 \times 0  )$, or  $(\RR^3, \RR^2 \times 0 \cup
0\times \RR^2)$, or $(\RR^2, \RR \times 0) \times \RR_+ $, or
$(\RR^2,  \RR \times 0
  \cup 0\times
\RR )\times \RR_+$,  where $ \RR_+=\{r\in \RR\,\vert\, r\geq 0\}$.
Points  of $  F -Br(F) $ having   neighborhoods of the first or
third type  are
 {\it flat}. Non-flat points  of $  F -Br(F)$ are called
   {\it double point} of $F$.   They form a   1-manifold    $d(F)$  with boundary
  $d(F)\cap \partial N$. The closure
  $\overline {d(F)} =d(F)\cup Br(F)$ is a
   compact 1-manifold with boundary $\partial d(F)\cup
   Br(f)$.

A simple surface $F $ in $N$ can be parametrized by an
  abstract surface $\widetilde  F$ obtained by blowing up  the double
  points of $F$. More precisely,   cutting out $F$
 along $d(F)$ we obtain a compact surface $F_{cut}$ and a projection $p: F_{cut} \to F$.
  For   $a\in d(F)$, the set $p^{-1}(a) \subset \partial F_{cut}$
 consists of 4 points adjacent to 4 branches of $F-d(F)$ near $a$. Moving
 around $d(F)$ in a neighborhood of $a$ in $N$ we can cyclically numerate these
 branches - and the corresponding points  of $p^{-1}(a)$ -
 by  the numbers  1,2,3,4. The permutation $1\leftrightarrow 3$, $2\leftrightarrow 4$
  defines an involution  on  $p^{-1}(a)$. This
 gives  a free involution on $p^{-1}(d(F))\subset \partial F_{cut}$
 commuting with   $p$. Identifying
 every point of $  p^{-1}(d(F))$  with its image under this involution, we transform
 $F_{cut}$ into a compact   surface, $\widetilde  F$. The mapping $p$ induces a mapping
 $\widetilde  F\to F$ denoted $\omega$. This is a parametrization of
 $F$ in the sense that the pre-image of each
 double  point of $F$ under $\omega$ consists of 2 points
 and the restriction of $\omega$
  to the complement of this pre-image is a homeomorphism onto $F-d(F)$.

 Suppose from now on that $N$ and  $ \widetilde  F$ are oriented and
 provide $\partial N, \partial \widetilde  F$ with induced orientations. Suppose also
 that
$\partial \widetilde  F$  is homeomorphic to a circle.
   The mapping $h=\omega\vert_{\partial \widetilde  F}:\partial
 \widetilde  F \to F \subset N$
 is a   (generic) loop on $\partial N$ and  $F\cap \partial N= h(\partial \widetilde  F)$.
   Let $\Join \, \subset \partial N$ be
 the set of double points of
$h$. We define an involution $\nu$ on $\Join$ as follows.
   Each point $x\in \, \Join $ is an endpoint of a component of $d(F)$. If this
component is compact, then it has another endpoint, $y  \in\, 
\Join$, and
  $\nu (x)=y $. Otherwise,  $\nu (x)=x$.

  Fix   a base point   $\ast\in \partial \widetilde  F$ such that $h(\ast)\notin \, \Join$.
  For any $x\in \,\Join$, consider the  path $h_x:[0,1]\to  \partial N$
beginning at $x$,  following along $h(\partial \widetilde  F)$ until
the first return to $x$ and not passing through $h(\ast)$. Set
$\sign (x) =+$ if  the pair of tangent vectors $(h'_x(0),
h'_x(1))$ is
  positively oriented in
    the tangent space of $x$ in $\partial N$  and   $\sign (x) =-$ in the opposite case.
  The path $h_x$  determines a loop $S^1\to \partial N$ whose homology
  class in $H_1(\partial N)$ is denoted
 $[h_x] $. For a subset $X$ of $\Join$, set $$[X]=\sum_{x\in X} \sign (x)\, [h_x] \in
H_1(\partial N).$$

    \begin{lemma}\label{lderfnprelitooodddd} Let  $\inc:H_1(\partial N)
\to H_1(N)$ be the inclusion homomorphism. For any    orbit $X $
of the involution $\nu$ on $\Join$, we have  $ \inc ([X] ) \in
\omega_\ast(H_1(\widetilde  F)) \subset H_1(N) $.
\end{lemma}

\begin{proof} We first  compute $\inc ([X])$ as follows.
 For any $x\in \, \Join$, consider the  path $\omega_x:[0,1]\to  \partial N$
beginning at $x$,  following along $\omega (\partial \widetilde 
F)=h(\partial \widetilde  F)$ until the first return to $x$ and such
that the pair of tangent vectors $(\omega'_x(0), \omega'_x(1))$ is
  positive   in
    the tangent space of $x$ in $\partial N$. In contrast to
    $h_x$, the path $\omega_x$ does not depend on the choice of
    the base point $\ast$.
  The path $\omega_x$  determines a loop $S^1\to \partial N$ whose homology
  class   $[\omega_x] \in H_1(\partial N)$ is equal to $[h_x]$ if $\sign (x) =+$ and to $[h]-[h_x]$ if $\sign (x) =-$.
   Therefore $[X]=\sum_{x\in X}   [\omega_x]$ modulo $[h] $. Since $h:\partial \widetilde  F \to \partial N$ extends to a mapping of $\widetilde  F $ to $N$, we have
   $\inc ([h])=0$. Therefore $\inc ([X])=\sum_{x\in X}
   \inc ([\omega_x])$.

Set $T=\omega^{-1}(d(F) \cup Br(F))\subset \widetilde F$.   A  local inspection shows
that $T$  is an embedded 1-manifold in $ \widetilde  F$ with
$\partial T= T \cap
\partial \widetilde  F= \omega^{-1} (\Join)$.
For any  point $a\in \omega^{-1}(d(F))\subset T $ there is exactly
one other point $b \in \omega^{-1}(d(F))$ such that
$\omega(a)=\omega(b)$.  The correspondence $a \leftrightarrow b$
extends by continuity to an involution $\Delta$  on $  T$ with
fixed-point set  $
    \omega^{-1}(  Br(F))$.    The mapping $\omega$
defines a homeomorphism   $\partial T/\Delta=\Join $. We identify
these two sets along  this homeomorphism.

  For
$a\in \partial T$,  let $I_a  $ be the
 component of $ T$ with endpoint  $a$ and  let  $\mu(a)\in
\partial T$     be its  other endpoint.  The formula $a\mapsto \mu_a$
 defines a fixed-point-free involution $\mu$ on $\partial T$.
  We claim that $\mu$ commutes with $\Delta\vert_{\partial T}$.
Indeed, if  $ \Delta(I_a)=I_a$, then   $\Delta$ exchanges the
endpoints of $I_a$ so that   $ \Delta=\mu $ on $\partial I_a$. (In
this case $\Delta$ must have a unique fixed point inside $I_a$.)
If $\Delta(I_a)\neq I_a$,  then  $\Delta(I_a)$ has the endpoints
$\Delta(a), \Delta(\mu(a) )$ so that $\mu(\Delta(a))=\Delta
(\mu(a))$. Since $\Delta \vert_{\partial T}$ and $\mu$ commute,
$\mu$ induces an involution on $\partial T/\Delta$. Under the
identification $\partial T/\Delta=\Join $, the latter  involution  coincides
with $\nu$.

We   now verify that $ \inc ([X] ) \in \omega_\ast(H_1(\widetilde  F))
 $ for any  orbit $X$ of $\nu$.  Pick    $x\in X$. The
  path $\omega_x$
  in $\partial N$ defined above is obtained (up to reparametrization)
  by restricting $\omega$ to
  an  arc
$\gamma_x\subset \partial \widetilde  F$ leading from a point $a$ to  a point $b$, 
where $\{a,b\} = \omega^{-1}(x)\subset   \partial T \subset
\partial \widetilde  F$.

Assume first that $ \nu(x)=x$. Then $\mu(a)\in \{a,b\}$ and since
$\mu(a)\neq a$, we have $\mu (a)=b$.  By the definition of
$\Delta$, we have $\Delta(a)=b$ and $ \Delta(b)=a$. Since $\Delta:
T \to T $ preserves the set $\partial I_a=\{a, b\}$, we have
$\Delta(I_a)=I_a$. Observe that the product of the path $\gamma_x
$   with the interval $ I_a \subset \widetilde  F$ oriented  from $ b$
to $a$ is a loop, $\rho$, in $\widetilde  F$. The loop  $\omega(\rho) $
in $N$ is a product of $\omega(\gamma_x)=\omega_x$  with the loop
$\omega \vert_{I_a}$. The latter  loop  is contractible in $N$
because it has the form $ \delta \delta^{-1}$ where $\delta $ is
the path in $  N$ obtained by restricting $ \omega $ to the arc in
$I_a$ leading from $b$ to the
  fixed point of $\Delta$  on $ I_a$.  Hence $\inc ([X])= \inc
([\omega_x])=[\omega(\rho)]\in \omega_*(H_1(\widetilde  F))$.

Suppose    that  $ \nu(x)\neq x$.  Inspecting the orientations of
the sheets of $F$ meeting along $\omega(I_a)$, we observe that the
path $\gamma_{\nu(x)} $ begins  at    $\mu(b)$ and terminates at $
\mu(a)$  (this was first pointed out by  Carter \cite{ca2}).
  Consider   the  loop $\rho= \gamma_x I_{b}
\gamma_{\nu(x)}   (I_{ a})^{-1} $ in $\widetilde  F$ beginning and
ending at $a$. Here the intervals   $ I_{b} ,   I_{a}$ are
oriented from $b$ to $\mu(b)$ and from $a$ to $\mu(a)$,
respectively. Then $\omega(\rho)$ is the product of the loop
$\omega_x$ beginning and ending at $x$, the path $\omega   (
I_{b})$ leading from $x$ to $ \nu(x)$, the   loop
$\omega_{\nu(x)}$ beginning and ending at $ \nu(x)$, and  the path
$(\omega  ( I_{ a}))^{-1}$ leading from $ \nu(x)$ to $x$. The
paths $\omega  ( I_{b})$,  $(\omega  ( I_{ a}))^{-1}$ are mutually
inverse  since $ I_{b} = \Delta( I_{a}) $ and $\omega
  \Delta =\omega  $. Hence
$$\inc ([X])= \inc ([\omega_x]  + [\omega_{\nu (x)}]) =[\omega(\rho)]\in
\omega_*(H_1(\widetilde  F)).$$
\end{proof}

 \begin{lemma}\label{78320ddd}
  Let
$b$ be the intersection form $H_1(\partial N)\times H_1(\partial
N) \to \ZZ$.   Let $X_1,...,X_t$ be the orbits of the involution
$\nu:\, \Join \to \Join$. Set $c_i=[X_i] \in H_1(\partial N)$ for
$i=1,...,t$ and $c_0= [h(\partial \widetilde  F)]=[\omega(\partial
\widetilde  F)]\in H_1(\partial N)$. Then  the rank of the $(t+1)\times
(t+1)$-matrix
 $(b(c_i ,
c_{j} ))_{i,j=0,1,...,t}$  is smaller than or equal to $4\ar  $
where $\ar=\ar(\widetilde  F)$ is the genus of $\widetilde  F$.
\end{lemma}
\begin{proof}
The  group   $H=H_1(\widetilde  F)$  is   isomorphic to $\ZZ^{2\ar}$.
Set $L= \inc^{-1} (\omega_\ast(H))\subset  H_1(\partial N)$. Since
the intersection form $b$ annihilates the  kernel of $\inc$,
\begin{equation}\label{zx}\rk\, (b\vert_L:L\times L \to \ZZ) \leq 2 \rk\,
\omega_\ast(H)
 \leq 2\rk  H=4 \ar \end{equation}
 where $\rk$ is the rank of a bilinear form or of an abelian
 group.
By Lemma \ref{lderfnprelitooodddd},     $c_i  \in L $ for
$i=1,...,t$.  Also $c_0 \in  L $ since $\inc (c_0)=0$.   The claim
of the lemma now follows  from Formula \ref{zx}.
\end{proof}

\subsection{Remark} \label{beremkrss}    {\it Generic surfaces} in
3-manifolds are defined as   simple surfaces but additionally
allowing triple points where the  surface looks like the union of  three 
coordinate planes in $\RR^3$. Although we shall not need it, note
that Lemmas \ref{lderfnprelitooodddd} and \ref{78320ddd} extend to
generic surfaces, cf.\ \cite{tu1}.

 \section{Proof of Lemma \ref{tphi}}\label{fdvmvmvmvmmsospmlensionf40}

 \subsection{Notation} \label{begins2prelimitoprss}
Consider  a bridge  in a nanoword  $(\A,w)$ over $\alpha$  formed by  a
factor
  $\nabla=(\B, (v_1\,\vert \, \cdots \, \vert \, v_k))$
    and an involution $\kappa:\widehat  k \to \widehat  k$ on
    $\widehat  k=   \{1,2,..., k\}$.
    Thus, $\B\subset \A$  and $w=x_1 v_1x_2v_2   \cdots
x_{k} v_k x_{k+1}$ where $x_1,x_2,..., x_{k+1}$ are words in the
$\alpha$-alphabet $\C=\A-\B$. The associated bridge move $m$
transforms $w$ in    the nanoword $(\C, x=x_1 x_2 \cdots
x_{k+1})$. Set $\iota=\iota_{\nabla, \kappa}:\B \to \B$ and
$\varepsilon=\varepsilon_{\nabla, \kappa}:\B \to \{0,1\}$.

Replacing each letter $C\in \C$ by its copy $C'$ we obtain  a
nanoword
 $ (\C'=\{C'\}_{C\in \C}, x')$  isomorphic to  $(\C,x)$. The
 nanoword $wx^-$ is isomorphic to $w(x')^-$.
 Consider   the $\alpha$-pairing of the latter  nanoword $$p_{w(x')^-}=( S=\A \cup \C' \cup
\{s\},s,  e_{w(x')^-} : S \times  S \to \pi(\alpha, \tau)).$$
  For $C\in \C$, set
 $\lambda_C=C+ C'\in  \ZZ  S$. For    $B\in \B$, consider
 the vector
 $\lambda_B\in \ZZ  S $ equal to $B$ if $B=B^\iota$ and equal to
 $B+ (-1)^{\varepsilon (B)} B^\iota$  if $B\neq B^\iota$.
 (Note that $\lambda_{B^\iota}= (-1)^{\varepsilon  (B)} \lambda_B$.)
Pick a set $\B_+\subset \B$ meeting every orbit of
 $\iota: \B \to \B $ in one  element.
  The set of vectors
\begin{equation}\label{BBB1}\lambda=\{\lambda_A \}_{ A\in \B_+\cup \C}
 \cup \{
s+s'\}\end{equation} is a filling  of $ p_{w(x')^-}$. It plays a
crucial role in the next lemma.

\begin{lemma}\label{l:e56bbbpp009c12}
If   $\alpha=\alpha_0=\{+,-\}$, $\tau=\tau_0$,  and
$\varphi_0:\pi(\alpha_0, \tau_0)\to \ZZ$ is the canonical
isomorphism, then $\ar(m)\geq \sigma_{\varphi_0}(\lambda)/2$.
\end{lemma}
                     \begin{proof} Applying  the   constructions
                  of   Theorem \ref{edrdcppmolbrad} to $w$, we obtain a pointed spinal
 loop $f:S^1\to \Sigma$ where $S^1=\RR\cup \{\infty\}$ and the origin of $f$ is the point $f(0)$. The self-crossings of $f$
 are labelled by elements of $\A$ bijectively. The first part
                    of the proof of Lemma \ref{l:edlocfc12} (till
                    the formula $I^2=\id$)
  applies  word for word, though here $\iota=\iota_{\nabla, \kappa}:\B \to \B$ and
$\varepsilon=\varepsilon_{\nabla, \kappa}:\B \to \{0,1\}$.  This
gives    a  surface $\Sigma_1\subset \Sigma$  and
 an orientation-reversing involution $I:\Sigma_1\to \Sigma_1$.
This involution  maps the sub-path of $f$ corresponding to $v_r$
onto itself for all $r=1,...,k$ such that $\kappa (r)=r$. If
$\kappa(r)\neq r$, then $I$ maps the sub-path of $f$ corresponding
to $v_r$ onto the sub-path of $f$ corresponding to $v_{\kappa
(r)}$ with reversed direction. We can define     surfaces $\Sigma_2=\overline {\Sigma-\Sigma_1}$ and $\Psi$ as well as  a 3-manifold $M$ as
in the proof of Lemma \ref{l:edlocfc12}. However, $\Psi$ and $M$ are  inadequate
for our aims. The problem is that the pieces of   $f(S^1)$ lying
on $\Sigma_2 $ may not form
a single loop in $\Psi$.   For example, for
$w=x_1v_1x_2v_2x_3$ and $\kappa=(12)$, this procedure gives  two
loops: one is glued from  the paths arising from $x_1, x_3$ (the
involution $I$ maps the head of the first path to the tail of the
second path) and another  loop is the image of the path arising
from $x_2$ (the involution $I$ permutes its endpoints). To
circumvent this problem, we modify our constructions as follows.

 Pick a small  positive number $\delta< 1/10$.
 Let $\R$ denote the set of  all $r=1,...,k$ such that $r<\kappa (r)$.  For $r\in \R$,
consider the ribbon $ R_{A,B}\subset \Sigma$ where $B\in \B$ is the
first letter of $v_r$ and $A\in \A$ is the preceding  letter in $w$ (we may have $A\in \B$ if there are no letters beween $v_{r-1} $ and $v_r$; if
$r=1$, then $A$ is the last letter of $w$). Let $D_r\subset R_{A,B}$ be the
rectangle    defined in the coordinates $(p,q)$ by  $D_r=[3/4-\delta, 3/4+\delta]\times [-1/10, 1/10]$. This rectangle meets
  $f(S^1)$ along the arc $[3/4-\delta, 3/4+\delta]\times 0$.
Consider also the ribbon   $R_{B',A'}\subset \Sigma $ where $B'\in \B$
  is the last letter of $v_{r}$ and $A'\in \A$ is the
  (cyclically) next letter of $w$. Let $D'_r\subset R_{B',A'}$ be the
rectangle   defined in the coordinates $(p,q)$ by  $D'_r=[1/4-\delta, 1/4+\delta]\times [-1/10, 1/10]$. This rectangle meets
  $f(S^1)$ along the arc $  [1/4-\delta, 1/4+\delta]\times 0$.
  Set $D=\cup_{r\in \R} (D_r\cup D'_r)\subset \Sigma-\Sigma_1$.
  We choose $\delta$     small  enough so that the origin
  $f(0)$ of $f$ does not belong to $D$ and moreover, the arc $f([-\delta, 0])$ is disjoint from $\Sigma_1$ and from $D$. (To ensure these properties one may need to deform 
  the coordinate  $p$ on the ribbon containing $f(0)$.) 

   We extend the involution $I$ on $\Sigma_1$ to an orientation reversing involution
$I'$ on the (disconnected) surface 
 $\Sigma_3=  \Sigma_1 \cup D$  which sends  a point with coordinates
 $(p,q)$ on $D_r$ to the point with coordinates $(1-p,q)$ on $D'_r$ for all $r\in \R$.
   Clearly,  $\Sigma_4=\overline {\Sigma -\Sigma_3} $  is a compact oriented
  surface. The set $Y'= \Sigma_3\cap \Sigma_4= \partial \Sigma_3
 \cap \partial \Sigma_4$ consists of $2k+4 \card (\R)$ disjoint closed  intervals,  each meeting $f(S^1)$  transversely in one   point.
It is easy to see that $I'(Y')=Y'$ and  the restriction of $I'$ to $Y'$  inverts the  orientation  on
$ Y'$ induced from  the one on~$\Sigma_4$.   Let
$\Psi'$ be the  compact  oriented
 surface obtained from $\Sigma_4$ by
identifying each point $y\in Y'$ with $I' (y) \in
Y'$. One may check that $\Psi'$ is obtained from $\Psi$ by adding $\card (\R)$ one-handles.

The pieces of $f(S^1)$ lying on $\Sigma_4$ glue together into
a single loop   $g':S^1\to \Psi'$.  The point $f(0)\in \Int \Sigma_4$
serves as the origin of $g'$.  Note that $f$ and $g'$ have the same
germ in their common origin $f(0)=g'(0)$.   The
  self-crossings of $g'$ are precisely the self-crossings of $f$ labelled by
the elements of $\C=\A-\B$. We  prefer to   label the self-crossings of $g'$
with elements of $\C'$ rather than the corresponding elements of~$\C$. The underlying nanoword of $g'$ is then the copy
 $ (\C', x')$  of   $(\C,x)$.

Let $N$ be the oriented 3-manifold obtained from  $\Sigma \times [0,1]$ by the identification $a\times 1 =
I'(a)\times 1$ for all $a\in \Sigma_3$.
    Denote  the projection $\Sigma \times
[0,1]\to N$ by $\eta'$. The embedding $ \Sigma_4\times 1
\hookrightarrow \Sigma\times 1$ composed with $\eta'$ yields an
inclusion $ \Psi' \hookrightarrow
\partial N$ whose image is disjoint from $\Sigma= \Sigma\times
0\subset \partial N$.     It is easy to check that $F=\eta'(f(S^1) \times [0,1])$ is a simple
surface in $N$ in the sense of Sect.\ \ref{2prelimitoprss}. Its
set of branch points is  $\{\eta'(V_B\times 1)\}_B$, where $B  $ runs over the letters in $ \B$ such that  $B=B^\iota$, and $V_B$ denotes
the self-crossing of $f$ labelled by $B$. The double points of $F$ are the points of type $\eta'(V_B, t)$, where $B\in \B$ and $t\in [0,1]$. The set $F\cap
\partial N$ consists of two   loops $f(S^1)\subset \Sigma$ and $ g' (S^1)\subset \Psi'$.

We now modify $N$ and $F$ to obtain a simple surface with
connected boundary. Consider a   2-disk $D_0\subset \Int \Sigma_4$
meeting $f(S^1)$ along the  arc $f([-\delta, 0])$. Set $$N_0=N- \eta'(\Int D_0 \times
[0,1]).$$ Then $N_0$ is a compact oriented 3-manifold
with $\partial N_0 \supset \Sigma \# ( -\Psi')$,  where $\#$ is the
connected sum of surfaces, and  the sign $-$ reflects the fact that
the orientation of $\Psi'$ induced from $N$ is opposite to the one
induced from $\Sigma$. The set $F_0=F\cap N_0 $ is obtained from
$F$ by removing an embedded band joining two components of $ F\cap
\partial N $ in the complement of branch points and double points. Clearly, $F_0$ is a simple surface in $N_0$ with the same
branch points and double points as $F$. Blowing up the double
points of $F_0$, we obtain  a parametrization $\omega: \widetilde  F_0
\to F_0$ by an abstract surface $\widetilde  F_0$. The construction of
$F, F_0$ implies that $ \widetilde  F_0$ is a compact connected
orientable surface   with boundary   homeomorphic
to $S^1$. The genus of  $\widetilde  F_0$ is easily seen to be equal to  the number of arches 
 $\ar(m)=\card (\R)$ of $m$.

 The loop $h=\omega\vert_{\partial
\widetilde  F_0}:\partial \widetilde  F_0\to
\partial N_0$ starts at   $f(0)$ (which serves as the origin) and goes along $f$ in $\Sigma $
till   $f(-\delta)$, then   along $\eta'(f(-\delta)\times [0,1])$
to $\eta'(f(-\delta)\times 1) $, then   along $(g')^{-1}$ in
$\Psi'$ till  $\eta'(f(0)\times 1)$, and finally     down  to
$f(0)$ along $\eta'(f(0)\times [0,1])$. The self-crossings of $h$
are those of $f$ and those of $g'$. They are bijectively labelled
by elements of $\A\cup \C' $. The self-crossing of $h$ labelled by
a letter $A\in \A\cup \C'$ is denoted $V_A$. The underlying
nanoword of $h$ is $(\A\cup \C', w(x')^{-})$.

We apply to $F_0$ and $h$ the definitions of Sect.\
\ref{2prelimitoprss}.  The
 involution $\nu$ on the set of self-crossings
 of $h$ permutes $V_C, V_{C'}$ for
all $C\in \C$ and sends $V_B $ to $V_{B^\iota}$ for $B\in \B$.  By
Lemma \ref{78320ddd},  the matrix, $K$, of the intersection form
on $H_1(\partial N_0)$ computed on the vectors $[h], \{[X]\}_X\in H_1(\partial N_0)$, where $X$ runs
over the orbits of~$\nu$,  satisfies  $$ \rk K \leq 4
\ar (\widetilde  F)=4 \ar (m).$$ We now compute the vectors $[X]$.
We shall write $[h_A]$ for the homology class
$  {[h_{V_A}]} \in H_1(\partial N_0)$, where $A\in \A\cup \C'$. Note that 
$\sign (V_A) =\vert A\vert$
for   $A\in \A$ and   $\sign (
V_{C'}) =\vert C\vert$
for   $C\in \C$;
 the latter equality follows from the fact  that $h$ goes
along $(g')^{-1}$ on $\Psi'$ and that the orientation on $\Psi'$
induced from $N$ is opposite to the one induced from~$\Sigma$. 
For the orbit $X=\{V_C, V_{C'}\} $  of $\nu$ with  $C\in \C$,  
\begin{equation}\label{ymin}   {[X]}= \pm ([h_C]+
[h_{C'}]).\end{equation}  For the orbit $X=\{V_B, V_{B^\iota}\} $, where
  $B\in \B$ and $B^\iota \neq B$,  
\begin{equation}\label{ymin+}  {[X]}= \vert B\vert [h_B]+ \vert B^\iota \vert   [h_{B^\iota}]=
\pm  ( [h_B]+ (-1)^{\varepsilon (B)} [h_{B^\iota}]).\end{equation}  For the orbit $X=\{V_B \} $, where $B=B^\iota \in \B$,  
\begin{equation}\label{ymin++}  {[X]}=\pm [h_B].\end{equation}

Consider now the filling $\lambda=\{\lambda_i\}_i$ of $ p_{w(x')^-}$
given by (\ref{BBB1}).  Here    $i$ runs over the subset $\B_+\cup \C$ of $\A $ plus one additional index numerating $s+s'$. We can apply
  Lemma
\ref{ltye4rtutu12}  to the loop $h$ representing the nanoword $w(x')^{-1}$. This lemma computes the   matrix
$({\varphi_0} e_{w(x')^-} (\lambda_i, \lambda_j))_{i,j}$ in terms of the intersection numbers of  (formal linear combinations of)   loops on $\partial N_0$.  The loops in question are $h$ and the  formal linear combinations of   loops   appearing  on the right-hand sides of Formulas \ref{ymin} -- \ref{ymin++}. Therefore  the matrix $({\varphi_0} e_{w(x')^-} (\lambda_i, \lambda_j))_{i,j}$ a
sub-matrix of the matrix $K$,  at least  up to multiplication of
rows and columns by $\pm 2$ and $\pm 1$.  Therefore the half-rank
$\sigma_{\varphi_0} (\lambda)$ of the former matrix   can not
exceed $(\rk  K)/2 \leq 2 \ar (m)$. Hence,
 $ \sigma_{\varphi_0} (\lambda)/2\leq \ar(m)$.
\end{proof}

\subsection{Proof of  Lemma \ref{tphi}}\label{2predddlimitoprss}  It suffices to verify that   $\ar(m)\geq  \sigma_\varphi
(\lambda)/2$, where $\lambda $ is the filling (\ref{BBB1}) of $
p_{w(x')^-}$. By assumption, $\varphi(\alpha)\subset  \alpha_0$.
 Pushing   forward a
                     nanoword
                     $v$ over $\alpha$ along $\varphi\vert_\alpha:\alpha\to \alpha_0$, we obtain a nanoword   over
                     $\alpha_0$ denoted by $v_0$.  Every filling
                     $\mu$ of the $\alpha$-pairing $p_v$ induces a
                     filling $\mu_0$ of the $\alpha_0$-pairing $p_{v_0}$ (actually, $\mu_0=\mu$ as
                     sets of vectors). By Section \ref{fi2ampsortnib99}, 
                     $\sigma_\varphi (\mu)=\sigma_{\varphi_0}
                     (\mu_0)$, where $\varphi_0:\pi(\alpha_0,
                     \tau_0)\to \ZZ$ is the
                     canonical isomorphism. We   apply this observation  to $v=w(x')^-$ and the
                     filling $\mu=\lambda $ of $p_{v}$. Here $v_0=w_0(x'_0)^-$ and  the induced filling $\lambda_0$ of
$p_{v_0}$ is  also given by Formula  \ref{BBB1}. The
bridge move $m:w\to x$
                     induces a bridge move $m_0:w_0\to x_0$ with
                     the same number of arches.   By  Lemma \ref{l:e56bbbpp009c12},
                     $$\ar (m)= \ar(m_0) \geq \sigma_{\varphi_0} (\lambda_0)/2 = \sigma_\varphi (\lambda)/2.$$

 \section{Further directions and open problems}

 1. Give a combinatorial proof of  Formula
  \ref{plo}  and      Lemma \ref{tphi}.   An incomplete  combinatorial approach to   Formula
  \ref{plo} is discussed in  the first version of this paper available in  arXiv:math/0511513.  A combinatorial proof  of Lemma \ref{tphi}  might enable one to
extend  Theorem \ref{tcdre53} to other
$\varphi$.

2. Compute the image of the homomorphism
 $p:{\mathcal{N}}_c \to {\mathcal{P}_{sk}}$ from Sect.\ \ref{fihomomlb99}. 
 
 3. Find further cobordism  invariants of nanowords.  

4. Is it true that for $\alpha$ consisting of only one element,
$\N_c=1$ ? At the moment of writing, nothing contradicts the stronger conjecture
that any two nanowords over a 1-letter alphabet are homotopic.

5. Classify words of small length, say  $\leq 10$, up to
cobordism.

6. A metamorphosis of nanowords over  $(\alpha_0,
 \tau_0)$ gives rise to a generic surface in a 3-manifold  
$N$ interpolating between two disjoint loops in  
$\partial N$. (Besides the constructions
above, one should observe that the third homotopy move naturally
gives rise to a generic surface with one triple point in a
3-dimensional cylinder.)  This defines a functor from the category of nanowords over  $(\alpha_0,
 \tau_0)$ and their metamorphoses to the category of spinal loops on surfaces and interpolating surfaces in 3-manifolds. In what sense this   is an equivalence of categories?   

7. One can model homotopy (resp.\ cobordism) of  surfaces in
3-manifolds to define    homotopy (resp.\   cobordism)  for
metamorphoses of nanowords.  Are  their interesting   invariants of metamorphoses preserved under these relations?

                     \end{document}